\documentclass[12pt,reqno]{amsart}

\usepackage{amssymb}
\usepackage{bbm}
\usepackage{mathrsfs}
\usepackage[all]{xy}
\usepackage{float}
\usepackage[abs]{overpic}
\usepackage{stackrel}

\newtheorem{theo}{Theorem}[section]
\newtheorem{lem}[theo]{Lemma}
\newtheorem{prop}[theo]{Proposition}
\newtheorem{coro}[theo]{Corollary}

\theoremstyle{definition}
\newtheorem{defi}[theo]{Definition}
\newtheorem{term}[theo]{Terminology}
\newtheorem{algo}[theo]{Algorithm}
\newtheorem{mconstr}[theo]{Main construction}
\newtheorem{constr}[theo]{Construction}
\newtheorem{strat}[theo]{Strategy}
\newtheorem{ttt}[theo]{}

\newtheorem*{nota_allg}{Notation}
\newtheorem*{comp_allg}{Computations}

\newtheorem*{proo_periods}{Proof of Theorem~\ref{periods_ints}}
\newtheorem*{proo_indep}{Proof of Theorem~\ref{indep_b}}

\newtheorem*{ass}{Assumptions}

\theoremstyle{remark}
\newtheorem{nota}[theo]{Notation}
\newtheorem{ex}[theo]{Example}
\newtheorem{exs}[theo]{Examples}
\newtheorem{rem}[theo]{Remark}
\newtheorem{rems}[theo]{Remarks}
\newtheorem{cau}[theo]{Caution}

\newcommand{\br}{ }
\newcommand{\brr}{, }

\renewcommand{\atop}[2]{\genfrac{}{}{0pt}{}{#1}{#2}}

\newcommand{\Gal}{\mathop{\text{\rm Gal}}\nolimits}
\newcommand{\cha}{\mathop{\text{\rm char}}\nolimits}
\newcommand{\supp}{\mathop{\text{\rm supp}}\nolimits}
\newcommand{\spann}{\mathop{\text{\rm span}}\nolimits}
\renewcommand{\div}{\mathop{\text{\rm div}}\nolimits}
\newcommand{\Kum}{\mathop{\text{\rm Kum}}\nolimits}
\newcommand{\Bl}{\mathop{\text{\rm Bl}}\nolimits}
\newcommand{\Spec}{\mathop{\text{\rm Spec}}\nolimits}
\newcommand{\Frob}{\mathop{\text{\rm Frob}}\nolimits}
\newcommand{\Div}{\mathop{\text{\rm Div}}\nolimits}
\newcommand{\Pic}{\mathop{\text{\rm Pic}}\nolimits}
\newcommand{\Aut}{\mathop{\text{\rm Aut}}\nolimits}
\newcommand{\End}{\mathop{\text{\rm End}}\nolimits}
\newcommand{\PGL}{\mathop{\text{\rm PGL}}\nolimits}
\renewcommand{\Re}{\mathop{\text{\rm Re}}\nolimits}
\renewcommand{\Im}{\mathop{\text{\rm Im}}\nolimits}
\newcommand{\Tr}{\mathop{\text{\rm Tr}}\nolimits}
\newcommand{\id}{\mathop{\text{\rm id}}\nolimits}
\newcommand{\im}{\mathop{\text{\rm im}}\nolimits}
\newcommand{\rk}{\mathop{\text{\rm rk}}\nolimits}
\newcommand{\pr}{\mathop{\text{\rm pr}}\nolimits}
\newcommand{\bl}{\mathop{\text{\rm bl}}\nolimits}

\newcommand{\tr}{\text{\rm tr}}
\newcommand{\et}{\text{\rm \'et}}
\newcommand{\Hg}{\text{\rm Hodge}}

\newcommand{\Ab}{{\text{\bf A}}}
\newcommand{\Pb}{{\text{\bf P}}}
\newcommand{\Tb}{{\text{\bf T}}}

\newcommand{\fraki}{{\mathfrak i}}
\newcommand{\frakt}{{\mathfrak t}}

\newcommand{\mi}{\mathrm{i}}

\newcommand{\frakX}{{\mathfrak X}}

\newcommand{\bbC}{{\mathbbm C}}
\newcommand{\bbD}{{\mathbbm D}}
\newcommand{\bbF}{{\mathbbm F}}
\newcommand{\bbN}{{\mathbbm N}}
\newcommand{\bbQ}{{\mathbbm Q}}
\newcommand{\bbR}{{\mathbbm R}}
\newcommand{\bbZ}{{\mathbbm Z}}

\newcommand{\calL}{{\mathscr{L}}}
\newcommand{\calO}{{\mathscr{O}}}
\newcommand{\calX}{{\mathscr{X}}}

\newcounter{ABC}
\newenvironment{ABC}{\begin{list}{\rm \Alph{ABC}) }%
{\usecounter{ABC} \leftmargin=0.0pt \labelsep=0.0pt %
\listparindent=0.0pt \labelwidth=0.0pt \parsep=\smallskipamount%
 \itemsep=0.0pt \topsep=0.0pt \partopsep=\smallskipamount}}{\end{list}}

\newcounter{abc}
\newenvironment{abc}{\begin{list}{\rm \alph{abc}) }%
{\usecounter{abc} \leftmargin=0.0pt \labelsep=0.0pt %
\listparindent=0.0pt \labelwidth=0.0pt \parsep=\smallskipamount%
 \itemsep=0.0pt \topsep=0.0pt \partopsep=\smallskipamount}}{\end{list}}

\newcounter{iii}
\newenvironment{iii}{\begin{list}{\rm \roman{iii}) }%
{\usecounter{iii} \leftmargin=0.0pt \labelsep=0.0pt %
\listparindent=0.0pt \labelwidth=0.0pt \parsep=\smallskipamount%
 \itemsep=0.0pt \topsep=0.0pt \partopsep=\smallskipamount}}{\end{list}}

\newenvironment{Aiii}{\begin{list}{\bf A.\roman{iii}) }%
{\usecounter{iii} \leftmargin=0.0pt \labelsep=0.0pt %
\listparindent=0.0pt \labelwidth=0.0pt \parsep=\smallskipamount%
 \itemsep=0.0pt \topsep=0.0pt \partopsep=\smallskipamount}}{\end{list}}

\newenvironment{Oiii}{\begin{list}{\bf 0.\roman{iii}) }%
{\usecounter{iii} \leftmargin=0.0pt \labelsep=0.0pt %
\listparindent=0.0pt \labelwidth=0.0pt \parsep=\smallskipamount%
 \itemsep=0.0pt \topsep=0.0pt \partopsep=\smallskipamount}}{\end{list}}
\makeatletter
\def\hsmash{\relax % \relax, in case this comes first in \halign
  \ifmmode\def\next{\mathpalette\mathhsm@sh}\else\let\next\makehsm@sh
  \fi\next}
\def\makehsm@sh#1{\setbox\z@\hbox{#1}\finhsm@sh}
\def\mathhsm@sh#1#2{\setbox\z@\hbox{$\m@th#1{#2}$}\finhsm@sh}
\def\finhsm@sh{\wd\z@\z@ \box\z@}
\makeatother

\def\rightend#1#2{{%
 \leavevmode\nobreak\hskip .5em plus 1fil
 \penalty600 \hskip 0pt plus -1filll
 \vadjust{}\nobreak\hskip 0pt plus 1filll%
 #1\parfillskip=#2\relax \par}}

\def\eop{\ifmmode\rule[-22pt]{0pt}{1pt}\ifinner\tag*{$\square$}\else\eqno{\square}\fi\else\rightend{$\square$}{0pt}\fi}

\thanks{}

\title[Real and complex multiplication via period integration]{Real and complex multiplication on
{\boldmath $K3$}
surfaces via period integration}

\begin{document}

\author{Andreas-Stephan Elsenhans}

\address{Institut f\"ur Mathematik\\ Universit\"at W\"urzburg\\ Emil-Fischer-Stra\ss e 30\\ D-97074 W\"urzburg\\ Germany}
\email{stephan.elsenhans@mathematik.uni-wuerzburg.de}
\urladdr{https://www.mathematik.uni-wuerzburg.de/computeralgebra/team/elsenhans-step\discretionary{}{}{}han-\discretionary{}{}{}prof-dr/}

\author[J\"org Jahnel]{J\"org Jahnel}

\address{\mbox{Department Mathematik\\ \!Univ.\ \!Siegen\\ \!Walter-Flex-Str.\ \!3\\ \!D-57068 \!Siegen\\ \!Germany}}
\email{jahnel@mathematik.uni-siegen.de}
\urladdr{http://www.uni-math.gwdg.de/jahnel}

\thanks{}

\date{October~24,~2021}



\begin{abstract}
We report on a new approach, as well as some related experiments, to construct families of
$K3$~surfaces
having real or complex multiplication. The~approach is based on an explicit description of the transcendental part of the cohomology in a topological way, using topological tori. Fundamental ideas include considering the period space of marked
$K3$~surfaces,
determining the periods by numerical integration, as well as tracing the modular curve by a numerical continuation~method.\vspace{-6mm}
\end{abstract}

\maketitle
\thispagestyle{empty}

\section{Introduction}

The endomorphism algebra of an elliptic curve
$X$
over
$\bbC$
is isomorphic either
to~$\bbZ$
or to an order in an imaginary quadratic field. The latter phenomenon is called complex multiplication. The~theory of CM elliptic curves is very rich, cf.~\cite[Chapter~II]{Si} or~\cite[Chap\-ter~3]{Co}, and their construction in an analytic setting is~classical.
The~whole theory generalises to higher dimensions, most obviously to abelian~varieties. There~are, however, natural generalisations of elliptic curves, other than abelian~varieties. One~such kind is provided by the surfaces of
type~$K3$.
Indeed,~elliptic curves may be characterised as being the curves with trivial canonical class, a property shared by
$K3$~surfaces.

Since~$K3$~surfaces
do not carry a natural group structure, the endomorphism algebra needs to be defined~cohomologically.
For~$X$
a
$K3$~surface,
$H := H^2(X, \bbQ)$
is a pure
\mbox{weight-$2$}
Hodge structure of
dimension~$22$,
being the direct sum of the image
of~$\Pic X \!\otimes_\bbZ\! \bbQ$
under the Chern class homomorphism and its orthogonal
complement~$T$,
the {\em transcendental part\/}
of~$H$.
The endomorphism algebra
$\End_\Hg(T)$
in the category of pure
\mbox{weight-$2$}
Hodge
structures is either
$\bbQ$,
or a totally real number field, or a CM~field \cite[Theorem~1.6.a)]{Za}, cf.~Proposition~\ref{Zarh}.i).
If~$E \supsetneqq \bbQ$
is totally real then
$X$
is said to have {\em real multiplication (RM)}.
If~$E$
is CM then one speaks of {\em complex multiplication (CM)}.

For example, the Kummer surface
$X := \Kum(E_1 \times E_2)$
attached to the product of two elliptic curves has complex multiplication if one of the elliptic curves~has. In~this case, CM is actually caused by an endomorphism
of~$X$.
In~\ref{CM_ex1} to~\ref{CM_ex4}, we give examples of
$K3$~surfaces
of lower Picard ranks having CM even due to an~automorphism.

On~the other hand, a Kummer surface does {\em not\/} inherit real multiplication from the underlying abelian surface~\cite[Remark~3.5.ii)]{EJ14}. Nonetheless,~the existence of
$K3$~surfaces
having real multiplication has been established for quite some time by analytic, i.e.\ Hodge-theoretic, means. It is known, for instance, that
$K3$~surfaces
having RM
by~$E$
exist for every real quadratic
field~$E$.
The~same is true for every cubic field that is totally~real. Note though that these methods do not provide explicit equations for the surfaces found to~exist. Cf.~the work \cite{vG} of B.~van~Geemen, in particular~\cite[Proposition~3.3]{vG}.

\subsubsection*{Dimensions of the RM and CM loci}
Theorem~\ref{RMCM_dim} is similar in spirit to van~Gee\-men's result. It covers, however, the CM case, too, and provides exact formulae for the dimensions of the RM~and CM loci in period~space.
These~loci are, strictly speaking, not even closed~sets. The~reason is that there is a semicontinuity phenomenon. If the general fibre of a family has RM or CM by a
field~$K$
then each special fibre has RM or CM by a field
containing~$K$.
Cf.~Corollary~\ref{sem_con}.

\subsubsection*{Our terminology concerning RM/CM}
In~order to cope with that situation, we use a particular terminology.
If
$\End_\Hg(T) \supseteq K$
then we say that the transcendental part
$T \subset H^2(X, \bbQ)$
is {\em acted upon\/}
by~$K$.
Thus,~if
$T$
is acted upon by
$K \supsetneqq \bbQ$
then
$X$
has RM or CM by a field that
contains~$K$.

\subsubsection*{Semicontinuity of\/
$\End_\Hg(T)$
in families}
Semicontinuity may thus be formulated as follows. If the general fibre of a family
$q\colon \frakX \to Y$
of complex
$K3$~surfaces
has RM or CM by a field
$K$
then, for every special
fibre~$\frakX_t$,
the transcendental part
$T_t \subset H^2(\frakX_t, \bbQ)$
is acted upon
by~$K$.
Moreover,~as Corollary~\ref{sem_con} shows, there is a countable union
$V \subset Y$
of analytic subsets such that
$\frakX_t$
actually has RM or CM
by~$K$
if and only if
$t \in Y \!\setminus\! V$.

\subsubsection*{Periods}
By~a marked
$K3$~surface, we mean a complex
$K3$~surface
$X$
together with an isomorphism
$i\colon \bbZ^{22} \to H^2(X, \bbZ)$,
i.e.\ equipped with a distinguished basis
$(c^1, \ldots, c^{22})$
of
$H^2(X, \bbZ)$.
%Associated~with the marking, one also has the basis
%$(c_1, \ldots, c_{22})$
%that is dual to
%$\smash{(c^1, \ldots, c^{22})}$
%with respect to the cup product pairing. In~this way,
A~marked
$K3$~surface
$(X,i)$
gives rise to the {\em period~point\/}
$$\Pi_{X,i} := ((c^1, [\omega]) : \cdots : (c^{22}, [\omega])) \in \Pb^{21}(\bbC) \,.$$
Here,
$[\omega] \in H^2(X, \bbC)$
denotes the class of a nowhere vanishing holomorphic
$(2,0)$-form.
The~form
$\omega$
is unique up to a constant factor
$\lambda \in \bbC^*$.
Having chosen a particular such form, one may consider the {\em period vector\/}
$((c^1, [\omega]), \ldots, (c^{22}, [\omega])) \in \bbC^{22}$,
which is usually denoted by
$\Pi_{X,i}$,
as~well. The~marking
$i$
is often suppressed from the notation, when there seems to be no danger of~confusion. The~coordinates of the period vector are usually referred to as {\em periods}.

\looseness-1
Moreover, a marking
$i\colon \bbZ^{22} \to H^2(X, \bbZ), e_k \mapsto c^k$,
induces a perfect symmetric bilinear pairing
on~$\bbZ^{22}$,
the pull-back of the cup product pairing. On~the other hand, let such a symmetric bilinear
pairing~$\kappa$
be given
on~$\bbZ^{22}$.
Then~a classical result (cf.~\cite[Chapter~IX]{Sh}) states that the period points of all marked
$K3$~surfaces
that induce
$\kappa$
on~$\bbZ^{22}$
form a (possibly void) open subset
$\Omega_\kappa$
of the quadric
$Q_\kappa$
in~$\Pb^{21}(\bbC)$,
defined
by~$\kappa$.

\subsubsection*{RM in the situation of Picard
rank~$16$}
In~a sufficiently general family
$q\colon \frakX \to Y$
of
$K3$~surfaces
that is generically of geometric Picard
rank~$16$,
and not containing an isotrivial subfamily, the surfaces that are acted upon by a {\em real\/} quadratic field form families over base
curves~$C \subset Y$.
For~arbitrary Picard rank, number fields of larger degree, or CM instead of~RM, the dimensions of the base varieties occurring are known, too. Cf.\ Theorem~\ref{RMCM_dim}, shown~below.

In~coincidence with this, several explicit
\mbox{$1$-dimensional}
families of RM~surfaces, but also isolated examples, have been found. There are two families, for which RM by
$\smash{\bbQ(\sqrt{2})}$
and~$\smash{\bbQ(\sqrt{5})}$,
respectively, is~proven. Cf.~Examples \ref{Qw2_fam} and~\ref{Qw5_fam} or \cite{EJ20}. For~conjectural examples, we refer to \cite[Conjectures~5.2]{EJ16}, cf.\ Remarks \ref{conj_RM}.

On~the side of the period space, the geometry of the corresponding base curves
$\Pi(C) \subset \Pb^{21}(\bbC)$
is very~simple. Concretely,~the restricted period space
$\Omega_{\kappa,16} \subset \Omega_\kappa$,
taking into account the fact that a fixed
\mbox{$16$-dimensional}
subspace of
$H^2(X, \bbQ)$
is contained in the algebraic part, is an open subset of a quadric
in~$\Pb^5(\bbC) \subset \Pb^{21}(\bbC)$.
And~$\Pi(C) \subset \Pb^5(\bbC)$
is just the intersection
of~$\Omega_{\kappa,16}$
with a projective subspace of
codimension~$3$.

\subsubsection*{Explicit construction of tori representing transcendental cohomology classes}
Unfortunately,~the period map
$\Pi\colon Y \to \Omega_{\kappa,16}$
is not too~explicit. In~order to deal with this problem, in Section~\ref{part_fam}, we describe a particular family of
$K3$~surfaces,
which is the one we work~with.

For~the surfaces in this family, we give a topological construction of {\em topological tori,} representing a generating system of the transcendental part of the cohomology. I.e., of continuous mappings
$\alpha\colon T^2 \to X$,
for
$T^2$
the
$2$-dimensional
torus. This~construction is fundamental to the whole~approach.

\begin{cau}
In what follows, we just write {\em torus\/} to mean a topological torus. The reader should, however, keep in mind that
$\alpha$
needs not even to be~injective.
\end{cau}

The tori we construct are not just continuous, but
almost~$C^1$,
cf.\ Definition~\ref{almost_C1}. This~gives us a link to the de Rham cohomology theory, so that the coordinates
$(c_k, [\omega])$
of the period vector turn out to be improper integrals, cf.\ Theorem~\ref{periods_ints}.
This~makes it possible to compute, for a single point
$t \in Y$,
i.e.\ a single
surface~$X = \frakX_t$,
the period point
$\smash{\Pi(t) = \Pi_{\frakX_t}}$
by numerical integration. In~addition, the directional derivatives of
$\Pi$
at~$t$
may be~calculated.

We~describe the construction of the tori in the course of Section~\ref{part_fam}. The~technical details are treated in the final~section.

\begin{rem}
There is a completely different approach to the calculation of periods, due to  E.\,C.\ Sert\"oz~\cite{Ser}. Sert\"oz's method is based on an understanding of the periods under deformation. I.e., he computes the directional derivatives of
$\Pi$
exactly, something we do not do, and derives a partial differential equation for
$\Pi$
from the outcome. This~works for nonsingular hypersurfaces, which includes
$K3$~surfaces
of
degree~$4$,
while our method is, at least at this moment, limited to double covers
of~$\Pb^2$,
ramified in a union of six~lines. Thus,~a fair comparison seems to be difficult to undertake.
\end{rem}

\subsubsection*{The base curve\/
$C$
parametrising RM~surfaces--Numerical continuation}
Given~an isolated example as mentioned above, numerical continuation methods (cf.\ \cite{AG}) may apply in order to extend to example to a
\mbox{$1$-dimensional}
family. One~determines a list of further surfaces in the same family, the period points of which lie
on~$\Pi(C)$
at high~precision. Then,~using numerical linear algebra (in particular the singular value decomposition), one can recover algebraic equations with small coefficients defining a curve that numerically contains the parameters of the surfaces~found.

This~approach should work rather~generally. We~illustrate it at one particular example, which is a
\mbox{$1$-dimensional}
family of
$K3$~surfaces
of
degree~$2$,
for which we report strong evidence for real multiplication
by~$\bbQ(\sqrt{13})$.
The~generic geometric Picard rank is
$16$.
Cf.~Example~\ref{qw13} for the isolated example that was known to us before and Example~\ref{result} the
\mbox{$1$-dimensional}
family~found.

\begin{nota_allg}
We follow standard mathematical conventions and use standard notation. In~particular,

\begin{iii}
\item[ $\bullet$ ]
The standard basis of the free
$R$-module
$R^n$
is denoted by
$(e_1, \ldots, e_n)$.
\item[ $\bullet$ ]
When
$f\colon M \to \bbR^n$
is a map from a set
$M$
to~$\bbR^n$,
we denote by
$f_1, \ldots, f_n$
the components
of~$f$.
I.e.,
$f = (f_1, \ldots, f_n)$.
\item[ $\bullet$ ]
When
$K$
is a field, we identify
$K = \Ab^1(K)$
with a subset
of~$\Pb^1(K)$,
according to the embedding
$K \to \Pb^1(K)$,
$x \mapsto (1:x)$.
Moreover,~we let
$\infty := (0:1)$.
\item[ $\bullet$ ]
For positive integers
$m < n$,
we denote by
$\Pb^m(\bbC) \subset \Pb^n(\bbC)$
the particular
\mbox{$m$-di}\-men\-sional
projective subspace
$\{(x_0:\cdots:x_n) \in \Pb^n(\bbC) \mid x_{m+1} = \cdots = x_n = 0\}$
of~$\Pb^n(\bbC)$.
\item[ $\bullet$ ]
For
$l_1,\ldots,l_r$
global sections of an invertible sheaf
$\calL \in \Pic X$
on a scheme
$X$,
we let
$V(l_1,\ldots,l_r)$
be the common vanishing locus of
$l_1,\ldots,l_r$.
\item[ $\bullet$ ]
When
$f\colon S \to X$
is a continuous map between topological spaces, we denote by
$f^*\colon H^i(X, \bbZ) \to H^i(S, \bbZ)$
the induced homomorphism in cohomology and by
$f_*\colon H_i(S, \bbZ) \to H_i(X, \bbZ)$
that in~homology. Moreover,~as in \cite[Chapter\,1]{Sp}, we use the notation
$f_\#\colon \pi_1(S,\cdot) \to \pi_1(X,\cdot)$
for the induced homomorphism between the fundamental~groups.
\item[ $\bullet$ ]
When
$f\colon S \to X$
is a continuous map between equidimensional, oriented, compact manifolds, we let
$f_!\colon H^i(S, \bbZ) \to H^{i+\dim X-\dim S}(X, \bbZ)$
denote the (cohomological) transfer map. For a definition, cf.~\cite[Chapter~VIII, Definition~10.5]{Dd}.
%\item[ $\bullet$ ]
%Following,~e.g., W.~Fulton~\cite{Fu}, we distinguish between a curve
%$D$
%on an algebraic surface
%$X$
%and the Weil divisor
%$[D] \in \Div X$
%it~induces.
\item[ $\bullet$ ]
By a manifold, we always mean a smooth
(i.e.\
$C^\infty$)
manifold.
\item[ $\bullet$ ]
When
$X$
is a complex
$K3$~surface,
we usually identify the Picard group
$\Pic X$
with its image in
$H^2(X, \bbZ)$,
under the Chern class homomorphism.
\end{iii}
\end{nota_allg}

\begin{comp_allg}
All computations are done with {\tt magma}~\cite{BCP} on one core of an Intel i7-7700 processor running at 3.6\,GHz.
\end{comp_allg}

\section{$K3$~surfaces}

\subsubsection*{The algebraic setting}
Let~$k$
be an algebraically closed field. Then
$K3$~surfaces
form one of the types of surfaces
over~$k$,
according to the Enriques--Kodaira classification of surfaces~\cite{CE,Sh}. 
This~classification is due to the Italian school of Algebraic Geometry, cf.~\cite{CE,Sev}, and today called the Enriques--Kodaira classification. More~recent treatments are given in~\cite{Sh}, \cite{Be}, and~\cite{BHPV}.
By~a {\em surface,} here, one means a
\mbox{$2$-dimensional},
connected, regular scheme that is projective
over~$k$.

More~concretely,
$K3$~surfaces
form one of the four types of surfaces of Kodaira dimension
$0$
and, among them, one of the two types of surfaces having trivial canonical~class. The~other such type are abelian surfaces.
A~surface
$X$
is characterised (\cite[Theorem~VIII.2]{Be} and \cite[Corollary~VIII.8.6]{BHPV}) to
be~$K3$
by the properties~that

\begin{iii}
\item
$K_X \cong \calO_X$
and
\item
$\pi_1^\et(X,\cdot) = 0$.
\end{iii}

\noindent
Moreover,~$H^2_\et(X, \bbZ_l)$
is free of
rank~$22$
and, of course,
$H^4_\et(X, \bbZ_l) \cong H^0_\et(X, \bbZ_l) \cong \bbZ_l$
as well as
$H^3_\et(X, \bbZ_l) \cong H^1_\et(X, \bbZ_l) = 0$
[under the assumption that
$l \neq \cha k$].

\begin{rems}
\begin{abc}
\item
The degree of a
$K3$~surface
embedded in ordinary projective space is automatically even. There~are examples in every even degree
$d \geq 4$
\cite[Proposition~VIII.15]{Be}.
\begin{iii}
\item
Smooth quartics
in~$\Pb^3$
are
$K3$~surfaces
of
degree~$4$.
\item
Smooth complete intersections of a quadric and a cubic
in~$\Pb^4$
are
$K3$~surfaces
of
degree~$6$.
\item
Smooth complete intersections of three quadrics
in~$\Pb^5$
are
$K3$~surfaces
of
degree~$8$.
\end{iii}
\item
Similarly,~double covers
of~$\Pb^2$
ramified in a smooth sextic curve are
$K3$~surfaces
of
degree~$2$.
\item
In~each of the cases mentioned in a) and~b), isolated singular points may be allowed, as long as all of them are of type
$AD\hspace{-0.3mm}E$~\cite[Section~8.2.7]{Do}.
Then~a
$K3$~surface
is obtained as the minimal desingularisation of the surface~described.
\end{abc}
\end{rems}

\subsubsection*{The setting of complex manifolds}

There is a theory, parallel to the classification just mentioned, for complex~surfaces. I.e.,~for
\mbox{$2$-dimensional}
compact complex~manifolds. In~this setting, a
{\em $K3$~surface\/}
is a complex surface such that

\begin{iii}
\item
there is a global holomorphic
\mbox{$(2,0)$-form}
$\omega \in \Omega^{2,0}(X)$
without zeros (or~poles), and
\item
$\pi_1(X,\cdot) = 0$.
\end{iii}

\noindent
A~complex
$K3$~surface
is automatically a K\"ahler manifold~\cite[Theorem~IV.3.1]{BHPV}. In~general, however, it is {\em not\/} projective so that this theory does not exactly mirror the algebraic theory of
$K3$~surfaces
in the case
that~$k = \bbC$,
but is more~general.

All
$K3$~surfaces
$X$
are mutually diffeomorphic~\cite[Corollary~VIII.8.6]{BHPV}.
Moreover,
$H^2(X, \bbZ)$
is free of
rank~$22$
and, of~course,
$H^4(X, \bbZ) \cong H^0(X, \bbZ) \cong \bbZ$
as well as
$H^3(X, \bbZ) \cong H^1(X, \bbZ) = 0$.
Consequently,
$H_2(X, \bbZ) \cong \bbZ^{22}$,~too.

\section{The period space}

\subsubsection*{Explicit cohomology classes}
A complex
$K3$~surface
$X$
has a natural orientation, defined by the complex~structure. Consider~a compact, oriented
\mbox{$2$-manifold}
$S$,
together with a continuous map
$\alpha\colon S \to X$.
Then~$\alpha$
defines a cohomology~class
$$c_\alpha := \alpha_!(1) \in H^2(X, \bbZ) \,,$$
for
$1 \in H^0(S, \bbZ)$
the canonical~generator.

\begin{defi}
\label{cohclass_geom}
We call the class
$c_\alpha \in H^2(X, \bbZ)$
the {\em cohomology class given by\/}
$S$,
{\em together with the continuous map\/}
$\alpha\colon S \to X$.
\end{defi}

\begin{rem}
The class
$c_\alpha \in H^2(X, \bbZ)$
is the same as the Poincar\'e dual of the homology class of the image of
$S$
in~$X$.
Cf.\ Lemma~\ref{homcohom}.b),~below.
\end{rem}

\begin{rems}
\begin{iii}
\item
Note that the continuous map
$\alpha\colon S \to X$
might be the embedding of a divisor, the embedding of a non-holomorphic submanifold, or not even an~embedding.
\item
In~the particular case that
$i_D\colon D \to X$
is the embedding of a divisor, we prefer the more classical notation
$[D] := c_{i_D} \in H^2(X, \bbZ)$.
\end{iii}
\end{rems}

\begin{lem}[Cohomology versus homology]
\label{homcohom}
\leavevmode

\begin{abc}
\item
There~is the canonical~isomorphism\/
$$\iota\colon H^2(X, \bbZ) \to H_2(X, \bbZ) \,, \quad u \mapsto u \cap z_X \,,$$
for\/
$z_X \in H_4(X, \bbZ)$
the fundamental~class.
\item
Let\/~$\alpha\colon S \to X$
be as above. Then,
for\/~$u \in H^2(X, \bbC)$
arbitrary,
$$\langle u \cup c_\alpha, z_X \rangle = \langle \alpha^*(u), z_S \rangle \,,$$
where\/
$z_S \in H_2(S, \bbZ)$
denotes the fundamental class
of\/~$S$.
Moreover,
$$\iota(c_\alpha) = \alpha_*(z_S) \,.$$
\end{abc}%\smallskip

\noindent
{\bf Proof.}
{\em
a)
This is a version of Poincar\'e duality \cite[Chapter~6, Theorem~3.12]{Sp}.\smallskip

\noindent
b)
The cohomological transfer map
$\alpha_!$
is characterised \cite[Chapter~VIII, Definition~10.5]{Dd} by the property that
$\alpha_!(a) \cap z_X = \alpha_* (a \cap z_S)$,
for every
$a \in H^0(S, \bbZ)$.
The~second assertion is just a particular case of~this.

For the first one, the claim is simply
$(u \cup c_\alpha) \cap z_X = \alpha_*(\alpha^*(u) \cap z_S)$.
But~the term on the left hand side is
$u \cap (c_\alpha \cap z_X) = u \cap \iota(c_\alpha) = u \cap \alpha_*(z_S)$,
due~to \mbox{\cite[Chapter~5, 6.18]{Sp}}, while that to the right is
$u \cap \alpha_*(z_S)$,
by \mbox{\cite[Chapter~5, 6.16]{Sp}}.%
}%
\eop
\end{lem}

\begin{rems}
\begin{iii}
\item
In what follows, we generally prefer cohomology versus homology, although there are situations, in which homology classes appear to be more~natural. For~instance, as
$K3$~surfaces
are simply connected, the Hurewicz isomorphism theorem~\cite[Chapter~7, Proposition~5.2]{Sp} shows that
$\smash{\pi_2(X,\cdot) \cong H_2(X,\bbZ)}$.
Thus,~every class in 
$H_2(X,\bbZ)$,
and hence every class in
$H^2(X, \bbZ)$,
may be represented by a spheroid
$S^2 \to X$.
Therefore, it can be represented by a torus
$\alpha\colon T^2 \to X$,~too.
\item
Lemma~\ref{homcohom}.a) does not hold for singular surfaces. Cf.~Proposition \ref{h_doubcov},~below.
\item
There is the {\em cup product pairing\/}
$(\cdot\,,\cdot) := \langle \cdot\cup\cdot, z_X\rangle\colon H^2(X, \bbZ) \times H^2(X, \bbZ) \to \bbZ$,
which is symmetric, bilinear, and~perfect. The~cup product pairing bilinearly extends to a~pairing
$$(\cdot\,,\cdot)\colon H^2(X, \bbZ) \times H^2(X, \bbC) \to \bbC \,.$$
\end{iii}\smallskip
\end{rems}

\subsubsection*{Marked
$K3$~surfaces}

\begin{defi}
By a {\em marked\/
$K3$~surface,}
we mean a complex
$K3$~surface
together with an isomorphism
$i\colon \bbZ^{22} \longrightarrow H^2(X, \bbZ)$.
\end{defi}

\begin{nota}
\label{mark}
\begin{iii}
\item
The {\em marking\/}
$i\colon \bbZ^{22} \to H^2(X, \bbZ)$
determines the cohomology~classes
$$c^k := i(e_k) \in H^2(X, \bbZ) \,,$$
for~$k = 1,\ldots,22$,
which form a basis
of~$H^2(X, \bbZ)$.
\item
There~is the dual basis
$(c_1, \ldots, c_{22})$
of
$H^2(X, \bbZ)$,
given by the condition that
$(c_k, c^j) = \delta_{kj}$,
for
$( .\,,. )\colon H^2(X, \bbZ) \times H^2(X, \bbZ) \to \bbZ$
the cup product pairing.

We denote the isomorphism
$\bbZ^{22} \to H^2(X, \bbZ), e_k \to c_k$,
by~$\check{\imath}$.
\item
Moreover,~the marking defines a perfect, symmetric pairing
on~$\bbZ^{22}$,
namely the pull-back of the cup product pairing via
$\check{\imath}$.
\end{iii}
\end{nota}

\noindent
Given any cohomology class
$u \in H^2(X, \bbC)$,
there is the {\em tautological decomposition}
$$u = (c^1, u) c_1 + \cdots + (c^{22}, u) c_{22} \,,$$
of~$u$
{\em with respect
to~$i$}.
For~instance, this applies to the {\em distinguished\/} cohomology class
$[\omega] \in H^2(X, \bbC)$
of the
$K3$~surface
$X$,
which is uniquely determined up to~scaling. Here~and in what follows, we denote
by~$\omega$
a nowhere vanishing holomorphic
\mbox{$(2,0)$-form}
on~$X$.

\begin{defi}
A~marked
$K3$~surface
$(X,i)$
gives rise to a point
$$\Pi_{X,i} := ((c^1, [\omega]) : \cdots : (c^{22}, [\omega])) \in \Pb^{21}(\bbC) \,,$$
which is called the {\em period point\/} of
$(X,i)$.
\end{defi}

\begin{lem}
\label{indep_lift}
Let\/~$X$
be a complex\/
$K3$~surface.
Then~the period point\/
$\Pi_{X,i}$
is completely determined by the composition\/
$\smash{\bbZ^{22} \stackrel{i}{\longrightarrow} H^2(X, \bbZ) \twoheadrightarrow H^2(X, \bbZ)/\Pic X}$.\smallskip

\noindent
{\bf Proof.}
{\em
Suppose that
$i$
and~$i'$
are two markings
on~$X$
differing only by elements
in~$\Pic X$.
I.e., that
$c'{}^k - c^k \in \Pic X$,
for
$k = 1,\ldots,22$.
Then
$$(c'{}^k, [\omega]) - (c^k, [\omega]) = (c'{}^k - c^k, [\omega]) = 0 \,,$$
since
$c'{}^k - c^k \in \Pic X$
is of type
$(1,1)$,
while
$[\omega]$
is of type
$(2,0)$.
The assertion~follows.
}
\eop
\end{lem}

\begin{theo}[I.\,R.~Shafarevich et al.]
\label{Shafa}
Let\/~$\kappa$
be a perfect, symmetric pairing
on\/~$\bbZ^{22}$.
Denote~the\/
$\bbC$-bilinear
pairing induced
by\/~$\kappa$
on\/~$\bbC^{22}$,
by\/~$\kappa$
as~well.

\begin{abc}
\item
Let\/
$(X_1,i_1)$
and\/
$(X_2,i_2)$
be marked
$K3$~surfaces
defining the pairing\/
$\kappa$
on\/~$\bbZ^{22}$
and having the same period~point.
Then\/
$(X_1,i_1)$
and\/
$(X_2,i_2)$
are~isomorphic.
\item
The set\/
$\Omega_\kappa$
of the period points of all marked
$K3$~surfaces
defining the
pairing\/~$\kappa$
is possibly empty. Otherwise, one has
$$\Omega_\kappa = \{(x_1:\cdots:x_{22}) \in \Pb^{21}(\bbC) \mid \kappa(x,x) = 0, \kappa(x,\overline{x}) > 0 \} \,.$$
This is an open subset of the quadric\/
$Q_\kappa \subset \Pb^{21}(\bbC)$,
defined
by\/~$\kappa$.
\end{abc}%\smallskip

\noindent
{\bf Proof.}
{\em This is \cite[Theorem~VIII.11.1]{BHPV}, together with \cite[Theorem VIII.\discretionary{}{}{}14.1]{BHPV}. Cf.~\cite[Chapter~IX]{Sh}.
}
\eop
\end{theo}

\subsubsection*{The relative situation}

Let~$q\colon \frakX \to Y$
be a holomorphic submersion of complex manifolds and suppose that every fibre is a
$K3$~surface.
Then~the higher direct image sheaf
$R^2 q_* \bbZ$
on~$Y$
is locally constant of
dimension~$22$.
Moreover,~due to Grauert's Theorem (\cite[Satz~5]{Gr}), one has a natural base change~isomorphism
$\smash{\iota_t^*\colon (R^2 q_* \bbZ)_t \stackrel{\cong}{\longrightarrow} H^2(\frakX_t, \bbZ)}$,
for
every~$t \in Y$.

\begin{defi}
\label{fam_marked}
By a {\em family of marked
$K3$~surfaces,}
one means a holomorphic submersion
$q\colon \frakX \to Y$,
every fibre of which is a 
$K3$~surface,
together with a system
$\fraki = (i_t)_{t\in Y}$
of markings satisfying the following~conditions.

\begin{iii}
\item
For each
$t\in Y$,
$i_t$
is a marking on the
fibre~$\frakX_t$.
\item
For~every
$v \in \bbZ^{22}$,
$t \mapsto (\iota_t^*)^{-1} \!\circ\! i_t(v) \in (R^2 q_* \bbZ)_t$
is a global section
of~$R^2 q_* \bbZ$.
I.e., the markings
$i_t$
trivialise
$R^2 q_* \bbZ$~globally.
\end{iii}
\end{defi}

%\begin{rems}
%\begin{iii}
%\item
%Suppose that the base
%manifold~$Y$
%is~connected. Then~the
%system~$\fraki$
%of markings is uniquely determined by just one
%making~$i_t$.
%\item
%Furthermore, assume that
%$Y$
%is simply connected, e.g.\ a~polydisc. Then, automatically,
%$R^2 q_* \bbZ \cong \bbZ^{22}$
%is constant of
%dimension~$22$.
%In~particular, every marking
%$i_t$
%on a fibre extends to the whole~family.
%\end{iii}
%\end{rems}

\begin{lem}
\label{period_holom}
Let\/~$q\colon (\frakX, \fraki) \to Y$
be a family of marked
$K3$~surfaces.
Then~the period mapping\/
$\Pi\colon Y \to \Pb^{21}(\bbC)$,
$t \mapsto \Pi_{\frakX_t, i_t}$,
is~holomorphic.\smallskip

\noindent
{\bf Proof.}
{\em
This~is easily verified by a direct calculation, cf.~\cite[Theorem~IV.4.3]{BHPV}.%
}%
\eop\smallskip
\end{lem}

\subsubsection*{The restricted period space}

\begin{lem}
\label{periods_Pic}
Let\/
$(X,i)$
be a marked\/
$K3$~surface,
$x = \Pi_{X,i} \in \Pb^{21}(\bbC)$
its period point, and\/
$\kappa$
the perfect, symmetric pairing
on\/~$\bbZ^{22}$
defined
by\/~$(X,i)$.
Denote~the pairing induced
on\/~$\bbC^{22}$
by\/~$\kappa$,
too. Then\/
$$\check{\imath}^{-1} \Pic X = \spann(x)^\perp \cap \bbZ^{22} \,,$$
the orthogonal complement being taken
in\/~$\bbC^{22}$
with respect
to\/~$\kappa$.\smallskip

\noindent
{\bf Proof.}
{\em
One has
$\Pic X \cong \{u \in H^2(X, \bbZ) \mid u \cup [\omega] = 0\}$,
according to the Lefschetz theorem on
\mbox{$(1,1)$-classes}~\cite[p.~163]{GH}.
Moreover,
\begin{eqnarray*}
\langle u \cup [\omega], z_X\rangle = \Big\langle\!\smash{\Big(\sum_k (c^k, u) c_k\Big) \cup \Big(\sum_j x_j c_j\Big)}, z_X\Big\rangle &=& \kappa\big((c^k, u)_{k=1,\ldots,22}, (x_j)_{j=1,\ldots,22}\big) \\
 &=& \kappa\big((c^k, u)_{k=1,\ldots,22}, x\big) \,.
\end{eqnarray*}
Thus,
$\check{\imath}^{-1} \Pic X = \{v \in \bbZ^{22} \mid \kappa(v,x) = 0\} = \spann(x)^\perp \cap \bbZ^{22}$,
as~required.
\eop
}
\end{lem}

\begin{coro}
\label{periods_quadric}
Let\/~$r \in \{1,\ldots,20\}$
and\/
$\kappa$
a perfect, symmetric pairing
on\/~$\bbZ^{22}$.
Denote~the pairing induced
on\/~$\bbC^{22}$
by\/~$\kappa$,
too. Then~the set\/
$\Omega_{\kappa,r}$
of the period points of all marked
$K3$~surfaces\/
$(X, i)$
such~that

\begin{iii}
\item
the cohomology classes\/
$c^{22-r+1}, \ldots, c^{22} \in H^2(X, \bbZ)$
are algebraic, i.e.\ contained in\/
$\Pic X \subset H^2(X, \bbZ)$,~and
\item
the pairing
on\/~$\bbC^{22}$
induced
by\/~$(X,i)$
is
exactly\/~$\kappa$,
\end{iii}

\noindent
is either void or one~has
\begin{equation}
\label{restr_per}
\Omega_{\kappa,r} = \{(x_1:\cdots:x_{22-r}:0:\cdots:0) \in \Pb^{21}(\bbC) \mid \kappa(x,x) = 0, \kappa(x,\overline{x}) > 0 \} \,.
\end{equation}
This is an open subset of a quadric\/
$Q_{\kappa,r} \subset \Pb^{21-r}(\bbC) \subset \Pb^{21}(\bbC)$.\smallskip

\noindent
{\bf Proof.}
{\em
In comparison with Theorem~\ref{Shafa}, 
$c^{22-r+1}, \ldots, c^{22} \in H^2(X, \bbZ)$
being algebraic is the only additional condition. By~Lemma~\ref{periods_Pic}, this means precisely~that
$$x \!\perp\! \check{\imath}^{-1}(c^{22-r+1}), \;\ldots,\; x \!\perp\! \check{\imath}^{-1}(c^{22})$$
with respect to
$\kappa$.
I.e., that
$\check{\imath}(x) \!\perp\! \spann(c^{22-r+1}, \ldots, c^{22})$
with respect to the cup product pairing. The~latter is equivalent to
$\check{\imath}(x) \in \spann(c_1, \ldots, c_{22-r})$,
which exactly means that
$x \in \spann(e_1,\ldots,e_{22-r})$.
\eop
}
\end{coro}

\begin{defi}
Let
$q\colon (\frakX, \fraki) \to Y$
be a family of marked
$K3$~surfaces
of the kind that
$i_t(e_{22-r+1}), \ldots, i_t(e_{22})$
are algebraic, for every
$t \in Y$.
Then~we call
$$\Pi\colon Y \to \Pb^{21-r}(\bbC) \,, \quad t \mapsto \Pi_{\frakX_t, \fraki_t} := ((i_t(e_1), [\omega]) : \cdots : (i_t(e_{22-r}), [\omega])) \,,$$
the {\em restricted period map}. In~the case that
$Y$
is a point, we speak of the {\em restricted period~point}.
\end{defi}

\begin{theo}[A basis of the transcendental part induces a class of markings]
\leavevmode
\label{mark_rel}

\noindent
Let\/~$r \in \{1,\ldots,20\}$,
$X$
a complex\/
$K3$~surface,
$P \subseteq \Pic X$
a subgroup that is co\-tor\-sion-free, and\/
$\iota\colon \bbZ^{22-r} \to H^2(X, \bbZ)/P$
an isomorphism.

\begin{abc}
\item
Then
\begin{itemize}
\item[{\rm i) }]
lifting\/~$\iota$
to a homomorphism\/
$\bbZ^{22-r} \to H^2(X, \bbZ)$
and\smallskip
\item[{\rm ii) }]
extending that to a homomorphism\/
$i\colon \bbZ^{22} \to H^2(X, \bbZ)$
by choosing an arbitrary basis\/
$(u_{22-r+1}, \ldots, u_{22})$
of\/~$P$
and
putting\/~$i(e_k) := u_k$,
for\/
$k = 22-r+1, \ldots, 22$,
\end{itemize}
provides a
marking\/~$i$
on\/~$X$.
\item
The
marking\/~$i$
yields a restricted period point\/
$\Pi_X \in \Pb^{21-r}(\bbC)$.
\item
The restricted period point\/
$\Pi_X \in \Pb^{21-r}(\bbC)$
in independent of the choices made in steps~a.i) and~a.ii).
\end{abc}%\smallskip

\noindent
{\bf Proof.}
{\em
The only assertion in~a) is that
$i$~is
an isomorphism, which is~clear. b)~follows, as the assumptions imply that
$i(e_{22-r+1}), \ldots, i_t(e_{22})$
are algebraic. Finally, c) is a direct consequence of Lemma~\ref{indep_lift}.
}
\eop
\end{theo}

\subsubsection*{Periods as integrals}

Suppose that, on a marked
$K3$~surface
$(X,i)$,
the cohomology class
$c^j$
is given by a compact, oriented
\mbox{$2$-manifold}
$S$,
together with a continuous map
$\alpha\colon S \to X$.
Then the
$j$-th
period~is
$$(c^j, [\omega]) = (c_\alpha, [\omega]) = \langle c_\alpha \cup [\omega], z_X\rangle = \langle \alpha^*([\omega]), z_S\rangle \,,$$
by Lemma~\ref{homcohom}.b). For~smooth
$\alpha$,
this is simply
$\int_S \alpha^*([\omega])$.
Thus, one might have the idea to avoid the topological machinery and to consider the periods just as integrals of the nowhere vanishing
$(2,0)$-form.

While we indeed use a representation as an integral for computational purposes, cf.\ Theorem~\ref{periods_ints} below, such a representation is not well suitable as a~definition. One~problem is that, in our situation, the map
$\alpha$
is usually non-smooth. Also,~and more crucially, the application of the periods we have in mind is towards real multiplication, which is defined in terms of cohomology.

\section{Explicit description of transcendental cohomology classes}
\label{part_fam}

\subsubsection*{The family}

We~consider double covers
of~$\Pb^2_\bbC$,
ramified in a union of six~lines. These~are given by
$X'\colon w^2 = l_1 \cdots l_6$,
for
$l_1, \ldots, l_6$
linear forms in three~variables. We~assume that no point of the plane is contained in three~lines. Then~there~are 15 singular points, each of which is isolated and conical, i.e.\ of
type~$A_1$.
The~minimal desingularisation
$X$
of~$X'$
is hence a
$K3$~surface.

Concretely,~the
$K3$~surface
$X$
is obtained
from~$X'$
by blowing up the 15 singular~points. This~yields 15 exceptional curves, each of which is a
\mbox{$(-2)$-curve}.
%as the singular points
%on~$X'$
%are~conical.
In~particular,
$\rk\Pic X \geq 16$.

\begin{nota}
\label{nota_sec4}
\begin{iii}
\item
For~$a = 1, \ldots, 6$,
we denote the coefficients of the linear form
$l_a$
by
$A_{a1}$,
$A_{a2}$,
and~$A_{a3}$.
I.e.,~$l_a = A_{a1} x + A_{a2} y + A_{a3} z$.
\item
We write
$x_{ij} := V(l_i,l_j) \in \Pb^2(\bbC)$,
for
$1 \leq i, j \leq 6$,
for the 15 double points of the branch~locus. Accordingly,~we denote the 15 exceptional curves by
$E_{ij}$,
$1 \leq i, j \leq 6$.
\item
We~use the notation
$\pi\colon X \to \Pb^2(\bbC)$
for the natural~map. By~definition,
$\pi$
factors
via~$X'$,
$$\pi\colon X \stackrel{\bl}{\longrightarrow} X' \stackrel{\pi_{X'}}{\longrightarrow} \Pb^2(\bbC) \, .$$
Similarly,~as is well-known,
$\pi$
factors via the blowing-up
of~$\Pb^2(\bbC)$
in the 15 double points of the branch~locus,
$$\pi\colon X \stackrel{\pi'}{\longrightarrow} \Bl_{x_{12},\ldots,x_{56}}(\Pb^2(\bbC)) \longrightarrow \Pb^2(\bbC) \,.$$
The~map
$\pi'$~is
a double cover, ramified in the union of six mutually disjoint projective~lines.
\item
Moreover,~$X$
carries the natural involution
$\zeta\colon X \to X$,
which is induced by the map
$\zeta'\colon X' \to X'$,
$(w;\,x:y:z) \mapsto (-w;\,x:y:z)$.
\item
We~put
$$P := V \cap H^2(X, \bbZ) \,,$$
for
$V \subset H^2(X, \bbQ)$
the subvector space spanned by the classes
$\smash{[E_{ij}]}$,
for
$1 \leq i,j \leq 6$,
together with
$\pi^*[l]$,
the inverse image of a general line
on~$\Pb^2$.

According~to this definition,
$P \subset H^2(X, \bbZ)$
is a saturated sublattice of
rank~$16$.
Let~us note explicitly that
$P$
is generated by the classes
$\smash{[E_{ij}]}$
and~$\pi^*[l]$
only up to finite~index.
\end{iii}
\end{nota}

\begin{rems}
\begin{iii}
\item
It~is a well-known fact from projective geometry that configurations of six lines
in~$\Pb^2$
have exactly four~moduli. Indeed,~one might normalise the six lines under the operation of
$\Aut(\Pb^2_\bbC) = \PGL_3(\bbC)$
in such a way that
$X'$
is given by an equation of the~form
\begin{equation}
\label{eq_norm}
X' = X'_{(a_0,b_0,c_0,d_0)}\colon w^2 = xyz(x+y+z)(x+a_0y+b_0z)(x+c_0y+d_0z) \,.
\end{equation}
%We~shall, however, not do so in the theory that~follows.
\item
The restricted period space corresponding to a system of generators
of~$P$
is an open subset of a
quadric~$Q \subset \Pb^5$,
which is in perfect coincidence with the fact that geometrically there are four~moduli.
\end{iii}
\end{rems}

\begin{lem}
\label{eigenv_invol}
Let\/~$X$
be a\/
$K3$~surface
of the family~described.

\begin{abc}
\item
Then~the
involution\/~$\zeta$
acts on\/
$P \!\otimes_\bbZ\! \bbQ \subset H^2(X, \bbQ)$
as the identity map and on\/
$P^\perp \subset H^2(X, \bbQ)$
as the multiplication
by\/~$(-1)$.
\item
In~particular,
$\zeta$
operates
on\/~$H^2(X, \bbZ)/P$
as the multiplication
by\/~$(-1)$.
\end{abc}\medskip

\noindent
{\bf Proof.}
{\em
a)
The
involution~$\zeta$
has a disjoint union of six projective lines as its fixed point~set. The~topological Euler characteristic of such a configuration amounts
to~$12$.
Therefore,~the Lefschetz trace formula~\cite[Theorem~8.5]{Ed}, cf.~\cite[Expos\'e~III, formule~(4.11.3)]{SGA5}, shows~that
$\Tr \zeta |_{H^2(X, \bbQ)} = 10$.
In~other words, the
eigenvalue~$1$
has
multiplicity~$16$,
while the
eigenvalue~$(-1)$
occurs with
multiplicity~$6$.

On~the other hand,
$\zeta$~clearly
fixes
$P$
pointwise, which shows that indeed
$P \!\otimes_\bbZ\! \bbQ$
is the eigenspace
for~$1$.
Finally,~the operation
of~$\zeta$
is compatible with the cup product pairing. Thus,
for~$x \in P$
and
$c \in H^2(X, \bbQ)$
in the eigenspace
for~$(-1)$,
one~finds
$$(x,c) = (\zeta^*x, \zeta^*c) = (x,-c) = -(x,c)$$
and hence
$c \in P^\perp$.
This~completes the proof of~a).\smallskip

\noindent
b)
directly follows from~a).
}
\eop
\end{lem}

\begin{lem}
\label{h_blow}
Let\/~$Y'$
be a complex space that is equidimensional of dimension two and has an isolated singular point\/
$y \in Y'$
of
type\/~$A_1$,
i.e.\ an ordinary double~point. Denote~by\/
$b\colon Y \to Y'$
the blowing-up of\/
$Y'$
at\/~$y$
and by\/
$E \subset Y$
the exceptional~curve.\smallskip

\noindent
Then~the homomorphism\/
$b_*\colon H_2(Y,\bbZ) \to H_2(Y',\bbZ)$
is surjective with kernel\/
$\spann_\bbZ(\iota([E]))$.
I.e.,
$b_*$
induces a natural~isomorphism
$$H_2(Y,\bbZ)/\spann_\bbZ(\iota([E])) \stackrel{\cong}{\longrightarrow} H_2(Y',\bbZ) \,. \smallskip$$

\noindent
{\bf Proof.}
{\em
There~is the following commutative diagram of Mayer-Vietoris exact sequences
$$
\xymatrixcolsep{2.7mm}
\xymatrix{
H_2(U(E) \!\setminus\! E, \bbZ) \ar@{->}[r] \ar@{=}[d] & H_2(Y \!\setminus\! E,\bbZ) \oplus H_2(U(E),\bbZ) \ar@{->}[r]\ar@{->}[d] & H_2(Y,\bbZ) \ar@{->}[r]\ar@{->}[d]^{b_*} & H_1(U(E) \!\setminus\! E, \bbZ) \ar@{=}[d] \phantom{\,.} \\
H_2(U(y) \!\!\setminus\!\! \{y\}, \bbZ) \ar@{->}[r]& H_2(Y' \!\setminus\!\! \{y\},\bbZ) \oplus H_2(U(y),\bbZ) \ar@{->}[r]& H_2(Y'\!,\bbZ) \ar@{->}[r]& H_1(U(y) \!\!\setminus\!\! \{y\}), \bbZ) \,.
}
$$
Here,~we write
$U(y)$
for the intersection
of~$Y'$
with a tiny ball around the point
$y$
and
$U(E) := b^{-1}(U(y))$
for its preimage
in~$Y$.
Clearly,~the natural map
$$b |_{Y \setminus E}\colon Y \!\setminus\! E \to Y' \!\setminus\! \{y\}$$
is a homeomorphism. Moreover,~one readily sees that
$U(y)$
is homotopy equivalent to an affine quadratic cone, and hence contractible, while
$U(E)$
is of the homotopy type
of~$\Pb^1(\bbC) \cong S^2$.
Finally,
$U(y) \!\setminus\! \{y\}$
has the homotopy type of a complex conic, i.e.\ that of
$\Pb^1(\bbC) \cong S^2$,~too.
In~particular, we find that
$H_1(U(y) \!\setminus\! \{y\}, \bbZ) = 0$.
The~commutative diagram therefore takes the form~below,
$$
\xymatrix{
H' \ar@{->}[r] \ar@{=}[d] & H'' \oplus \iota([E]) \!\cdot\!\bbZ \ar@{->}[r]\ar@{->>}[d] & H_2(Y,\bbZ) \ar@{->}[r]\ar@{->}[d]^{b_*} & 0 \phantom{\,.} \\
H' \ar@{->}[r]& H'' \ar@{->}[r]& H_2(Y',\bbZ) \ar@{->}[r]& 0 \,.
}$$
The~assertion follows from this by a simple diagram~chase.
}
\eop
\end{lem}

\begin{prop}[Homology of the double cover]
\label{h_doubcov}
Let\/~$X$
be a\/
$K3$~surface
of the family described and\/
$X'$
the underlying double cover
of\/~$\Pb^2(\bbC)$.

\begin{abc}
\item
Then~the homomorphism\/
$\bl_*\colon H_2(X,\bbZ) \to H_2(X',\bbZ)$
is surjective with kernel\/
$\spann_\bbZ(\iota([E_{12}]), \ldots, \iota([E_{56}]))$.
I.e.,
$\bl_*$
induces a natural~isomorphism
$$H_2(X,\bbZ)/\spann_\bbZ(\iota([E_{12}]), \ldots, \iota([E_{56}])) \stackrel{\cong}{\longrightarrow} H_2(X',\bbZ) \,. $$
\item
In~particular, the natural homomorphism\/
$H^2(X, \bbZ) \twoheadrightarrow H^2(X, \bbZ)/P$
factors~via
$$H^2(X, \bbZ) \stackrel{\iota}{\longrightarrow} H_2(X,\bbZ) \stackrel{\bl_*}{\longrightarrow} H_2(X',\bbZ) \,.$$
\end{abc}%\medskip

\noindent
{\bf Proof.}
{\em
Assertion~a) follows applying Lemma~\ref{h_blow} repeatedly, 15~times in~total.
b)~is a direct consequence of~a).
}
\eop
\end{prop}

\begin{rem}
Proposition~\ref{h_doubcov} shows, in particular, that
$H_2(X',\bbZ)$
has~torsion. Indeed, the cohomology class
$[E_{12}] + \cdots + [E_{16}] + 5\pi^*[l] = [\div(G)] \in H^2(X, \bbZ)$,
for
$G := l_2 \cdots l_6 - l_1^5$,
is divisible
by~$2$.
The~reason is that,
modulo~$G$,
the equation
$w^2 = l_1 \cdots l_6$
of the surface goes over into
$w^2 = l_1^6$.
Thus,~the divisor
$\div(G)$
splits into two components that are interchanged under the
involution~$\zeta$
and the claim follows from Lemma~\ref{eigenv_invol}.a). Consequently,
$[E_{13}] + \cdots + [E_{16}] - ([E_{23}] + \cdots + [E_{26}]) \in H^2(X, \bbZ)$
is divisible
by~$2$,
too, and the same is true for its image
under~$\iota$
in~$H_2(X,\bbZ)$.

The~nontrivial torsion classes
in~$H_2(X',\bbZ)$
are certainly worth being considered more closely, but, in this article, we do not have any applications for~them.
\end{rem}

\begin{coro}
\label{lift_coh_modP}
Let\/~$\alpha\colon T^2 \to X'$
be a torus and\/
$\smash{\underline\alpha, \underline{\underline\alpha}\colon T^2 \to X}$
continuous lifts
of\/~$\alpha$.
Then\/~$\smash{c_{\underline\alpha} - c_{\underline{\underline\alpha}} \in P \subset H^2(X, \bbZ)}$
(cf.\ the notation introduced in \ref{nota_sec4}.v).\smallskip

\noindent
{\bf Proof.}
{\em
According~to Proposition~\ref{h_doubcov}.b), it is sufficient to verify that
$$\smash{\bl_*(\iota(c_{\underline\alpha})) = \bl_*(\iota(c_{\underline{\underline\alpha}})) \in H_2(X',\bbZ) \,.}$$
But this is clear, since
$\iota(c_{\underline\alpha}) = \underline\alpha_*(z_{T^2})$,
by Lemma~\ref{homcohom}.b), and therefore
$$\bl_*(\iota(c_{\underline\alpha})) = \bl_*(\underline\alpha_*(z_{T^2})) = (\bl \circ \underline\alpha)_*(z_{T^2}) = \alpha_*(z_{T^2}) \,,$$
while exactly the same holds for
$\smash{\bl_*(\iota(c_{\underline{\underline\alpha}}))}$.
}
\eop
\end{coro}

%\subsubsection*{Transcendental\/
%$2$-cycles---Toroids
%representing cohomology classes}~%

\subsubsection*{Assumptions on the branch locus}

\begin{ass}[on the branch locus]
\begin{Oiii}
\item
We~assume once and for all that the branch locus is the union of six {\em real\/} lines
$V(l_1), \ldots, V(l_6)$
[no three of which have a point in common]
in~$\Pb^2_\bbC$.
I.e.,~that
$A_{a1}, A_{a2}, A_{a3} \in \bbR$,
for~$a = 1, \ldots, 6$.
\item
We~suppose that any two of the~vectors
$$(A_{11}, A_{12}), \ldots, (A_{61}, A_{62})$$
are linearly~independent. I.e.,~it is assumed that the open chart given
by~``$z=1$'', which is the one we are going to work with, contains the 15 double points
$x_{ij}$,
for
$1 \leq i, j \leq 6$.
\end{Oiii}
\end{ass}\vskip-\medskipamount

The linear forms
$l_1, \ldots, l_6$
then go over into the affine-linear maps
$l'_1, \ldots, l'_6$,
\begin{eqnarray*}
l_a'(x,y) = l_a(x,y,1) &=& A_{a1} x + A_{a2} y + A_{a3} \\
&=:& \widetilde{l}_a(x,y) + A_{a3} \,.
\end{eqnarray*}

\subsubsection*{Particular\/
$1$-manifolds
in\/~$\Pb^2(\bbR)$}

The~starting point of our construction are compact, connected
\mbox{$1$-manifolds}
$\Gamma$,
embedded
into~$\Pb^2(\bbR)$,
of the behaviour indicated in the two figures~below. Observe~that, in either case, the
\mbox{$1$-manifold}
meets the branch locus
$V(l_1 \cdots l_6)$
only in its double~points.

\begin{figure}[H]
\centerline{
\begin{overpic}[scale=1.0]{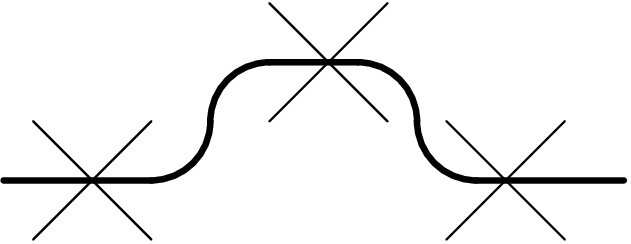}
 \put(-4,-9){$V(l_1)$}
 \put(65,-9){$V(l_2)$}
 \put(110,49){$V(l_3)$}
 \put(170,49){$V(l_4)$}
 \put(193,-9){$V(l_5)$}
 \put(262,-9){$V(l_6)$}
\end{overpic}}\vspace{.5cm}
\caption{A deformed line}
\label{Defo_line}\vspace{-.9cm}
\end{figure}

\begin{figure}[H]
\centerline{
\begin{overpic}[scale=1.0]{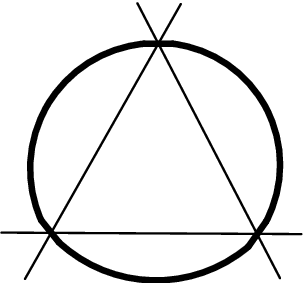}
 \put(-5,-9){$V(l_1)$}
 \put(125,-9){$V(l_2)$}
 \put(146,20){$V(l_3)$}
\end{overpic}}\vspace{.5cm}
\caption{A curve encircling a triangle}
\label{Curve_enc}\vspace{-.5cm}
\end{figure}

\begin{term}
We~call a
\mbox{$1$-manifold}
as in Figure~\ref{Defo_line} a {\em deformed line\/} and one of the type indicated in Figure~\ref{Curve_enc} a {\em curve encircling a triangle}.
We~allow as well compact
\mbox{$1$-manifolds}
that are of a similar shape as the one in Figure~\ref{Curve_enc}, but encircle a quadrangle or pentagon instead of a~triangle. Note~here that every compact, connected
\mbox{$1$-manifold}
is diffeomorphic
to~$S^1$.
\end{term}

\subsubsection*{Technical assumptions on the
\mbox{$1$-manifolds}}

Consider a single compact
\mbox{$1$-manifold}
in~$\Pb^2(\bbR)$
of one of the kinds described. We~may certainly assume that it is specified in parametrised~form.
In~the case of a deformed line, the affine part allows a parametrisation given by a
$C^1$-map
$$\gamma\colon \bbR \to \bbR^2 \subset \Pb^2(\bbR) \,.$$
Similarly, a curve encircling a polygon may be parametrised by a
$C^1$-map
$$\overline\gamma\colon \Pb^1(\bbR) \to \bbR^2 \subset \Pb^2(\bbR) \,.$$
In~this situation, we denote the restriction
$\overline\gamma |_\bbR$
by~$\gamma$.

For reasons of simplicity of the presentation, let us make the following assumptions on the
\mbox{$1$-manifolds}
we consider, as well as on their~parametrisations.

\begin{ass}[on the
$1$-manifolds]
\begin{Aiii}
\item
\label{A0}
In~the case that
$\Gamma$
is a curve encircling a polygon, the parametrisation is of the kind that
$\overline\gamma(\infty)$
is not a branch~point.
\item
\label{A1}
Suppose that the branch
point~$x_{ij}$
is met by the
\mbox{$1$-man}\-i\-fold~$\Gamma$,
and let
$x_{ij} = \gamma(t_{ij})$.
Then the parametrisation is {\em linear\/}
near~$t_{ij}$.
More~precisely,
$$\gamma(t) = (t-t_{ij})b_{ij} + x_{ij} \,,$$
for some
$b_{ij} \in \bbR^2$
and~$t \in [t_{ij}-1, t_{ij}+1]$.
We~require, moreover,~that
$$b_{ij} \not\in \bbR(A_{i2}, -A_{i1}) \cup \bbR(A_{j2}, -A_{j1}) \,.$$
I.e.,~that the direction
of~$\Gamma$
near
$x_{ij}$
differs from the directions of the lines
$V(l_i)$
and~$V(l_j)$
of the branch locus that meet
at~$x_{ij}$.
\item
\label{A2}
In~addition, for all branch points that are met
by~$\Gamma$,
the intervals
$[t_{ij}-1, t_{ij}+1]$
are mutually disjoint.
\item
\label{A3}
Finally,~for a deformed line, we assume that the map
$\gamma$
is of the shape %the parametrisation is such that
\begin{equation}
\label{def_line_formula}
\gamma(t) = tb + \gamma_0(t)
\end{equation}
for all
$t \in \bbR$,
where
$b \in \bbR^2 \!\setminus\! \{(0,0)\}$
is a certain vector and the
\mbox{$C^1$-map}
$\gamma_0\colon \bbR \to \bbR^2$
has compact~support. We~require, moreover,~that
$$b \not\in \bbR(A_{12}, -A_{11}) \cup \cdots \cup \bbR(A_{62}, -A_{61}) \,.$$
I.e.,~that the direction
of~$\Gamma$
near infinity differs from the directions of the lines
$V(l_1), \ldots, V(l_6)$
forming the branch locus.
\end{Aiii}
\end{ass}

\begin{nota}
\label{ramif_pts}
\begin{iii}
\item
We~let
$x_{i_1,j_1}, \ldots, x_{i_n,j_n} \in V(l_1 \cdots l_6) \subset \Pb^2(\bbR)$,
for
$n = 3$,
$4$,
or~$5$,
be the points of the branch locus that are met
by~$\Gamma$.
\item
Moreover,~we denote by
$t_1 = t_{i_1,j_1}, \ldots, t_n = t_{i_n,j_n} \in \bbR$
the parameters of these points. I.e.\ the reals such that
$\gamma(t_1) = x_{i_1,j_1}$,
\ldots,
$\gamma(t_n) = x_{i_n,j_n}$.
\end{iii}\smallskip
\end{nota}

\subsubsection*{The main construction---$2$-tori
from
$1$-manifolds}

\begin{nota}
Write~$\Tb := \Pb^1(\bbR) \times \Pb^1(\bbR) \supset \bbR \times \bbR = \bbR^2$,
for the obvious~compactification. This~is our standard
\mbox{$2$-torus}.
\end{nota}

\begin{defi}
By~a {\em liftable torus,} in this article, we mean a torus
$\alpha_\bullet \colon \Tb \to X'$
that lifts through the blowing-down map
$\bl\colon X \to X'$
between the
$K3$~surface
$X$
and the singular double cover
$X'$
of~$\Pb^2$.
I.e., we require the existence of a continuous map
$\alpha\colon \Tb \to X$
to the
$K3$~surface,
such that
$\bl \circ \alpha = \alpha_\bullet$.
\end{defi}

\begin{mconstr}[\mbox{$2$-dimensional}
tori from
\mbox{$1$-manifolds}]
\label{constr_main}
\leavevmode
\begin{abc}
\item
There~are two liftable tori, i.e.\ continuous maps
$$\alpha_\Gamma, {\widetilde\alpha}_\Gamma\colon \Tb \longrightarrow X' \,,$$
which differ only by the
involution~$\zeta$,
associated with each deformed
line~$\Gamma$,
as~above.
\item
Similarly,~there are two liftable tori
$$\alpha_{\Gamma,b}, {\widetilde\alpha}_{\Gamma,b}\colon \Tb \longrightarrow X' \,,$$
differing only by the
involution~$\zeta$,
associated with each
curve~$\Gamma$
encircling a polygon, together with a~vector
\begin{equation}
\label{b_generic}
b \in \bbR^2 \setminus \big( \bbR(A_{12}, -A_{11}) \cup \cdots \cup \bbR(A_{62}, -A_{61}) \big) \,.
\end{equation}
\end{abc}
\end{mconstr}

\begin{ttt}
\label{const_main_detail}
The tori are obtained by the following purely topological construction.

\begin{iii}
\item
(Adding an imaginary direction)
Extend the map
$\gamma\colon \bbR \to \bbR^2 \subset \Pb^2(\bbR)$
to a continuous map in two variables by~putting
\begin{equation}
\label{constr_Stephan}
\gamma'\colon \bbR^2 \longrightarrow \bbC^2 \subset \Pb^2(\bbC), \qquad (t,u) \mapsto \gamma(t) + \mi ub \,.
\end{equation}
Here,~in the case of a deformed line,
the vector~$b$
is taken from~A.\ref{A3}, while, in the case of a curve encircling a polygon,
$b$~is
part of the input~data.
The~map
$\gamma'$
allows a continuous prolongation
$$\alpha'\colon \Tb \to \Pb^2(\bbC) \,.$$

\item
(Lifting to the double cover)
The map
$\alpha'$
allows a continuous~lift
$$\alpha''\colon \Tb \to X'$$
to the double~cover. There~is, in fact, a second such lift,
$\widetilde\alpha''$,
differing
from~$\alpha''$
by the
involution~$\zeta$.
\item
(Widening the holes)
Use~a self-map
$\Psi\colon \Tb \to \Tb$
that induces a homeomorphism
$$\Psi |_{\Tb \setminus \big( \overline{U_\frac12(t_1,0)} \cup \ldots \cup \overline{U_\frac12(t_n,0)} \big)} \colon \Tb \setminus \big( \overline{U_\frac12(t_1,0)} \cup \ldots \cup \overline{U_\frac12(t_n,0)} \big) \stackrel{\cong}{\longrightarrow} \Tb \setminus \{(t_1,0), \ldots, (t_n,0)\}$$
and sends
$\smash{\overline{U_\frac12(t_i,0)}}$
constantly to
$(t_i,0)$,
for~$i=1,\ldots,n$,
to define
$$\alpha''' := \alpha'' \!\circ\! \Psi\colon \Tb \to X' \,.$$
This~is the liftable torus
$\alpha_\Gamma$
or
$\alpha_{\Gamma,b}$~desired.
\end{iii}
\end{ttt}

\begin{rems}
\begin{iii}
\item
The cohomology class
$c_\alpha \in H^2(X, \bbZ)$
is dependent on the lift
$\alpha\colon \Tb \to X$
chosen. Two~different lifts lead to cohomology classes differing by a summand
from~$P$.
I.e.,~the image
$\overline{c}_\alpha \in H^2(X, \bbZ)/P$
depends only on the torus
$\alpha_\Gamma$
or
$\alpha_{\Gamma,b} \colon \Tb \to X'$
lifted. In~fact,
$\overline{c}_\alpha$
depends only on the homotopy class of
$\alpha_\Gamma$
or
$\alpha_{\Gamma,b}$,~respectively. This~is a direct consequence of Corollary~\ref{lift_coh_modP}.

Moreover,~by Lemma~\ref{eigenv_invol}.b), the classes
$\overline{c}_\alpha$
and
$\overline{c}_{\widetilde\alpha} \in H^2(X, \bbZ)/P$
only differ by~sign.
\item
Step~iii) in the construction above could be eliminated at the cost of a somewhat less elegant description of the class
$\overline{c}_\alpha \in H^2(X, \bbZ)/P$,
which works as follows.

Take~the homology class
$\alpha''_*(z_\Tb) \in H_2(X',\bbZ)$
and put
$\overline{c}_\alpha \in H^2(X, \bbZ)/P$
to be its image, according to the homomorphism described in Lemma~\ref{h_doubcov}.b).
\item
The~two theorems below are the principal results on the construction just described. As~their proofs are a bit technical, we postpone them to the final~section. The~same applies to the claims made in~\ref{const_main_detail}. Cf.,~in particular, Propositions~\ref{lift} and \ref{ext_boundary}.
\end{iii}
\end{rems}

\begin{theo}[Independence
of~$b$
in the case of a curve encircling a polygon]
\label{indep_b}
Let\/~$\Gamma$
be a curve encircling a polygon~and
$$\underline{b}, \underline{\underline{b}} \in \bbR^2 \setminus \big( \bbR(A_{12}, -A_{11}) \cup \ldots \cup \bbR(A_{62}, -A_{61}) \big)$$
any two~vectors. Moreover,~let\/
$\smash{\underline\alpha\colon \Tb \to X}$
be any continuous lift
of\/~$\smash{\alpha_{\Gamma,\underline{b}}\colon \Tb \to X'}$
and\/
$\smash{\underline{\underline\alpha}\colon \Tb \to X}$
any continuous lift
of\/~$\smash{\alpha_{\Gamma,\underline{\underline{b}}}\colon \Tb \to X'}$.\smallskip
%and
%$\smash{\widetilde{\underline{\underline\alpha}}\colon \Tb \to X}$
%any continuous lift
%of\/~$\smash{\alpha_{\Gamma,\widetilde{\underline{\underline{b}}}}\colon \Tb \to X'}$.\smallskip

\noindent
Then\/~$\smash{\overline{c}_{\underline\alpha} \in H^2(X, \bbZ)/P}$
coincides either with\/
$\smash{\overline{c}_{\underline{\underline\alpha}}}$
or with\/
$\smash{-\overline{c}_{\underline{\underline\alpha}} \in H^2(X, \bbZ)/P}$.
\end{theo}

\begin{theo}[Periods as improper integrals]
\label{periods_ints}
Let\/~$X$
be the\/
$K3$~surface
obtained as the minimal desingularisation of the double cover\/
$X'\colon w^2 = l_1\cdots l_6$,
for\/
$l_1, \ldots, l_6$
real linear~forms. Assume that no three of the linear forms have a zero in~common.
On~$X$,
fix the global holomorphic\/
\mbox{$(2,0)$-form}~$\omega$,
determined by formula\/~(\ref{global_form}), below.\smallskip

\noindent
Let,~moreover,
$\Gamma$
be a compact\/
\mbox{$1$-manifold}
in\/~$\Pb^2(\bbR)$
being a deformed line or a curve encircling a~polygon. We~assume that\/
$\Gamma$
is parametrised by\/
$\gamma\colon\bbR \to \bbR^2 \subset \Pb^2(\bbR)$
and that Assumptions~A.i) to~A.iv) are~fulfilled.

\begin{iii}
\item
Let\/~$\Gamma$
be a deformed line. Then, for every continuous lift\/
$\alpha\colon\Tb \to X$
of the torus\/
$\alpha_\Gamma\colon\Tb \to X'$,
\begin{equation}
\label{int_def_line}
(c_\alpha, [\omega]) = \mi \int\limits_{-\infty}^\infty \! \int\limits_{-\infty}^\infty \!\frac{\dot{\gamma_1}(t) b_2 - \dot{\gamma_2}(t) b_1}{\sqrt{l_1\big(\gamma(t)+\mi ub\big) \cdots l_6\big(\gamma(t)+\mi ub\big)}} \,du\,dt \,,
\end{equation}
where the vector\/
$b \in \bbR^2 \setminus \bbR(A_{12}, -A_{11}) \cup \ldots \cup \bbR(A_{62}, -A_{61})$
is taken from~A.iv).
\item
Let\/~$b \in \bbR^2 \setminus \bbR(A_{12}, -A_{11}) \cup \ldots \cup \bbR(A_{62}, -A_{61})$
be any vector and\/
$\Gamma$
a curve encircling a~polygon.
%any vector not in\/
%$\bbR(A_{12}, -A_{11}) \cup \ldots \cup \bbR(A_{62}, -A_{61})$.
Then,~for every continuous lift\/
$\alpha\colon\Tb \to X$
of the torus\/
$\alpha_{\Gamma,b}\colon\Tb \to X'$,
\begin{equation}
\label{int_curve_enc}
(c_\alpha, [\omega]) = \mi \int\limits_{-\infty}^\infty \! \int\limits_{-\infty}^\infty \!\frac{\dot{\gamma_1}(t) b_2 - \dot{\gamma_2}(t) b_1}{\sqrt{l_1\big(\gamma(t)+\mi ub\big) \cdots l_6\big(\gamma(t)+\mi ub\big)}} \,du\,dt \,.
\end{equation}
\end{iii}
In~either case, the square root in the integrand is meant to be the same as that in the
form\/~$\omega$.
\end{theo}

\begin{rems}
\begin{iii}
\item
In the case of a deformed line, the outer integral in~(\ref{int_def_line}) is actually~proper. Indeed,~formulae (\ref{constr_Stephan}) and~(\ref{def_line_formula}) together imply that the
map~$\alpha'$,
and hence the
lift~$\alpha''$,
too, is holomorphic near every point
$(t,u) \in \Tb$
with~$t \not\in \supp \gamma_0$.
Furthermore,~since
$\omega$~is
a
\mbox{$(2,0)$-form}
on~$X$
and
$\Tb$
is of real dimension two, holomorphicity enforces the pull-back
$(\alpha'' |_{\bbR^2 \setminus \{\frakt_1,\ldots,\frakt_n\}} ){}^* \omega'$
to be the null form in a neighbourhood
of~$(t,u)$.
Finally,~the integrand in~(\ref{int_def_line}) exactly corresponds to this differential~form. Thus,~the outer integral actually ranges only
over~$\supp \gamma_0$. 
\item
On~the other hand, a curve encircling a polygon is a curve
in~$\bbR^2$
of finite~length. It~is convenient to reparametrise it by a compact interval, so that in this case the outer integral in~(\ref{int_curve_enc}) becomes proper,~too.
\end{iii}
\end{rems}

\subsubsection*{The global holomorphic differential form}

\begin{prop}
Let\/~$X$
be the minimal desingularisation of the double cover\/
$X'\colon w^2 = l_1\cdots l_6$,
for\/
$l_1, \ldots, l_6$
linear forms in the variables\/
$x$,~$y$,
and\/~$z$.
Assume that no three of the six linear forms have a zero in~common.\smallskip

\noindent
Then,~for any linear
form\/~$l$
that is not just a scalar multiple of\/
$x$
or\/~$y$
and defines a nonsingular curve
on\/~$X'$,
\begin{equation}
\label{global_form}
\omega' := \frac{d(\frac{x}l) \wedge d(\frac{y}l)}{\frac{w}{l^3}}\
\end{equation}
is a differential form
on\/~$X'$,
whose pull-back\/
$\omega$
to\/~$X$
is a global holomorphic\/
$(2,0)$-form
without zeroes or~poles.\smallskip

\noindent
{\bf Proof.}
{\em
On~the affine plane, given by
$l \neq 0$,
the functions
$\smash{\xi := \frac{x}l}$
and
$\smash{\eta := \frac{y}l}$
form a system of coordinates. In~particular, there are functions
$\lambda_1, \ldots, \lambda_6$
of total
degree~$1$
such that
$\smash{\lambda_i(\xi, \eta) = l_i(\frac{x}l, \frac{y}l, \frac{z}l)}$,
for~$i = 1, \ldots, 6$.
The~corresponding affine part of the double
cover~$X'$
is thus given by
$\psi^2 = \lambda_1 \cdots \lambda_6$,
for~$\smash{\psi := \frac{w}{l^3}}$.

Hence,~$\omega' = d\xi \wedge d\eta / \psi$.
This~immediately shows that, outside
$V(l_1 \cdots l_6)$,
which is the branch locus, and
$V(l)$,
which is not on our chart,
$\omega$
has neither zeroes, nor~poles.

Next,~consider a point that lies on exactly one of the lines
$V(\lambda_1), \ldots, V(\lambda_6)$.
Then~at least one of the partial derivatives
$\smash{\frac{\partial{(\lambda_1 \cdots \lambda_6)}}{\partial\xi}}$
and
$\smash{\frac{\partial{(\lambda_1 \cdots \lambda_6)}}{\partial\eta}}$
is nonzero. Without restriction, assume that
$\smash{\frac{\partial{(\lambda_1 \cdots \lambda_6)}}{\partial\eta} \neq 0}$.
Then~$\psi$
and~$\xi$
form a local system of coordinates. Moreover,~the equation
$\psi^2 = \lambda_1 \cdots \lambda_6$
yields
$\smash{2\psi\,d\psi = \frac{\partial{(\lambda_1 \cdots \lambda_6)}}{\partial\xi} d\xi + \frac{\partial{(\lambda_1 \cdots \lambda_6)}}{\partial\eta} d\eta}$,
which shows that
$\smash{2\psi\,d\xi \wedge d\psi = \frac{\partial{(\lambda_1 \cdots \lambda_6)}}{\partial\eta} d\xi \wedge d\eta}$.
Thus,~$d\xi \wedge d\eta$
is indeed divisible
by~$\psi$
and the quotient
$\omega' = d\xi \wedge d\eta / \psi$
is nonzero near the point~considered.

Finally,~let us write
$W := \div \omega \in \Div X$.
The~calculations above prove that
$W$
has the~form
$$W = n [L] + \!\!\sum_{1\le i<j\le 6}\!\!\!\! n_{ij} [E_{ij}] \,,$$
for~$L := V(l)$.
We~have to show that, in fact, all the coefficients~vanish.

For~this, let us note that
$W$~is,
by construction, a canonical~divisor. Thus,~the adjunction formula shows that
$\smash{(W+[E_{ij}])[E_{ij}] = 2g_{E_{ij}}-2 = -2}$,
i.e.\
$\smash{W[E_{ij}] = 0}$,
for~$1 \leq i < j \leq 6$.
Consequently,
$\smash{n_{ij} = -\frac12 W[E_{ij}] = 0}$.
Moreover,~$L$
is a double cover of a projective line ramified over six points, and hence its genus
is~$g_L = 2$.
Here,~adjunction shows that
$(W+[L])[L] = 2g_L-2 = 2$.
As~$[L][L] = 2$,
this yields
$W[L] = 0$
and
$\smash{n = \frac12 W[L] = 0}$.
The assertion~follows.
}
\eop
\end{prop}

\section{Explicit description of transcendental cohomology classes---\\Computation of periods}

\subsubsection*{The main algorithms}

\begin{ttt}
For clarity of the presentation, let us present our main algorithms in a somewhat idealised setting.
Imagine for a moment that one could calculate
\mbox{$2$-dimensional}
integrals {\em exactly,} as complex numbers. Then, given an arbitrary
$K3$~surface
$X = X_{(a_0,b_0,c_0,d_0)}$
(cf.~formula~(\ref{eq_norm})),
we could run the following~algorithm.
\end{ttt}

\begin{algo}[Explicit description of the transcendental part of
$H^2(X, \bbZ)$
--- Idealised setting]%
\leavevmode
\label{trans_part}

\begin{iii}
\item
Set up a list of
\mbox{$1$-manifolds}
in~$\Pb^2(\bbR)$
consisting of

\begin{itemize}
\item[$\bullet$ ]
for each triangle, quadrangle, or pentagon formed in the affine plane given by
``$z=1$''
by the lines of the branch locus, one curve encircling it,
\item[$\bullet$ ]
through each triple of double points of the branch locus that may be connected by a deformed line within the affine plane given by
``$z=1$'',
one deformed line.
\end{itemize}

\item
For every
\mbox{$1$-manifold}
listed, construct a torus
$\alpha\colon \Tb \to X$,
as described in~\ref{constr_main}.
\item
\label{item_iii}
For every torus
$\alpha$
constructed, calculate
$(c_\alpha, [\omega]) \in \bbC$
using Theorem~\ref{periods_ints}, formulae (\ref{int_def_line}) and~(\ref{int_curve_enc}).
\item
\label{item_iv}
Detect all
\mbox{$\bbZ$-linear}
dependencies among these numbers, and select six (linear combinations of) tori
$\alpha_1, \ldots, \alpha_6$
that yield a
\mbox{$\bbZ$-basis}
for the complex numbers obtained.

In~the case that the numbers
$(c_\alpha, [\omega])$
form a
\mbox{$\bbQ$-vector}
space of dimension less than six, do the calculations on a surface
$X_{(a,b,c,d)}$,
for a quadruple
$(a,b,c,d)$
near
$(a_0,b_0,c_0,d_0)$.
\item
\label{item_v}
Detect the
\mbox{$\bbQ$-linear}
dependency among the products
$\smash{(c_{\alpha_i}, [\omega]) \cdot (c_{\alpha_j}, [\omega])}$,
for
$1 \leq i \leq j \leq 6$.
I.e., determine the projective quadric
$Q_{\kappa,16} \subset \Pb^5(\bbC)$
in the restricted period~space.

Again, do this on a neighbouring surface that is generic, if~necessary.
\item
Determine the quadratic form
$\kappa$
exactly, not only up to scaling, from the fact that the self-product is known for a deformed~line. (Cf.\ Lemma~\ref{scal}.)
\item
Verify that
$c_{\alpha_1},\ldots,c_{\alpha_6}$
indeed form a generating system of
$H^2(X,\bbZ)/P$,
not only up to finite~index. (Cf.\ Theorem~\ref{tori_generate}.)
\end{iii}
\end{algo}

\begin{rems}
\begin{iii}
\item
By~Theorem~\ref{mark_rel}, the result that
$c_{\alpha_1},\ldots,c_{\alpha_6}$
generate
$H^2(X,\bbZ)/P$
implies that these elements give rise to a class of markings
on~$X$,
for which the numbers
$(c_{\alpha_i}, [\omega])$
are indeed the coordinates of the period~vector.
\item
In practice, of course, in step~\ref{item_iii}, one has to use numerical integration methods working at a finite precision. Moreover, instead of the exact linear algebra in step~\ref{item_iv}, the singular value decomposition is playing the main~role.
\item
Algorithm~\ref{cup_prod} below gives a more thorough description of step~\ref{item_v}.
\end{iii}
\end{rems}

\begin{ttt}
Furthermore,~given a
$K3$~surface
$X = X_{(a_0,b_0,c_0,d_0)}$
having real multiplication (cf.\ Section \ref{RMCM_Hodge}), one could run the algorithm~below.
\end{ttt}

\begin{algo}[Tracing the modular curve --- Idealised setting]%
\leavevmode
\label{trac_fict}

\begin{iii}
\item
Detect the linear relations between the coordinates of the period points encoding real multiplication. Cf.\ Theorem~\ref{RMCM_dim} for~details.
\item
Using~Newton iteration, find surfaces
$X_{(a,b,c,d)}$,
for quadruples
$(a,b,c,d)$
near
$(a_0,b_0,c_0,d_0)$,
the coordinates of the period points of which fulfil the same relation, not exactly but at high precision.
\end{iii}
\end{algo}

\begin{rem}
If~one had enough of these quadruples (exactly) then it would just be linear algebra to find an algebraic curve through these quadruples. Again, the singular value decomposition plays the main role in our adaption of Algorithm~\ref{trac_fict} to~practice.
\end{rem}

\subsubsection*{Numerical integration}

\begin{strat}
In~order to numerically calculate the
\mbox{$2$-d}i\-men\-sional
integrals occurring in Theorem~\ref{periods_ints}, our strategy is roughly as~follows.

\begin{iii}
\item
We decompose the domain of the outer integral into finitely many intervals
$[a_{k-1}, t_k]$
and
$[t_k, a_k]$,
each of which contains none of the
$t_1, \ldots, t_n$
as an inner point and exactly one as an~endpoint. For~integration, each such interval is treated individually, the results being added together at the very~end.
\item
For~each concrete interval, the inner integral needs to be computed only over
$[0,\infty) \subset \bbR$.
Indeed,~the sign change
$u \mapsto (-u)$
transforms the integrand into its complex conjugate or minus its complex conjugate, depending on the sign of
$(l_1 \cdots l_6) \!\circ\! \gamma$
on the~interval.
\item
We~decompose the two-dimensional domains of integration emerging into three subdomains each, as indicated in Figure~\ref{integral}~below.
For~each subdomain, we determine the integral using Fubini's Theorem in a rather naive manner. The~inner integrals, for the computation, are taken along the line segments and rays, as shown in the~figure.
\item
To~compute the inner integrals of the integral over the triangle, we use the Gau{\ss}-Legendre~method.

For~the computation of the inner integrals that are improper, we decompose the ray into finitely many intervals
$I_1,\ldots,I_N$
of lengths increasing by a factor of
$\smash{\frac{|I_{i+1}|}{|I_i|} = q = \frac32}$
and the remaining~ray. The~substitution
$\smash{u' := 1/u}$
transforms the integral over the latter into a definite~one. We~then use the Gau{\ss}-Legendre method for each of the individual integrals and add the results~together.

\begin{figure}[H]
\centerline{
\raisebox{.7mm}{
\begin{overpic}[scale=0.6]{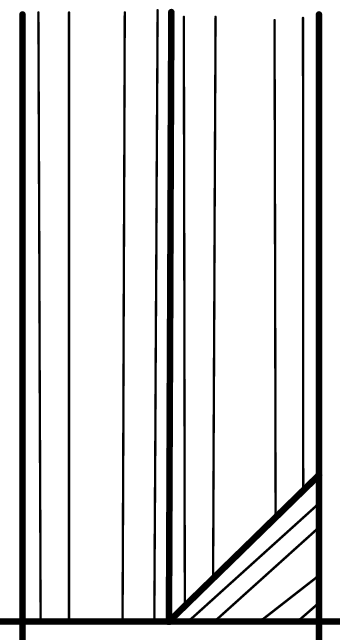}
 \put(0,-9){$a_{k-1}$}
 \put(88,-9){$t_k$}
 \put(-13,173){$u \!\!\uparrow$}
 \put(101,2){$t \!\!\rightarrow$}
\end{overpic}}\hspace{3cm}

\begin{overpic}[scale=0.6]{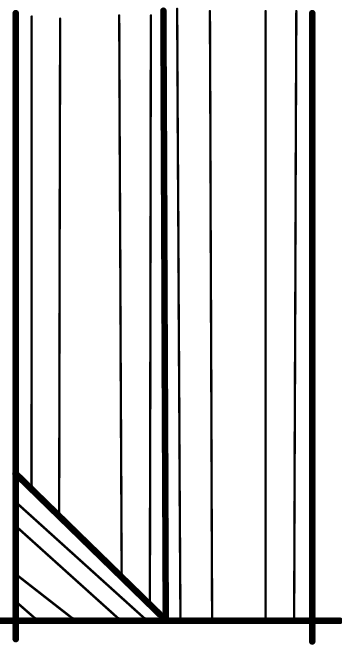}
 \put(0,-9){$t_k$}
 \put(88,-9){$a_k$}
 \put(-13,175){$u \!\!\uparrow$}
 \put(101,4){$t \!\!\rightarrow$}
\end{overpic}
}\vspace{.3cm}
\caption{The two-dimensional domains of integration}
\label{integral}\vspace{-.8cm}
\end{figure}

Thereby,~we work with adaptive stepsizes, as~follows. First,~we parametrise the~ray. This~results in a hyperelliptic integral of the~form
$$\int_{u_0}^\infty \!\!\frac1{\sqrt{(a_1 i u + b_1) \cdots (a_6 i u + 
b_6)}} \,du \,,$$
for
$u_0 \geq 0$
and~$a_1,\ldots,a_6,b_1,\ldots,b_6 \in \bbR$.
Note~that all singularities of the integrand are located on the imaginary line, while the domain of integration is a ray contained in the~reals.

We~compute the smallest and the largest among the absolute values of all the singularities of the~integrand. Denote~these by
$S$
and~$M$.
The~first interval
$I_1$,
starting
at~$u_0$,
is then taken to be of length
$\max(u_0, S)$.
We~then increase the lengths of the intervals by the factor of
$\smash{q = \frac32}$
in each step, until we reach a value of
$u_1 > 4M$.
Then,~on the remaining ray, coordinates are changed as described~above.

The aim of this strategy is to make sure that, for each interval of integration, the ratio between its length and the distance to the next singularity is bounded by a constant. This~suffices to make sure that the Gauss-Legendre method provides the expected~precision.

For~example, applying the change of coordinates, the final ray is turned into the interval
$\smash{[0,\frac{1}{u_1}]}$,
whereas the smallest absolute value of a singularity becomes
$\smash{\frac{1}{M}}$.
As~$\smash{\frac1{u_1} < \frac1{4M}}$,
the resulting integral can be computed numerically without any further 
splitting of the
interval~$\smash{[0,\frac{1}{u_1}]}$.
\item
For the outer integrals, we again use the Gau{\ss}-Legendre method.
\end{iii}
\end{strat}

\begin{rems}
\begin{iii}
\item
One~might want to choose the degrees for the Gau{\ss}-Legendre method differently for the two levels, the computations of the inner and the outer~integrals. However, we never did~so.
\item
Our~experience is that the approximations obtained start to get useful for
degree~$30$
[i.e.\ order~$60$],
and that the calculations start to become slow at
degree~$100$,
while
degree~$300$
appears to be some kind of acceptance limit for usability in~practice.
\item
To illustrate this, let us present some data related to Example~\ref{qw13}, below. Working on a suitably chosen affine chart, there are four deformed lines,
four curves encircling a triangle, and six curves encircling a quadrangle occurring. The running times to compute the corresponding 14 integrals at various precisions, as well as the largest relative errors in comparison with degree 200, are summarised in the table~below. We~did the numerical calculations in floats of as many decimal digits as the degree~chosen.

\begin{table}[H]
\begin{center}
\begin{tabular}{|c||c|c|}
\hline
Degree & Running time & Largest relative error \\\hline\hline
\phantom{0}20 & \phantom{00}20.99\,s & $3.301361055\cdot10^{-18}$ \\\hline
\phantom{0}30 & \phantom{00}45.95\,s & $6.544716595\cdot10^{-26}$ \\\hline
\phantom{0}50 & \phantom{0}138.58\,s & $2.767253174\cdot10^{-40}$ \\\hline
100           & \phantom{0}640.45\,s & $2.562917026\cdot10^{-77}$ \\\hline
150           &           1644.72\,s & $3.312576296\cdot10^{-114}$ \\\hline
200           &           3234.76\,s & --- \\\hline
\end{tabular}
\end{center}\vskip0mm
\caption{Running times to compute the 14 integrals}\vskip-5mm
\end{table}
\end{iii}
\end{rems}

\subsubsection*{Transcendental\/
$2$-cocycles}
Given a concrete
$K3$~surface
of the type aforementioned, we usually apply the main construction~\ref{constr_main} to all curves encircling a polygon and all deformed lines occurring on a fixed affine~chart. This pro\-vides by far more than six representatives of classes
in~$H^2(X, \bbZ) / P$.

\begin{algo}[Determining the cup product
on~$P^\perp \!\otimes_\bbZ\! \bbQ$,
up to scaling]%
\label{cup_prod}%
\leavevmode

\noindent
Let~$X = X_{(a_0,b_0,c_0,d_0)}$
be a
$K3$~surface
given as the minimal desingularisation of a double cover of the form~(\ref{eq_norm}), for
$a_0, b_0, c_0, d_0 \in \bbR$.
Moreover,~let
$\alpha_1, \ldots, \alpha_n\colon \Tb \to X$
be tori of the kind that the classes
$\overline{c}_{\alpha_1}, \ldots, \overline{c}_{\alpha_n} \in H^2(X, \bbZ)/P$
form a system of generators.

\begin{iii}
\item
Choose~open neighbourhoods
$\bbD \cong U(a_0) \ni a_0$,
\ldots,
$\bbD \cong U(d_0) \ni d_0$
in such a way that, for every
$(a,b,c,d) \in U := U(a_0) \times \cdots \times U(d_0)$,
no three of the resulting six lines
in~$\Pb^2_\bbC$
have a point in~common. Then~the cohomology classes
$c_{\alpha_1}, \ldots, c_{\alpha_n}$
uniquely extend to the whole family
$\frakX \to U$
of the
$K3$~surfaces
$X_{(a,b,c,d)}$.

Moreover, choose
$N$
surfaces
$X_1, \ldots, X_N$
at random from the family and write down the corresponding holomorphic
\mbox{$2$-forms}
$\omega_1, \ldots, \omega_N$.
(We~usually work with
$N=50$.)
\item\label{sing_val}
Set up the matrix
$M := (c_{\alpha_j}, [\omega_i])_{1\leq i\leq N, 1\leq j\leq n}$
using numerical integration and calculate the singular value decomposition
of~$M$.
Six singular values should be numerically nonzero. The others give rise to linear relations among the cohomology classes
$c_{\alpha_1}, \ldots, c_{\alpha_n}$.
\item
Choose a basis
$c_1, \ldots, c_6$
of the free
\mbox{$\bbZ$-module}
spanned by
$c_{\alpha_1}, \ldots, c_{\alpha_n}$
modulo the relations~found.
\item
Build from
$M$
the
\mbox{$N \!\times\! 15$-matrix}
$F := (c_{j_1}, [\omega_i]) (c_{j_2}, [\omega_i])_{1\leq i\leq N, 1\leq j_1 \leq j_2 \leq 6}$.
Then determine an approximate solution of the corresponding homogeneous linear system of equations, using the QR-factorisation
of~$F$.
I.e., detect the one linear relation between the 15 products that is approximately fulfilled for
all~$i$.
This~solution vector describes the symmetric, bilinear pairing~desired.
\end{iii}
\end{algo}

\begin{rems}
\begin{iii}
\item
As~it relies only on the projective quadric
$Q_{\kappa,16} \subset \Pb^5(\bbC)$
defined by the symmetric, bilinear 
form~$\kappa$
(cf.~Corollary~\ref{periods_quadric}, in particular formula~(\ref{restr_per})) and not on
$\kappa$~itself,
our method is inherently limited to determining the cup product pairing,
or~$\kappa$,
up to~scaling.
\item
Algorithm~\ref{cup_prod} works with the restricted period space. Therefore, it detects only the restriction of the cup product pairing
to~$P^\perp \!\otimes_\bbZ\! \bbQ \subset H^2(X, \bbQ)$.
\item
Every~class in
$H^2(X, \bbZ) / P$
has a unique representative
in~$P^\perp \!\otimes_\bbZ\! \bbQ$,
the orthogonal~projection
$$\pr\colon H^2(X, \bbZ) / P \longrightarrow P^\perp \!\otimes_\bbZ\! \bbQ$$
being~injective. (Note that
$P \subset H^2(X, \bbZ)$
is a sublattice not of discriminant
$\pm1$,
so the image
of~$\pr$
indeed contains non-integral~classes.) Thus,~we actually compute the pairing
$(\pr(\cdot), \pr(\cdot))\colon H^2(X, \bbZ) \times H^2(X, \bbZ) \to \bbQ$,
up to~scaling.
\item
The algorithm yields the matrix describing the cup product pairing, up to scaling, with respect to
$(c^1, \ldots, c^6)$,
but provides the dual basis
$(c_1, \ldots, c_6)$.
The~corresponding matrices are, in fact, inverse to each~other. Indeed,
$$\smash{\sum_k} (u, c_k) (c^k, v) = \Big( u, \smash{\sum_k} (c^k, v) c_k \Big) = (u, v)$$
and therefore
$\sum_k (c_i, c_k) (c^k, c^j) = \delta_{ij}$.
\end{iii}
\end{rems}

\begin{rem}
Often, one may detect the cup product pairing (up to scaling) relying only on a single surface (instead
of~$N = 50$),
using an {\tt LLL}-based approach. Indeed,~it seems plausible that the cohomology classes of the tori constructed are rather small within the transcendental lattice. Thus, between their periods, being irrational numbers known at high precision, there should be no further small quadratic relations besides that coming from the cup product~pairing. 
\end{rem}

\begin{lem}%[Scaling]
\label{scal}
Let\/~$\alpha\colon \Tb \to X$
be a torus constructed from a deformed line, as in Figure~\ref{Defo_line}.
Then\/~$(\pr(c_\alpha), \pr(c_\alpha)) = -1$.\smallskip

\noindent
{\bf Proof.}
{\em
The expression
$(\pr(c_\alpha), \pr(c_\alpha))$
is purely cohomological, and therefore homotopy invariant. Up~to homotopy, the three double points may be assumed to be~collinear. Without restriction, the deformed
line~$\Gamma$
is then actually a line, parametrised by
$\gamma(t) = tb$,
i.e.\
$\gamma_0(t) \equiv 0$.

In~this situation, the main construction~\ref{constr_main} yields~that
$$\alpha' |_{\bbR^2} = \gamma' \colon \bbR^2 \longrightarrow \bbC^2 \subset \Pb^2(\bbC)$$
is given by
$(t,u) \mapsto (t + \mi u)b$.
Thus,~$\alpha' |_{\bbR^2}$
is a holomorphic embedding onto an affine line contained
in~$\Pb^2(\bbC)$.
In~other words, the torus
$\alpha'$,
as described in (\ref{constr_Stephan}) and~(\ref{proj_ext}), parametrises, generically one-to-one, the complex projective line
$g' \subset \Pb^2(\bbC)$
through the three double~points.

The inverse image
of~$g'$
in~$X'$
splits into two complex algebraic curves both being rational as the branch locus consists entirely of double~points. Thus,~the lift
$\alpha''\colon \Tb \to X'$
parametrises, generically one-to-one, one such
component~$g''$.
And
$\alpha_\Gamma\colon \Tb \to X$
parametrises the strict transform
$g \subset X$
of~$g''$.
Consequently,
$c_\alpha = [g]$.

One~has
$[g]^2 = -2$,
$[g] \!\cdot\! \pi^*[l] = 1$,
and
$[g][E_1] = [g][E_2] = [g][E_3] = 1$,
for
$E_1$,
$E_2$,
and~$E_3$
the exceptional points over the double points met
by~$\Gamma$.
Clearly,~$[g]$
is perpendicular to the classes of the three other exceptional~curves. As~for algebraic classes, the cup product pairing coincides with the intersection pairing \cite[\S0.4]{GH}, using the facts that
$[E_i]^2 = -2$
and
$\pi^*[l] \!\cdot\! \pi^*[l] = 2$,
one finds that
$\pr(c_\alpha) = D$,
for
$D$
the
\mbox{$\bbQ$-divisor}
$$\textstyle D := [g] - \frac12 \pi^*[l] + \frac12 [E_1] + \frac12 [E_2] + \frac12 [E_3] \,.$$
Finally,
$D^2 = -1$
is easily seen by a direct~calculation.
}
\eop
\end{lem}

\begin{rem}
We~use the observation made in Lemma~\ref{scal} for~scaling.
\end{rem}

\begin{theo}
\label{tori_generate}
Let\/~$X$
be the minimal desingularisation of a double cover
of\/~$\Pb^2_\bbC$,
ramified over a union of six real lines, such that no three of them have a point in~common. Then~the classes of the tori, as described in Proposition~\ref{lift} below, always generate the whole
of\/~$H^2(X, \bbZ) / P$.\smallskip

\noindent
{\em
{\bf Proof} (depending on numerical integration){\bf.}
There are four essentially distinct configurations of six real lines
in~$\Pb^2$,
no three of which have a point in common~\cite[Chapter~VII, Section 8.1]{Yo}. We~made an experiment for each~case.
}%
\eop
\end{theo}

\begin{rem}
Our approach to find example configurations for each of the four types was invariant-theoretic and used only a small proportion of the theory developed in~\cite{Yo}. Some~details are as~follows.

The~admissible configurations of six real~lines correspond under projective duality to arrangements of six real points
$(x_1:y_1:z_1), \ldots, (x_6:y_6:z_6) \in \Pb^2$,
no three of which are~collinear. I.e.,~such that each of the twenty~determinants
$$\det
\left(
\begin{array}{ccc}
x_{i_1} & x_{i_2} & x_{i_3} \\
y_{i_1} & y_{i_2} & y_{i_3} \\
z_{i_1} & z_{i_2} & z_{i_3}
\end{array}
\right) ,
$$
for
$1 \leq i_1 < i_2 < i_3 \leq 6$,
is nonzero.
Thus,~the vector of signs in
$(\pm1)^{10}$
of the ten products
$$\det
\left(
\begin{array}{ccc}
x_{i_1} & x_{i_2} & x_{i_3} \\
y_{i_1} & y_{i_2} & y_{i_3} \\
z_{i_1} & z_{i_2} & z_{i_3}
\end{array}
\right)
\cdot
\det
\left(
\begin{array}{ccc}
x_{i_4} & x_{i_5} & x_{i_6} \\
y_{i_4} & y_{i_5} & y_{i_6} \\
z_{i_4} & z_{i_5} & z_{i_6}
\end{array}
\right) ,
$$
for
$\{i_1,\ldots,i_6\} = \{1,\ldots,6\}$,
$i_1 < i_2 < i_3$,
$i_4 < i_5 < i_6$,
and
$i_1 < i_4$
is an invariant of the configuration, up to the action
of~$S_6$
permuting the six points and up to that
of
$(\bbZ/2\bbZ)^6$,
by sign changes. However,~the latter group acts via the summation character
$+\colon (\bbZ/2\bbZ)^6 \to \bbZ/2\bbZ$.
I.e., there is either no change at all or all ten signs are changed~simultaneously.

Working with a few explicit example configurations, one readily finds four distinct types, having invariant vectors as~follows.

\begin{table}[H]
\begin{center}
\begin{tabular}{|c|c|}
\hline
Type & Multisets associated with the invariant vectors \\\hline\hline
0    & \rule{0pt}{11pt}$\{1^{10}\}$, $\{1^6, (-1)^4\}$,  $\{1^4, (-1)^6\}$, $\{(-1)^{10}\}$ \\\hline
1    & \rule{0pt}{11pt}$\{1^9, (-1)\}$, $\{1^7, (-1)^3\}$,  $\{1^5, (-1)^5\}$, $\{1^3, (-1)^7\}$, $\{1, (-1)^9\}$\\\hline
2    & \rule{0pt}{11pt}$\{1^8, (-1)^2\}$, $\{1^6, (-1)^4\}$,  $\{1^4, (-1)^6\}$, $\{1^2, (-1)^8\}$\\\hline
3    & \rule{0pt}{11pt}$\{1^7, (-1)^3\}$, $\{1^3, (-1)^7\}$\\\hline
\end{tabular}
\end{center}\vskip0mm
\caption{The admissible configurations of six real lines
in~$\Pb^2$}\vskip-4mm
\end{table}
\end{rem}

\begin{rems}
\begin{iii}
\item
We~constructed our example configurations systematically. We~assumed without restriction that the first four points are
$(1\!:\!0\!:\!0)$,
$(0\!:\!1\!:\!0)$,
$(0\!:\!0\!:\!1)$,
and~$(1\!:\!1\!:\!1)$.
The~complement
in~$\Pb^2(\bbR)$
of the six lines through these points breaks into
$\smash{\frac{6\cdot5}2 + 1 - 4\cdot\frac{2\cdot1}2 = 12}$
chambers. Thus, for the choice of the fifth pint, there are only twelve essentially different cases. Moreover,~each time, the complement of the ten lines through the five points obtained is a disjoint union of
$\smash{\frac{10\cdot9}2 + 1 - 5\cdot\frac{3\cdot2}2 = 31}$
chambers~\cite[Chapter~VII, Section 8.2]{Yo}. This~means that all cases are covered by
$31 \!\cdot\! 12 = 372$
explicit configurations.

Calculating~the invariants, we produced a partition of these 372~configurations into four~subsets. We~then checked that two configurations in the same subset were indeed always equivalent to each other, under the operation
of~$S_6$,
followed by a renormalisation of the first four points to
$(1\!:\!0\!:\!0)$,
$(0\!:\!1\!:\!0)$,
$(0\!:\!0\!:\!1)$,
and~$(1\!:\!1\!:\!1)$,
plus a naive relocation of points within the same~chamber.
\item
This~means that our approach actually provides a new proof for the fact that there are exactly four types of configurations of six real lines
in~$\Pb^2$,
no three of which have a point in common. Our~proof is essentially computational and thus rather different from the one presented in~\cite{Yo}.
\item
The group
$S_6$,
of course, does not give rise to all permutations
in~$S_{10}$,
but to a subgroup
of~$S_{10}$
of index
$5040$.
We~ignored the information that is perhaps carried by this~subgroup and just took care of the multisets that are associated with the invariant vectors.
\item
Plotting an example of each type and counting triangles, quadrangles, pentagons, and hexagons, respectively, one finds that the four types correspond, in this order, to the types O, I, II, and~III in the terminology of M.\ Yoshida~\cite[Chapter~VII, Table in Section 8.1]{Yo}.

\end{iii}
\end{rems}\pagebreak[4]

\section{Real and complex multiplication. Hodge structures}
\label{RMCM_Hodge}

\subsubsection*{The concepts}

\begin{defi}[Deligne]
\begin{iii}
\item
A (pure
$\bbQ$-)%
{\em Hodge structure\/} of weight
$i$
is a finite dimensional
$\bbQ$-vector
space
$V$
together with a decomposition
$$V_\bbC := V \!\otimes_\bbQ\! \bbC = H^{0,i} \oplus H^{1,i-1} \oplus \ldots \oplus H^{i,0} \, ,$$
having the property that
$\overline{H^{m,n}} = H^{n,m}$
for every
$m,n \in \bbZ_{\ge0}$
such that
$m+n=i$.

A~{\em morphism\/}
$f\colon V \to V'$
{\em of (pure\/
$\bbQ$-)
Hodge structures\/} is a
$\bbQ$-linear
map such that
$f_\bbC\colon V_\bbC \to V'_\bbC$
respects the decompositions.
\item
A {\em polarisation\/} on a pure
$\bbQ$-Hodge
structure
$V$
of even weight is a nondegenerate symmetric bilinear form
$(\cdot,\cdot)\colon V \times V \to \bbQ$
such that its
$\bbC$-bilinear
extension
$(\cdot,\cdot)\colon V_\bbC \times V_\bbC \to \bbC$
satisfies the two conditions~below.
\begin{iii}
\item[ $\bullet$ ]
One has
$(x,y) = 0$
for all
$x \in H^{m,n}$
and
$y \in H^{m',n'}$
such that
$m \neq n'$.
\item[ $\bullet$ ]
The inequality
$\mi^{m-n} (x,\overline{x}) > 0$
is true for every
$0 \neq x \in H^{m,n}$.
\end{iii}
\end{iii}
\end{defi}

\begin{rem}
The weight
$i$
Hodge structures form an abelian category \cite[2.1.11]{De71}. Furthermore,~the subcategory of polarisable Hodge structures is semisimple \cite[Lemme 4.2.3.i)]{De71}. I.e.,~every sub-Hodge structure of a polarisable Hodge structure is a direct summand.
Thus, every polarisable Hodge structure can be written as a direct sum of indecomposable subobjects, which are called {\em primitive\/} Hodge structures.
\end{rem}

\begin{defi}[Zarhin]
A {\em Hodge structure of
$K3$~type\/}
is a primitive polarisable Hodge structure of
\mbox{weight~$2$}
such
that~$\dim_\bbC H^{2,0} = 1$.
\end{defi}

\begin{exs}
Let~$X$
be a compact complex manifold that is~K\"ahler.

\begin{iii}
\item
Then
$H^j(X,\bbQ)$
is naturally a polarisable pure
\mbox{$\bbQ$-Hodge}
structure of
weight~$j$.
\item
For~$X$
a
$K3$~surface,
the {\em transcendental part\/}
$T := (\Pic X \!\otimes_\bbZ\! \bbQ)^\perp \subset H^2(X, \bbQ)$
is a~Hodge structure of
$K3$~type.
Indeed,~if
$T$
would split off a direct
summand~$T'$
then, without restriction, one had
$T'{}^{2,0} = 0$.
But~this yields, according to the Lefschetz
$(1,1)$-Theorem,
$T'$
that must be contained in the algebraic part, a~contradiction.
\end{iii}
\end{exs}

\begin{prop}[Zarhin]
\label{Zarh}
Let\/~$T$
be a Hodge structure of
$K3$~type.

\begin{iii}
\item
Then\/
$E := \End_\Hg(T)$
is either a totally real field or a CM~field.
\item
Suppose~that\/
$T$
is equipped with a polarisation
$(\cdot,\cdot)$.
Then~every\/
$\varphi \in E$
operates as a self-adjoint~mapping, i.e.\
$(\varphi(x), y) = (x, \overline\varphi(y))$
for\/
$\overline{\phantom{\iota}}$
the identity map in the case that\/
$E$
is totally real and the complex conjugation in the case that\/
$E$
is a CM~field.
\end{iii}%\smallskip

\noindent
{\bf Proof.}
{\em
Assertion~i) is \cite[Theorem~1.6.a)]{Za}, while assertion~ii) is \cite[Theorem 1.5.1]{Za}.
}
\eop
\end{prop}

\subsubsection*{Real and complex multiplication--Our terminology}

\begin{defi}[RM and CM]
Let~$T \subset H^2(X, \bbQ)$
be the transcendental part of the cohomology of a complex
$K3$~surface~$X$.

\begin{iii}
\item
If
$\End_\Hg(T) \supsetneqq \bbQ$
is a totally real field then
$X$
is said to have {\em real multiplication}.
\item
If~$\End_\Hg(T)$
is CM then we speak of {\em complex multiplication}.
\end{iii}
\end{defi}

\begin{rem}
For~complex multiplication, some authors, e.g.\ L.~Taelman \cite{Tae}, require, in addition, that
$\dim_{\End_\Hg(T)} T = 1$.
However,~for our purposes, the definition above appears to be more~practical.
\end{rem}

\begin{rem}
It is a rather common phenomenon that only a subfield of the endomorphism field
$\End_\Hg(T)$
is~known. Moreover, if the generic member of a family of
$K3$~surfaces
has real or complex multiplication by a certain field
$K$
then for particular members the endomorphism field may well be larger, cf.~Corollary \ref{sem_con},~below.\smallskip

\noindent
We~introduce the terminology below in order to cope with such situations.
\end{rem}

\begin{defi}[$T' \supseteq T$
being acted upon by some field]
\label{act_upon}
\leavevmode

\noindent
Let~$X$
be a
$K3$~surface
and
$T' \subseteq H^2(X, \bbQ)$
be a subvector space such that
$T' \supseteq T$,
for
$T \subset H^2(X, \bbQ)$
the transcendental part. We~say that
$T'$
is {\em acted upon\/} by a
field~$K$,
if there is a ring homomorphism
$K \to \End_\Hg(T')$
taking
$1$
to~$\id_{T'}$.
\end{defi}

\begin{rems}
\begin{iii}
\item
In the situation of Definition~\ref{act_upon},
$T'$
is automatically a sub-Hodge structure
of~$H^2(X, \bbQ)$.
\item
If~$T'$
is acted upon by
$K$
then
$T$
is clearly acted upon
by~$K$.

Moreover,~if~$T$
is acted upon
by~$K \supsetneqq \bbQ$
then
$X$
has real or complex multiplication by
$K$
or an extension field
of~$K$.
\end{iii}
\end{rems}

\subsubsection*{Arithmetic effects caused by real or complex multiplication}

\begin{defi}
Let
$X$
be an algebraic
$K3$~surface
defined over a
field~$k$
that is finitely generated
over~$\bbQ$.
Then~we say that
$X$
has {\em real multiplication\/} or {\em complex multiplication\/} by a number
field~$E$
if, for a certain embedding
$k \hookrightarrow \bbC$
of fields, the complex
$K3$~surface
$X(\bbC) = (X \times_{\Spec k} \Spec\bbC)(\bbC)$~has.
\end{defi}

\begin{rem}
This~is independent of the choice of an embedding
$k \hookrightarrow \bbC$
\cite[Corollary~4.2]{EJ20}.
The endomorphism field itself is independent, for
$\dim_\bbQ T < 7$
\cite[Corollary~4.2.b)]{EJ20} or when
$k$
is primary \cite[Corollary~4.4]{EJ20}.
\end{rem}

For~a
$K3$~surface
over a number field, the properties of RM and~CM cause a particular {\em arithmetic\/} behaviour. We summarise the known effects~below.

\begin{nota}
Let~$X$
be a
$K3$~surface
over~$\bbQ$
and
$p$
a prime~number. We~say that
$X$
has {\em good reduction\/}
at~$p$,
if there exists a proper model
$\calX$
of~$X$
over~$\bbZ$,
the reduction
$X_p := \calX_p$
modulo~$p$
of which is again a
$K3$~surface.
In~this case,
for~$f \in \bbN$,
we let
$\chi_{p^f} \in \bbQ[Z]$
be the characteristic polynomial of
$\smash{\Frob^f \in \Gal(\overline\bbF_{\!p}/\bbF_{\!p})}$
on
$\smash{H^2_\et((X_p)_{\overline\bbF_{\!p}}, \bbQ_l(1))}$
\cite[Th\'eor\`eme~1.6]{De74}.

We~factorise
$\smash{\chi_{p^f}}$
completely in the form
$$\chi_{p^f}(Z) = \chi^\tr_{p^f}(Z) \cdot \prod_{i=1}^d (Z - \zeta_{k_i}^{e_i}) \,,$$
for
$k_1,\ldots,k_d \in \bbN$.
I.e.\ in such a way that
$\smash{\chi^\tr_{p^f} \in \bbQ[Z]}$
does not have any further zeroes being roots of~unity.

According~to the Tate conjecture \cite{MP,KM},
$\smash{\chi^\tr_{p^f}}$
is the characteristic polynomial of
$\smash{\Frob^f}$
on the transcendental part of
$\smash{H^2_\et(X_{\overline\bbF_{\!p}}, \bbQ_l(1))}$.
Let~us note, in~particular, that
$\smash{\deg\chi^\tr_{p^f} = 22 - \rk\Pic X_{\overline\bbF_{\!p}}}$.
\end{nota}

\begin{prop}[Some arithmetic effects of real or complex multiplication]
\label{arith_eff}
\leavevmode\\
Let\/
$X$
be a\/
\mbox{$K3$~surface}
over\/~$\bbQ$.
Suppose that\/
$X$
has real or complex multiplication by a number
field\/~$E$.
Let,~moreover,
$p$
be a prime of good~reduction.

\begin{ABC}
\item
Then\/~$\rk\Pic (X_p)_{\overline\bbF_{\!p}}- \rk\Pic X(\bbC)$
is divisible
by\/~$[E\!:\!\bbQ]$.
\item
Assume~that\/
$\smash{E \supseteq \bbQ(\sqrt{\delta})}$,
for a non-square\/
$\delta \in \bbQ^*$.
\begin{abc}
\item
Then~one of the following is~true.
Either\/~$\smash{\chi^\tr_p}$
splits over\/
$\smash{\bbQ(\sqrt{\delta})}$.
Or,~for a certain\/
$f > 0$,
the polynomial\/
$\smash{\chi^\tr_{p^f}}$
is a~square.
\item
If~the
prime~\/~$p$
is inert
in\/~$\smash{\bbQ(\sqrt{\delta})}$
then the reduction
$X_p$\/
is non-ordinary. I.e., one has that\/
$\#X_p(\bbF_{\!p}) \equiv 1 \pmod p$.
\end{abc}%
\end{ABC}\smallskip

\noindent
{\bf Proof.}
{\em
A)
This~result appears between the lines of F.\ Charles' famous article~\cite{Ch}, but was not explicitly stated. A proof is provided in \cite[Lemma~6.2]{EJ20}.\smallskip

\noindent
B)
Assertion~B.a) is shown in \cite[Theorem~4.9]{EJ14}, while B.b) is \cite[Corollary 4.13.i)]{EJ14}.
}
\eop
\end{prop}

\begin{prop}[Sufficient criterion for real or complex multiplication]
\label{suff_crit}
\leavevmode\\
Let\/~$a,D \in \bbZ$
be such that\/
$\gcd(a,D) = 1$
and\/~$X$
a\/~$K3$~surface
over\/~$\bbQ$.
Suppose~that\/
$\#X_p(\bbF_{\!p}) \equiv 1 \pmod p$
for every good
prime\/~$p \equiv a \pmod D$.
Then\/~$X$
has real or complex~multiplication.\smallskip

\noindent
{\bf Proof.}
{\em
This~is \cite[Lemma~6.1]{EJ14}.
}
\eop
\end{prop}

\begin{rem}
A relative version of Proposition~\ref{suff_crit} is established in~\cite[Theorem~3.5]{EJ20}.
\end{rem}

\subsubsection*{Some known examples I}

\begin{ex}
\label{CM_ex1}
Let~$f_4$
be a homogeneous quartic form in three variables that defines a regular curve
$C \subset \Pb^2$.
Then~the fourfold cover
$X\colon w^4 = f_4(x,y,z)$
ramified
at~$C$
is a
$K3$~surface
of
degree~$4$
and Picard rank at
least~$8$.
It~has an automorphism, given~by
$$I\colon (w\!:\!x\!:\!y\!:\!z) \mapsto (\mi w\!:\!x\!:\!y\!:\!z) \,.$$
Moreover,~there is a
\mbox{$14$-dimensional}
subvector space
$T' \subset H^2(X, \bbQ)$
containing~$T$
that is acted upon
by~$\bbQ(\sqrt{-1})$.
I.e.,~$X$
has complex multiplication by
$\bbQ(\sqrt{-1})$
or by a field that properly contains
$\bbQ(\sqrt{-1})$.\smallskip

\noindent
{\bf Proof.}
The automorphism
$J := I \!\circ\! I$
is of
order~$2$
and has the
curve~$C$
as its fixed point~set. The~genus
of~$C$
is~$3$,
hence the topological Euler characteristic is equal
to~$(-4)$.
Therefore, the Lefschetz trace formula~\cite[Theorem~8.5]{Ed}, cf.~\cite[Expos\'e~III, formule~(4.11.3)]{SGA5}, shows~that
$\Tr J |_{H^2(X, \bbQ)} = -6$.
In~other words,
$J |_{H^2(X, \bbQ)}$
has the eigenvalue
$1$
with
multiplicity~$8$,
while the eigenvalue
$(-1)$
occurs with multiplicity~$14$.

On the other hand,
$X$
is a double cover of the degree two del Pezzo
surface~$X'$,
given by
$w^2 = f_4(x,y,z)$.
One~has
$\Pic X' \cong \bbZ^8$
and the natural homomorphism
$\pi^*\colon \Pic X' \to \Pic X$
is an injection, as it doubles all intersection numbers.
Thus,~$J$
acts nontrivially, with only eigenvalue
$(-1)$,
on the 
\mbox{$14$-d}i\-men\-sional
orthogonal complement of
$\im(\pi^*\colon H^2(X', \bbQ) \to H^2(X, \bbQ))$.
%In~particular, it does so on the transcendental part
%$T \subset H^2(X, \bbQ)$.
\eop
\end{ex}

\begin{ex}
Let~$f_2$
and~$f_3$
be homogeneous forms in four variables such that the subscheme
$X \subset \Pb^4$,
given~by
$f_2(x,y,z,u) = 0$
and~$v^3 = f_3(x,y,z,u)$
is~regular.
Then~$X$
is a
$K3$~surface
of
degree~$6$.
It~has an automorphism, given~by
$$I\colon (v\!:\!x\!:\!y\!:\!z\!:\!u) \mapsto (\zeta_3 v\!:\!x\!:\!y\!:\!z\!:\!u) \,.$$
Moreover,~there is a
\mbox{$20$-dimensional}
subvector space
$T' \subset H^2(X, \bbQ)$
containing~$T$
that is acted upon
by~$\bbQ(\sqrt{-3})$.\smallskip

\noindent
{\bf Proof.}
The automorphism
$I$
is of
order~$3$
and has the~curve
$$C\colon \quad f_2(x,y,z,u) = f_3(x,y,z,u) = 0 \,,\; v=0 \,,$$
as its fixed point~set. This~is a canonical curve of
genus~$4$,
such that the Lefschetz trace formula yields
$\smash{\Tr J |_{H^2(X, \bbQ)} = -8}$.
In~other words, the eigenvalue
$1$
has
multiplicity~$2$,
while the eigenvalues
$\zeta_3$
and~$\smash{\zeta_3^{-1}}$
both occur with
multiplicity~$10$.

On the other hand,
$X$
is a threefold cover of the space quadric, given by
$X'\colon f_2(x,y,z,u) = 0$.
One~has
$\Pic X' \cong \bbZ^2$
and the natural homomorphism
$\Pic X' \to \Pic X$
is an injection, tripling all intersection numbers.
Thus,~$J$
acts nontrivially, with eigenvalues only
$\zeta_3$
and~$\smash{\overline\zeta_3}$,
on the 
\mbox{$20$-d}i\-men\-sional
orthogonal complement of the image of
$H^2(X', \bbQ)$
in~$H^2(X, \bbQ)$.
%In~particular, it does so on the transcendental part
%$T \subset H^2(X, \bbQ)$.
\end{ex}

\begin{ex}
Let~$f_6$
be a homogeneous sextic form in three variables that defines a regular curve
$C \subset \Pb^2$.
Then~the double cover
$X\colon w^2 = f_6(x,y,z)$
ramified
at~$C$
is a
$K3$~surface
of
degree~$2$.
Suppose that
$f_6(x,y,z) = \zeta_3 f_6(y,z,x)$.
Then~$X$
has an automorphism, given~by
\begin{align*}
I\colon (w, x\!:\!y\!:\!z) &\mapsto (\zeta_6 w, y\!:\!z\!:\!x) \,.
\end{align*}
Moreover,~there is a
\mbox{$14$-dimensional}
subvector space
$T' \subset H^2(X, \bbQ)$
containing~$T$
that is acted upon
by~$\bbQ(\sqrt{-3})$.
I.e.,~$X$
has complex multiplication either by
$\bbQ(\sqrt{-3})$
or by a field that contains
$\bbQ(\sqrt{-3})$.
Furthermore,~$\rk\Pic X \geq 8$.\smallskip

\noindent
{\bf Proof.}
The automorphism
$J := I \!\circ\! I \!\circ\! I$
of~$X$
is of
order~$2$
and has the
curve~$C$
as its fixed point~set. The~genus
of~$C$
is~$10$,
hence the topological Euler characteristic is equal
to~$(-18)$.
Therefore, the Lefschetz trace formula shows~that
$\Tr J |_{H^2(X, \bbQ)} = -20$.
This~means that
$J$
acts nontrivially, with only eigenvalue
$(-1)$,
on the 
\mbox{$21$-d}i\-men\-sional
orthogonal complement of the inverse image of a general line on
$\Pb^2$.
Consequently,~$I$
acts on this space with eigenvalues at most
$(-1)$,
$\zeta_6$,
and~$\zeta_6^{-1}$.

On~the other hand, the fixed point set
of~$I$
is
$\{(0, 1\!:\!1\!:\!1), (0, 1\!:\!\zeta_6^2\!:\!\zeta_6^4), 0, 1\!:\!\zeta_6^4\!:\!\zeta_6^2)\}$,
which yields that
$\Tr I |_{H^2(X, \bbQ)} = 1$.
Consequently, each of the eigenvalues
$(-1)$,
$\zeta_6$,
and~$\zeta_6^{-1}$
occurs with
multiplicity~$7$.
The~direct sum of the eigenspaces for
$\zeta_6$,
and~$\zeta_6^{-1}$
has the property required. Finally,~the assertion on the Picard rank follows from \cite[Theorem 1.4.1]{Za}.
\eop
\end{ex}

\begin{ex}
\label{CM_ex4}
Let~$f_5$
be a homogeneous quintic form in three variables that defines a regular curve
$C \subset \Pb^2$.
Suppose that
$C \cap V(z)$
consists of five distinct points.
Then~the minimal
desingularisation~$X$
of the double cover
$X'\colon w^2 = z f_5(x,y,z)$
is a
$K3$~surface
of Picard rank at
least~$6$.
Suppose that
$f_5(x,y,-z) = f_5(x,y,z)$.
Then~$X$
has an automorphism, given~by
$$I\colon (w, x\!:\!y\!:\!z) \mapsto (i w, x\!:\!y\!:\!-z) \, .$$
Moreover,~there is a
\mbox{$16$-dimensional}
subvector space
$T' \subset H^2(X, \bbQ)$
containing~$T$
that is acted upon
by~$\bbQ(\sqrt{-1})$.
I.e.,~$X$
has complex multiplication either by
$\bbQ(\sqrt{-1})$
or by a field that contains
$\bbQ(\sqrt{-1})$.\smallskip

\noindent
{\bf Proof.}
Here,~the automorphism
$J := I \!\circ\! I$
of~$X$
is of
order~$2$.
Its~fixed point set is the union of the strict transforms of
$C$
and~$V(z)$
in~$\Pb^2$,
blown up in the five points
of~$C \cap V(z)$.
Note~that this is a disjoint union. As~the genera of
$C$
and~$V(z)$
are
$6$
and~$0$,
the topological Euler characteristic
of~$C \cup V(z)$
is equal
to~$(-10) + 2 = (-8)$.
Thus,~the Lefschetz trace formula implies
$\Tr J |_{H^2(X, \bbQ)} = -10$.
This~means that
$J$
acts nontrivially, with only eigenvalue
$(-1)$,
on the 
\mbox{$16$-d}i\-men\-sional
orthogonal complement of the span of the five exceptional curves and the inverse image of a general line on
$\Pb^2$.
\eop
\end{ex}

\begin{rem}
More generally, when a
$K3$~surface
$X$
has an automorphism of finite order that operates nontrivially
on~$H^{2,0}(X)$
then
$X$
has complex multiplication. More precisely,
$T \subset H^2(X, \bbQ)$
is acted upon by a cyclotomic~field. Such
$K3$~surfaces
have been intensively studied, there is even a classification for the case that the order is a prime~number. The~interested reader is advised to consult the article \cite{AST} and the references given~therein.
\end{rem}

\subsubsection*{Some known examples II}\leavevmode

\noindent
Real multiplication tends to be more complicated than complex multiplication. In~particular, no examples are known, in which real multiplication occurs due to an~automorphism.
Instead, the examples presented below were found searching excessively for
$K3$~surfaces
showing the arithmetic abnormalities predicted by Proposition~\ref{arith_eff}.

\begin{ex}[A family of
$K3$~surfaces
for which RM
by~$\bbQ(\sqrt{2})$
is established]
\label{Qw2_fam}
Let
$\smash{q\colon \frakX \!\to\!\! B}$,
for~$\smash{B \!:=\! \Spec \bbQ[T\!,\!\frac1{T(T^2-2)(T^2+2)(T^2-4T+2)(T^2+4T+2)}] \!\subset\! \Ab^1_\bbQ}$,
be the family of
$K3$~surfaces,
the fibre
at~$t \in B$
of which is the minimal desingularisation of the double cover
of~$\Pb^2$,
given~by
\begin{eqnarray*}
w^2 & = & \textstyle
 [(\frac18 t^2 \!-\! \frac12 t \!+\! \frac14)y^2 + (t^2 \!-\! 2t \!+\! 2)yz + (t^2 \!-\! 4t \!+\! 2)z^2] \\[-1mm]
 & & \textstyle\hspace{4mm}
 [(\frac18 t^2 \!+\! \frac12 t \!+\! \frac14)x^2 + (t^2 \!+\! 2t \!+\! 2)xz + (t^2 \!+\! 4t \!+\! 2)z^2] [2x^2 + (t^2 \!+\! 2)xy + t^2y^2] \, .
\end{eqnarray*}
\begin{iii}
\item
Then the geometric generic fibre
$\frakX_{\overline\eta}$
of~$q$
is of Picard
rank~$16$.
\item
For~every
$\theta \in B(\bbC)$,
the transcendental part
$\smash{T \subset H^2(\frakX_\theta(\bbC), \bbQ)}$
of the cohomology of the fibre
$\frakX_\theta(\bbC)$
of the holomorphic submersion
$q(\bbC)\colon \frakX(\bbC) \to B(\bbC)$
is acted upon
by~$\smash{\bbQ(\sqrt{2})}$.
\item
Let the complex point
$\theta \in B(\bbC)$
be of the kind that the fibre
$\frakX_\theta(\bbC)$
of~$q(\bbC)$
has Picard
rank~$16$.
Then
$\frakX_\theta(\bbC)$
has real multiplication
by~$\smash{\bbQ(\sqrt{2})}$.
\end{iii}\smallskip

\noindent
{\bf Proof.}
This is \cite[Example~5.1]{EJ20}. The~family was presented in~\cite{EJ14} for the first~time.
\eop
\end{ex}

\begin{rem}
For~$\theta \in \bbQ$,
the field of definition of the Picard group is
$\smash{\bbQ(\sqrt{2})}$,~too.
\end{rem}

\begin{ex}[A family of
$K3$~surfaces
for which RM
by~$\bbQ(\sqrt{5})$
is established]
\label{Qw5_fam}
Let
$\smash{q\colon \frakX \!\to\!\! B}$,
for
$\smash{B \!:=\! \Spec \bbQ[T\!,\!\frac1{(T-1)(T^4-T^3+T^2-T+1)}] \!\subset\! \Ab^1_\bbQ}$,
be the family of
$K3$~surfaces,
the fibre
at~$t \in B$
of which is the minimal desingularisation of the double cover
of~$\Pb^2$,
given~by
\begin{eqnarray*}
w^2 &=& y(x - 2(t\!-\!1)y - tz) \\[-1mm]
&& (x^4 + x^3y - x^3z + x^2y^2 - 2x^2yz + x^2z^2 + xy^3 - 3xy^2z - 2xyz^2 - xz^3 + y^4 \nonumber \\[-1mm]
&& \hspace{9.2cm} {}+ y^3z + y^2z^2 + yz^3 + z^4) \, .
\end{eqnarray*}
\begin{iii}
\item
Then the generic fibre
$\frakX_\eta$
of~$q$
is of geometric Picard
rank~$16$.
\item
For~every
$\theta \in B(\bbC)$,
the transcendental part
$\smash{T \subset H^2(\frakX_\theta(\bbC), \bbQ)}$
of the cohomology of the fibre
$\frakX_\theta$
is acted upon
by~$\smash{\bbQ(\sqrt{5})}$.
\item
Let the complex point
$\theta \in B(\bbC)$
be of the kind that the fibre
$\frakX_\theta$
has geometric Picard
rank~$16$.
Then
$\frakX_\theta$
has real multiplication
by~$\smash{\bbQ(\sqrt{5})}$.
\end{iii}\smallskip

\noindent
{\bf Proof.}
This is \cite[Example~1.5]{EJ20}.
\eop
\end{ex}

\begin{rem}
Here,~for
$\theta \in \bbQ$,
the field of definition of the Picard group is
$\smash{\bbQ(\zeta_5)}$.
\end{rem}

\begin{rems}
\label{conj_RM}
\begin{iii}
\item
Several further families have been found, for the generic members of which there is strong evidence for real multiplication by specific fields, including
$\smash{\bbQ(\sqrt{2})}$,
$\smash{\bbQ(\sqrt{3})}$,
and~$\smash{\bbQ(\sqrt{5})}$.
Furthermore,~there are isolated examples that conjecturally have complex multiplication by
$\bbQ(\sqrt{-1}, \zeta_7 + \zeta_7^{-1})$,
$\bbQ(\sqrt{-1}, \zeta_9 + \zeta_9^{-1})$,
and
$L(\sqrt{-1})$,
respectively, for
$L \subset \bbQ(\zeta_{19})$
the unique cubic~subfield. Explicit~equations are given in~\cite[Conjecture 5.2.b)]{EJ16}.

For~each of these surfaces, it has been verified that the arithmetic effects caused by real or complex multiplication (cf.\ Proposition~\ref{arith_eff}) occur for all primes
$p < 500$.
This~required fast algorithms for point counting on algebraic surfaces over finite fields, cf.~\cite{Ha} or~\cite{EJ16} for the main ideas behind~them. By~now, an implementation of the method we used is available to the public, as the {\tt magma} intrinsic {\tt WeilPolynomialOfDegree2K3Surface}.
\item
Moreover,~there are a few isolated examples
for~$\bbQ(\sqrt{13})$.
For~instance, there is the~following.
\end{iii}
\end{rems}

\begin{ex}
\label{qw13}
Let~$X$
be the minimal desingularisation of the double cover
of~$\Pb^2$,
given~by
$$w^2 = (x - 4z)(5x - 9y - 8z)f_4(x,y,z) \,,$$
for~$f_4 := x^4 - 2x^3y - 5x^2y^2 - 26x^2z^2 + 6xy^3 + 104xyz^2 + 9y^4 - 130y^2z^2 + 52z^4$.\smallskip

\noindent
\begin{iii}
\item
Then~$X$
is a
$K3$~surface
over~$\bbQ$
of geometric Picard
rank~$16$.
\item
There~is strong evidence that
$X(\bbC)$~has
real multiplication
by~$\bbQ(\sqrt{13})$.
\end{iii}\medskip

\noindent
{\bf Proof} of i).
The quartic
form~$f_4$
is the norm of a linear form over the cyclic quar\-tic number~field
$$\textstyle K = \bbQ(\sqrt{13 - 3 \sqrt{13}})$$
of
conductor~$8 \cdot 13 = 104$.
Thus,~the branch locus is the union of six lines in general~position. This shows that
$\smash{\rk\Pic X_{\overline\bbQ} \geq 16}$,
while an upper bound
of~$17$
may, once again, be obtained using van Luijk's~method. Let~us note that
$X$
has bad reduction only at the primes
$2$,
$3$,
and~$13$.

In order to verify that
$\smash{\rk \Pic X_{\overline\bbQ} \neq 17}$,
a modification of the method described in \cite{EJ11}~applies. Indeed,~the
$\smash{\Gal(\overline\bbQ/\bbQ)}$-representation
$\smash{H^2_\et((X_p)_{\overline\bbF_{\!p}}, \overline\bbQ_l(1))}$
splits into a direct summand of
dimension~$16$
corresponding to the obvious rank 16 part of
$\smash{\Pic X_{\overline\bbQ}}$
and a
complement~$V_6$.
Geometric Picard
rank~$17$
would cause a free
\mbox{$\bbZ$-module}
of rank one being contained 
in~$V_6$
that is acted upon
by~$\smash{\Gal(\overline\bbQ/\bbQ)}$.
However, the characteristic polynomial of
$\smash{\Frob_5}$
on~$V_6$
turns out to be
$\smash{(Z^2 + 1)(Z^4 - \frac25 Z^2 + 1)}$.
Thus,~there is no eigenvalue
$\pm1$,
which is a~contradiction.
\eop\smallskip

\noindent
ii)
Evidence for real multiplication is as~follows. For~every good
prime~$p < 500$,
the characteristic polynomial
$\smash{\chi^\tr_p}$
splits over
$\smash{\bbQ(\sqrt{13})}$
or
$\smash{\chi^\tr_{p^2}}$
is a square, as predicted by Proposition~\ref{arith_eff}.B.a). Moreover,~for every good prime
$p \equiv 2,5,6,7,8,11 \pmod {13}$
up
to~$1000$,
one has
$\smash{\#X(\bbF_{\!p}) \equiv 1 \pmod p}$,
cf.~Proposition~\ref{arith_eff}.B.b). In~fact,
$\smash{\#X'(\bbF_{\!p}) = p^2+p+1}$
for~$X'$
the double cover
of~$\Pb^2$
underlying~$X$.
\end{ex}

\begin{rem}
For~the surface in Example~\ref{qw13}, the Picard group is defined
over
$\smash{K = \bbQ(\sqrt{13 - 3 \sqrt{13}})}$.
\end{rem}

\subsubsection*{Periods of
$K3$~surfaces
having real or complex multiplication}

\begin{theo}
\label{RMCM_dim}
Let\/~$r \in \{1,\ldots,20\}$
be an integer and\/
$\kappa$
a perfect pairing
on\/~$\bbZ^{22}$.

\begin{abc}
\item[{\rm a) (Periods of the families generically having RM.)} ]
\leavevmode\\
Let\/~$K$
be a totally real number field of
degree\/~$d$.
Then~there is an at most countable
union\/~$M_{K,\kappa,r} \subseteq Q_{\kappa,r}$
of\/ {\em quadrics of dimension}
$\smash{\frac{22-r}{d}-2}$
such that the following~holds.

\looseness-1
Let\/~$x \in \Omega_{\kappa,r} \subset Q_{\kappa,r}$
be the period point of a marked\/
$K3$~surface
$(X, i)$,
for which\/
$c_k \in H^2(X, \bbQ)$
(cf.~\ref{mark}.ii)) is algebraic,
for\/~$k = 22-r+1, \ldots, 22$,
and the Picard rank
of\/~$X$
is exactly\/~$r$.
Then\/~$T \subset H^2(X, \bbQ)$
is acted upon
by\/~$K$
if and only if\/
$x \in M_{K,\kappa,r}$.
\item[{\rm b) (Periods of the families generically having CM.)} ]
\leavevmode\\
Let\/~$K$
be a CM field of
degree~$d$.
Then~there is an at most countable union\/
$M_{K,\kappa,r} \subseteq Q_{\kappa,r}$
of\/ {\em projective subspaces of dimension}
$\smash{\frac{22-r}{d}-1}$
such that the following is~true.

\looseness-1
Let\/~$x \in \Omega_{\kappa,r} \subset Q_{\kappa,r}$
be the period point of a marked\/
$K3$~surface
$(X, i)$,
for which\/
$c_k \in H_2(X, \bbQ)$
(cf.~\ref{mark}.ii)) is algebraic,
for\/~$k = 22-r+1, \ldots, 22$,
and the Picard rank
of\/~$X$
is exactly\/~$r$.
Then\/~$T \subset H^2(X, \bbQ)$
is acted upon
by\/~$K$
if and only if\/
$x \in M_{K,\kappa,r}$.
\end{abc}%\smallskip

\noindent
{\bf Proof.}
{\em
By~Lemma \ref{periods_Pic}, one has that
$\Pic X = \big(i(x)\big)^\perp \supseteq \spann(c^1, \ldots, c^{22-r})^\perp$.
However,~the assumption
$\rk\Pic = r$
ensures equality, such that the transcendental lattice
of~$X$~is
$T = (\Pic X \!\otimes_\bbZ\! \bbQ)^\perp = \spann(c^1, \ldots, c^{22-r})$.\smallskip

\noindent
a)
Being acted upon by the totally real
field~$K$
means that
$K$~operates
\mbox{$\bbQ$-linearly}
on~$T$,
keeping the distinguished cohomology class
$x_1 c^1 + \cdots + x_{22-r} c^{22-r}$
as a simultaneous~eigenvector.

For~this, take a primitive element
$u \in K$.
To~fix the operation
of~$K$,
one simply has to choose a self-adjoint endomorphism
$U\colon T \to T$,
whose minimal polynomial is the same as that
of~$u$.
If~this is possible then
$d \mid (22-r)$.
Moreover,~there are certainly only countably many choices. Thus,~let us consider
$U$
as being~fixed.

Then~$U_\bbR$
has
$d$~eigenvalues,
the numbers
$\sigma_1(u), \ldots, \sigma_d(u)$,
for
$\sigma_i\colon K \to \bbR$
the real embeddings, and the eigenspaces
$E_1, \ldots, E_d$
are of dimension
$\smash{\frac{22-r}{d}}$~each.
Since~$U$
is self-adjoint, the eigenspaces are perpendicular to each~other. In~particular, they are non-degenerate quadratic~spaces.

To~complete the proof, let us note that, for
$T$~being
acted upon
by~$K$,
the vector
$x_1 c^1 + \cdots + x_{22-r} c^{22-r}$
needs to be contained in
$E_i \!\otimes_\bbR\! \bbC$,
for one of the eigenspaces.\smallskip

\noindent
b)
Take a primitive element
$u \in K$.
Similarly~to the above, for
$T$
being acted upon by the CM field
$K$,
one has to choose an endomorphism
$U\colon T \to T$,
whose minimal polynomial is the same as that
of~$u$.
Again,~if this is possible then
$d \mid (22-r)$,
and, again, there are only countably many choices, so that we may consider
$U$
as being~fixed.

Then~$U_\bbC$
has
$d$~eigenvalues,
the numbers
$\sigma_1(u), \ldots, \sigma_d(u)$,
for
$\sigma_i\colon K \to \bbC$
the complex embeddings, and the eigenspaces
$E_1, \ldots, E_d$
are of dimension
$\smash{\frac{22-r}{d}}$~each.
And,~clearly, for
$T$
being acted upon
by~$K$,
the vector
$x_1 c^1 + \cdots + x_{22-r} c^{22-r}$
needs to be contained in one of the
eigenspaces~$E_i$.
This~provides us with projective subspaces of dimension
$\smash{\frac{22-r}{d}-1}$,
as~desired.

Moreover,~$U$
must fulfil, together with the linear map
$V$,
associated with
$\overline{u}$,
the self-adjointness relation
$(U_\bbC(x), y) = (x, V_\bbC(y))$.
In~particular,
for~$x \in E_i$,
this~yields
$$\sigma_i(u) (x, x) = (U_\bbC(x), x) = (x, V_\bbC(x)) = \overline{\sigma_i(u)} (x, x)$$
and hence
$(x, x) = 0$,
as 
$\sigma_i(u)$
is certainly non-real, when
$K$~is
a CM~field and
$u$
a primitive~element. This~shows that the projective subspaces
$\Pb(E_i)$,
associated with the
eigenspaces~$E_i$,
are already contained in the
quadric~$Q_{\kappa,r}$.
}
\eop
\end{theo}

\begin{coro}[Semicontinuity]
\label{sem_con}
Let\/~$r \in \{1,\ldots,20\}$
be an integer and\/
$\kappa$
a perfect pairing
on\/~$\bbZ^{22}$.

\begin{abc}
\item[{\rm a) (RM)} ]
Let\/~$K$
be a totally real number field of
degree\/~$d$.
Assume that\/
$\smash{\frac{22-r}{d} \geq 3}$.
Then~there is an at most countable
union\/~$V_{K,\kappa,r} \subset M_{K,\kappa,r}$
of analytic subsets such that the following is~true.

Let\/~$x \in \Omega_{\kappa,r} \subset Q_{\kappa,r}$
be the period point of a marked\/
$K3$~surface
$(X, i)$,
for which\/
$c_k \in H^2(X, \bbQ)$
(cf.~\ref{mark}.ii)) is algebraic,
for\/~$k = 22-r+1, \ldots, 22$.
Then\/~$X$
has real multiplication
by\/~$K$
and\/
$\rk \Pic X = r$
if and only if\/
$x \in M_{K,\kappa,r} \setminus V_{K,\kappa,r}$.
\item[{\rm b) (CM)} ]
Let\/~$K$
be a CM field of
degree\/~$d$.
Assume that\/
$\smash{\frac{22-r}{d} \geq 2}$.
Then~there is an at most countable
union\/~$V_{K,\kappa,r} \subset M_{K,\kappa,r}$
of analytic subsets such that the following is~true.

Let\/~$x \in \Omega_{\kappa,r} \subset Q_{\kappa,r}$
be the period point of a marked\/
$K3$~surface
$(X, i)$,
for which\/
$c_k \in H^2(X, \bbQ)$
(cf.~\ref{mark}.ii)) is algebraic,
for\/~$k = 22-r+1, \ldots, 22$.
Then\/~$X$
has complex multiplication
by\/~$K$
and\/
$\rk \Pic X = r$
if and only if\/
$x \in M_{K,\kappa,r} \setminus V_{K,\kappa,r}$.
\end{abc}%\smallskip

\noindent
{\bf Proof.}
{\em
We~present the proof only for~a), as that for~b) is completely~analogous.
There are two ways for the property stated
on~$X$
to~fail. Either
$\rk \Pic X > r$
or
$\rk \Pic X = r$,
but
$\End_\Hg(T) \supsetneqq K$.\smallskip

\noindent
{\em Case~1:\/}
$\rk \Pic X > r$.

\noindent
This means that, besides the linear combinations of
$c_{22-r+1},\ldots,c_{22}$,
a further cohomology class
$c \in H^2(X, \bbZ)$
is~algebraic. According~to Lemma~\ref{periods_Pic}, that is equivalent to
$x \in (i^* c)^\perp$,
which clearly defines an analytic subset. As~there are only countably many possibilities
for~$c$,
this case indeed contributes
to~$V_{K,\kappa,r}$
a countable union of analytic~subsets.\smallskip

\noindent
{\em Case~2:\/}
$\rk \Pic X = r$
and
$\End_\Hg(T) \supsetneqq K$.

\noindent
If~$\End_\Hg(T) = K' \supsetneqq K$
then
$x \in M_{K',\kappa,r}$,
which, according to Theorem~\ref{RMCM_dim}.a), defines a countable union of analytic subsets
of~$M_{K,\kappa,r}$.
As~there are, up to isomorphism, only countably many number fields, the assertion~follows.
}
\eop
\end{coro}

\begin{rem}
This~shows, in particular, the~following. As~long as the Picard rank remains unchanged, the endomorphism field
$\End_\Hg(T)$
cannot shrink under specialisation. This~fact actually has been obtained before and in more generality \cite[Corollary~4.8]{EJ20}. The restrictive assumption on the Picard ranks is in fact~unnecessary.
\end{rem}

\section{Tracing the preimage of a curve in the period space}

\begin{rem}
As~a particular case of Theorem~\ref{RMCM_dim}.a), we see that in a sufficiently general family of
$K3$~surfaces
of Picard
rank~$16$,
not containing an isotrivial subfamily, those surfaces that are acted upon by a {\em real\/} quadratic number field form families over
curves~$C \subset Y$,
cf.~\cite[Example~3.4]{vG}.
\end{rem}

\begin{strat}
\label{str}
Let~$X$
be an isolated example of a
$K3$~surface
that has real multiplication by a quadratic
field~$\smash{\bbQ(\sqrt{d})}$.
Assume~that
$X$
is given as the minimal desingularisation of a double cover of the form~(\ref{eq_norm}). The~strategy below describes how to find the
\mbox{$1$-dimensional}
family of RM surfaces,
$X$
belongs~to.

\begin{iii}
\item
Run~Algorithm~\ref{trans_part}
on~$X$.
The~resulting elements
$c_{\alpha_1}, \ldots, c_{\alpha_6} \in H^2(X, \bbZ)/P$
yield a class of markings
on~$X$,
as described in Theorem~\ref{mark_rel}.

As~step~\ref{item_v}, Algorithm~\ref{trans_part} includes running Algorithm~\ref{cup_prod}. In~particular, open neighbourhoods
$\bbD \cong U(a_0) \ni a_0$,
\ldots,
$\bbD \cong U(d_0) \ni d_0$
are chosen in such a way that, for every
$$(a,b,c,d) \in U := U(a_0) \times \cdots \times U(d_0) \,,$$
no three of the resulting six lines
in~$\Pb^2_\bbC$
have a point in~common. Thus,~any
marking~$i$
from the class above extends to the whole family
over~$U$.
There~is the associated restricted period~map
$$\Pi\colon U \longrightarrow \Pb^5(\bbC), \quad (a,b,c,d) \mapsto \Pi_{X_{(a,b,c,d)}, i_{(a,b,c,d)}} \,,$$
which is independent of the choice
of~$i$,
cf.\ Theorem~\ref{mark_rel}.c).
%For~every concrete point
%$(a,b,c,d) \in \bbD^4$,
%the value
%of~$\Pi$,
%as well as those of the partial derivatives, may be calculated by numerical~integration.
\item
Calculate the period point
of~$X = X_{(a_0,b_0,c_0,d_0)}$,
based on Theorem~\ref{periods_ints}, using numerical integration, and identify the three linear relations between the six periods that encode real multiplication.
These~define, together with the quadric induced by the cup product pairing, a conic
$C \subset \Pb^5(\bbC)$
in the restricted period~space.
\item\label{num_cont}
Trace~the curve
$\Pi^{-1}(C) \subset U \cong \bbD^4$
using a numerical continuation method~\cite{AG}.
\item
Use~the singular-value decomposition in order to find algebraic relations between the coordinates of the points~found. Control,~by using Gr\"obner bases, that they indeed define an irreducible algebraic~curve.
\end{iii}\smallskip

\noindent
We provide a few details on how this strategy was implemented in Remarks~\ref{det_exp}, below.
\end{strat}

\subsubsection*{The result}

\begin{ex}
\label{result}
Consider the family of double
covers\/~$\smash{X'_{(a,b,c,d)}}$
of\/~$\Pb^2$,
given~by
$$w^2 = (x + ay + bz)(x + cy + dz)f_4(x,y,z) \,,$$
for\/~$f_4 := x^4 - 2x^3y - 5x^2y^2 - 26x^2z^2 + 6xy^3 + 104xyz^2 + 9y^4 - 130y^2z^2 + 52z^4$.

\begin{abc}
\item[{\rm a.i) }]
Then~the branch locus is the union of six lines, which are in general position for a generic choice
of\/~$(a,b,c,d) \in \bbC^4$.
\item[{\rm a.ii) }]
Consider the closed subscheme\/
$C \subset \Ab^4$,
given by the~equations
\begin{eqnarray*}
0 &=& 630\,272a - 11\,421bd^5 + 411\,400bd^3 - 871\,552bd - 272\,976c^2d^2 + 315\,136c^2 \\[-1mm]
   && {} + 98\,982cd^4 - 3\,508\,064cd^2 + 2\,205\,952c + 233\,496d^4 - \!6\,409\,856d^2 \!+ 4\,411\,904\,, \\[-.5mm]
0 &=& 78\,784bc - 243bd^4 + 37\,040bd^2 + 110\,528b - 5808c^2d + 2106cd^3 - 319\,792cd \\[-1mm]
   && {} + 4968d^3 - 714\,688d\,,\\[-.5mm]
0 &=& 243bd^6 - 8960bd^4 + 29\,952bd^2 - 26\,624b + 5808c^2d^3 - 11\,648c^2d - 2106cd^5 \\[-1mm]
   && {} + 76\,432cd^3 - 144\,768cd - 4968d^5 + 140\,608d^3 - 259\,584d\,, \\[-.5mm]
0 &=& 2c^3 + 28c^2 - 3cd^2 + 98c - 8d^2 + 104\,.
\end{eqnarray*}
Then\/~$C$
is a geometrically irreducible, nonsingular curve of
genus\/~$1$.
\item[{\rm b) }]
Moreover,~there is strong evidence that, for~generic\/
$(a,b,c,d) \in C(\bbC)$,
the
$K3$~surface
$\smash{X_{(a,b,c,d)}}$
obtained as the minimal desingularisation
of~$\smash{X'_{(a,b,c,d)}}$
is of Picard
rank\/~$16$
and has real multiplication
by\/~$\bbQ(\sqrt{13})$.
\end{abc}%\smallskip

\noindent
{\bf Proof} {\rm of~a).}
i)
Putting~$(a,b,c,d) := (0,-4,-9,-8)$,
the surface from Example~\ref{qw13} is~obtained.\smallskip

\noindent
ii)
This~is easily obtained by a calculation in any computer algebra~system.
%We~used {\tt magma} for this~purpose.
The~curve
$C$
is the result of Strategy~\ref{str}, taking the surface from Example~\ref{qw13} as the starting~point.
\eop\smallskip

\noindent
{\bf Evidence} {\rm for~b).}
For~every prime
$p < 500$
and every
$(a_0,b_0,c_0,d_0) \in C(\bbF_{\!p})$
of the kind that
$X_{(a_0,b_0,c_0,d_0)}$
is a nonsingular surface
over~$\bbF_{\!p}$,
the characteristic
polynomial~$\smash{\chi^\tr_p}$
either splits over
$\smash{\bbQ(\sqrt{13})}$
or
$\smash{\chi^\tr_{p^2}}$
is a square, as predicted by Proposition~\ref{arith_eff}.B.a). Moreover,~for every prime
$p \equiv 2,5,6,7,8,11 \pmod {13}$
up
to~$1000$
and every
$(a_0,b_0,c_0,d_0) \in C(\bbF_{\!p})$,
one has that
$\smash{\#X_{(a_0,b_0,c_0,d_0)}(\bbF_{\!p}) \equiv 1 \pmod p}$,
cf.~Proposition~\ref{arith_eff}.B.b). To~be more precise,
$\smash{\#X'_{(a_0,b_0,c_0,d_0)}(\bbF_{\!p}) = p^2+p+1}$.
\eop
\end{ex}

\begin{rem}
The~genus~$1$
curve~$C$
has
\mbox{$\bbQ$-rational}
points. Taking~any of them as the origin, the Mordell--Weil group
of~$C$
is isomorphic
to~$\bbZ$.
\end{rem}

\begin{rems}[Some details on the experiment]
\label{det_exp}
\begin{iii}
\item
When~running Algorithm~\ref{cup_prod}, we used 14 transcendental cohomology classes, which were represented by tori and obtained as explained in Section~\ref{part_fam}.
For numerical integration, the Gau{\ss}-Legendre method of degree
$30$
[i.e.\
order~$60$]
was applied, using floats of 30~digits.

In~step~\ref{sing_val} of Algorithm~\ref{cup_prod}, we found six singular values within a factor
of~$100$,
while the next one was smaller by nine orders of~magnitude.
In~the basis chosen, the cup product form found
on~$P^\perp$
had coefficients only from
$\smash{\{\pm1, \pm\frac12, 0\}}$,
up to errors less than
$10^{-10}$.
\item\looseness-1
In step~\ref{num_cont} of Strategy~\ref{str}, in the language of \cite{AG}, we applied a predictor--cor\-rec\-tor method. More precisely, we used the Euler predictor~\cite[Section~2.2]{AG}, followed by Newton~corrector steps. We~did not care too much about rounding errors, as we worked with floats of high~precision.
Neither did we implement a step length adaptation, but worked, as simply as possible, with a constant step length. For numerical integration, the Gau{\ss}-Legendre method of degree
$100$
[i.e.\
order~$200$]
was used.

Based~on this, we determined a sample of 101 points
on~$\Pi^{-1}(C) \subset \bbD^4$,
each with a numerical precision of 80~digits. Due to the constant stepsize, these points are essentially equidistant. There are further particularities, caused by the limitations of our approach. First of all, all the points are real. Moreover, they have a limited distance from the starting point, i.e.\ the parameters of the surface from Example~\ref{qw13}, because we are forced to stop sampling when approaching the first singular surface along the path.
\item
Polynomials of degree
$\leq \!3$
in four variables form a vector space of
dimension~$35$.
When~looking for cubic relations between the 101 points found, we ended up with 25 singular values in the range from
$1714$
to~$6.08 \!\cdot\! 10^{-41}$,
the other ten being less
than~$10^{-80}$.
Thus,~the curve sought is contained in an intersection of ten cubics
in~$\Ab^4$.
The~equations given form a Gr\"obner basis for the ideal generated by~them.
\end{iii}
\end{rems}

\section{Explicit description of transcendental cohomology classes---\\Proofs of the main results}

\subsubsection*{Almost $C^1$-maps}
In~the application below, it turns out to be convenient to work with the following technical condition on maps between smooth manifolds.

\begin{defi}
\label{almost_C1}
A continuous map
$\varphi\colon S \to X$
between two smooth manifolds
is called {\em
almost\/~$C^1$},
if there exists a auxiliary\/
$C^1$-map
$\iota\colon S' \to S$
being a homeomorphism from another smooth manifold
$S'$
that satisfies the following two~conditions.

\begin{iii}
\item
The map
$\varphi \!\circ\! \iota\colon S' \to X$
is~$C^1$.
\item
For a suitable Lebesgue null set
$N \subset S'$,
the restriction
$\varphi |_{S \setminus \iota(N)}\colon S \setminus \iota(N) \to X$
is~$C^1$.
\end{iii}
\end{defi}

\begin{exs}
\begin{iii}
\item
Every
$C^1$-map
between smooth manifolds is
almost~$C^1$.
\item
The map
$\varphi\colon \bbR \to \bbR, x \mapsto \sqrt[3]{x}$,
is
almost~$C^1$.
\end{iii}
\end{exs}

\begin{rems}
\begin{iii}
\item
One~has that
$\iota(N)$
is a Lebesgue null set, according to~\cite[section 16.22, probl\`eme~1.c)]{Di3}. Thus, being almost
$C^1$
implies being 
$C^1$
outside of a null~set.
\item
In~particular, if
$\varphi\colon S \to X$
is almost
$C^1$
and
$\eta$
a
$C^1$
differential form
on~$X$
then the pull-back
$\varphi^* \eta$
is defined as a
\mbox{$C^1$-form}
on~$S$,
outside of a Lebesgue null~set.
\item
Rather~generally, for a differential form
$\xi$
that is defined
on~$S$,
up to a Lebesgue null
set~$\underline{N}$,
we put
$\int_S \xi = \int_{S \setminus \underline{N}} \xi$,
as soon as the integral to the right is~existing. 
\end{iii}
\end{rems}

From now on, as in the sections above,
$X$
always denotes a complex
$K3$~surface.

\begin{prop}\looseness-1
\label{cup_prod_integral}
Let~the cohomology class\/
$c_\varphi \in H^2(X, \bbZ)$
be given by a compact oriented\/
\mbox{$2$-manifold}\/~$S$,
together with an almost\/
$C^1$-map\/~$\varphi\colon S \to X$
(cf.~Definition~\ref{cohclass_geom}). Moreover,~let\/
$w \in H^2(X, \bbC)$
be represented by a closed, smooth\/
\mbox{$2$-form}\/~$\eta$.
Then~the extended cup product pairing may be evaluated as the\/
\mbox{$2$-dimensional}
integral
$$(c_\varphi,w) = \int_S \varphi^*\eta \,.$$

\noindent
{\bf Proof.}
{\em
{\em First step.}
The case that
$\varphi$
is~$C^\infty$.

\noindent
One has
$(c_\varphi,w) = \langle w \cup c_\varphi, z_X\rangle = \langle \varphi^* (w), z_S\rangle$,
according to the definition of the cup product pairing and the first claim of Lemma~\ref{homcohom}.b). Interpreting~the term to the right in the de Rham cohomology theory, one indeed has
$\langle \varphi^* (w), z_S\rangle = \int_S \varphi^*\eta$.\smallskip

\noindent
{\em Second step.}
The case that
$\varphi$
is~$C^1$.

\noindent
There~is a
$C^1$
homotopy
$H\colon S \times [0,1] \to X$
connecting
$\varphi$
with a
$C^\infty$-map~$\smash{\underline\varphi}$,
cf.~the proof of \cite[Proposition~17.8]{BT}. As~homotopic maps induce the same homomorphism on cohomology, one has
$\smash{c_\varphi = c_{\underline\varphi}}$
and, therefore,
$\smash{(c_\varphi,w) = (c_{\underline\varphi},w) = \int_S \underline\varphi^*\eta}$,
according to the first~step.

Moreover, the homotopy formula, cf.~\cite[formule~(24.2.4.2)]{Di9}, shows that
$$\varphi^*\eta - \underline\varphi^*\eta = j_0^* H^*\eta - j_1^* H^*\eta = d(LH^*\eta) + L(dH^*\eta) \,,$$
for~$L\colon E^1(S \times [0,1]) \to E^1(S)$
the operator of integration along the fibre. Here,
$dH^*\eta = 0$,
as the form
$\eta$
is~closed. And~hence
$\smash{\int_S \varphi^*\eta - \int_S \underline\varphi^*\eta = \int_S d(LH^*\eta) = 0}$,
due to Stokes' Theorem~\cite[formule~(24.14.2.1)]{Di9}.\smallskip

\noindent
{\em Third step.}
The general case that
$\varphi$
is
almost~$C^1$.

\noindent
Let~$\iota\colon S' \to S$
be an auxiliary map
for~$\varphi$,
as in Definition~\ref{almost_C1}.
Since~$\iota$
is a ho\-me\-o\-morphism, one has
$c_\varphi = \varphi_!(1) = \varphi_!(\iota_!(1)) = (\varphi \!\circ\! \iota)_!(1) = c_{\varphi \circ \iota}$.
Consequently,
$$(c_\varphi,w) = (c_{\varphi \circ \iota}, w) = \int_{S'} (\varphi \!\circ\! \iota)^* \eta \,,$$
according to the step~before. The~integral to the right is the same as
$$\int_{S' \setminus N} \!(\varphi \!\circ\! \iota)^* \eta = \int_{S' \setminus N} \!\iota^* (\varphi^* \eta) \,.$$
Note~here that
$\varphi$
is~$C^1$
on~$\iota(S' \!\setminus\! N) = S \setminus \iota(N)$.
Since~$\iota$
is a
\mbox{$C^1$-map}
and bijective, the last integral coincides with
$\int_{S \setminus \iota(N)} \varphi^* \eta = \int_S \varphi^* \eta$,
as~required.
}
\eop
\end{prop}

\subsubsection*{The main construction--Technical details}
Our idea for the proof of Theorem~\ref{periods_ints} is to apply Proposition~\ref{cup_prod_integral}. This~means that we have to provide a model of~$\alpha_\Gamma$,
respectively~$\alpha_{\Gamma,b}$,
allowing a lift
to~$X$
that is not only continuous, but
almost~$C^1$.\medskip

\paragraph{{\bf Step~1.}~%
\mbox{$2$-dimensional}
tori from compact\/
\mbox{$1$-manifolds} -- Adding an imaginary~direction}
\leavevmode\smallskip

\noindent
As~before, we extend the
$C^1$-map
$\gamma\colon \bbR \to \bbR^2 \subset \Pb^2(\bbR)$
to a
$C^1$-map
in two variables by~putting
$$\gamma'\colon \bbR^2 \longrightarrow \bbC^2 \subset \Pb^2(\bbC), \qquad (t,u) \mapsto \gamma(t) + \mi ub \,.$$
Then
$\smash{\lim_{t\to\pm\infty} \gamma'(t, u)}$
exists
in~$\Pb^2(\bbC)$,
for every
$u \in \bbR$.
Indeed,~this is obvious in the case of a curve encircling a polygon, when simply
$\smash{\!\!\lim\limits_{t\to\pm\infty} \!\!\gamma'(t, u) = \overline\gamma(\infty) + \mi ub}$.
On~the other hand, for a deformed line, we have
\begin{equation}
\label{constr_Stephan_defl}
\gamma'(t,u) = \gamma_0(t) + (t \!+\! \mi u)b \,,
\end{equation}
$\gamma_0$
being of compact~support. This~shows that, independently
of~$u$,
the limit is, on
$\Pb^2(\bbC)$,
the point on the line at infinity, in the direction of the
vector~$b$.
Thus,~in either case,
$\gamma'$~defines
a continuous
map~$\smash{\underline\gamma'}$
from
$\Pb^1(\bbR) \times \bbR$
to~$\Pb^2(\bbC)$.

Moreover,
$\smash{\lim_{u\to\pm\infty} \underline\gamma'(t, u)}$
exists
in~$\Pb^2(\bbC)$
and is independent
of~$t \in \Pb^1(\bbR)$.
In~all cases, it is the point on the infinite line
of~$\Pb^2(\bbC)$
in the direction of the vector
$\mi b \sim b$.
Therefore,~$\gamma'$
actually provides a continuous~map
\begin{equation}
\label{proj_ext}
\alpha'\colon \Tb = \Pb^1(\bbR) \times \Pb^1(\bbR) \to \Pb^2(\bbC) \,.
\end{equation}

\begin{cau}
Here,~in the case of a deformed~line, for the imaginary
direction~$b$,
we take the real direction of the deformed line near~infinity. In~the case of a curve encircling a polygon, the imaginary direction
$b \in \bbR^2 \setminus \big( \bbR(A_{12}, -A_{11}) \cup \cdots \cup \bbR(A_{62}, -A_{61}) \big)$
is to be specified~later. The~actual choice turns out to be irrelevant, cf.~Theorem~\ref{indep_b}.
\end{cau}

\begin{constr}[An unusual differential structure on the
\mbox{$2$-torus}]
We~equip the topological space
$\Pb^1(\bbR) \!\times\! \Pb^1(\bbR)$
with the structure of a smooth
\mbox{$2$-manifold},
in a way that fits our purposes. On~the open subset
$\Pb^1(\bbR) \!\times\! \Pb^1(\bbR) \setminus \{(\infty,\infty)\}$,
we take the natural
$C^\infty$
differential structure, given as an open submanifold of the product. We~extend this structure to the whole
of~$\Pb^1(\bbR) \!\times\! \Pb^1(\bbR)$,
by adding the chart
\begin{eqnarray*}
\varphi_{\infty,\infty}\colon \hspace{3cm} O &\longrightarrow& \bbR^2 \\
((t_0:t_1), (u_0:u_1)) &\mapsto& \textstyle
\left\{
\begin{array}{cl}
(0,0) & {\rm ~if~} (t_0,u_0) = (0,0) \,, \\
\frac{({\frac{t_0}{t_1}}, {\frac{u_0}{u_1}})}{\sqrt[3]{(\frac{t_0}{t_1})^2 + (\frac{u_0}{u_1})^2}} & {\rm ~otherwise} \,,
\end{array}
\right.
\end{eqnarray*}
for
$O := \Pb^1(\bbR) \!\times\! \Pb^1(\bbR) \!\setminus\! (\Pb^1(\bbR) \!\times\! \{0\} \cup \{0\} \!\times\! \Pb^1(\bbR)) \subset \Pb^1(\bbR) \!\times\! \Pb^1(\bbR)$.
It~is obvious that the chart
$\varphi_{\infty,\infty}$
is compatible with those on
$$\Pb^1(\bbR) \!\times\! \Pb^1(\bbR) \setminus \{(\infty,\infty)\} \,,$$
so that altogether they form a
$C^\infty$-atlas
for~$\Pb^1(\bbR) \!\times\! \Pb^1(\bbR)$.

We~denote the
\mbox{$C^\infty$-manifold}
obtained in this way
by~$\Tb'$.
The~identity
$\id\colon\Tb' \to \Tb$
is then a
\mbox{$C^\infty$-map}
and a homeomorphism. The~restriction
$$\id\colon \Tb' \setminus \{(\infty,\infty)\} \to \Tb \setminus \{(\infty,\infty)\}$$
is a diffeomorphism.
\end{constr}

\begin{lem}
\label{C1_toP2_toro}
The map
$$\alpha'\colon \Tb \longrightarrow \Pb^2(\bbC) \,,$$
as defined in~(\ref{proj_ext}), is
almost\/~$C^1$.
In~the case of a curve encircling a polygon,
$\alpha'$
is\/~$C^1$.
Otherwise,
$\alpha'$~is\/
$C^1$
except at the
point\/~$(\infty,\infty)$.\smallskip

\noindent
{\bf Proof.}
{\em
{\em First case.}
$\Gamma$
is a deformed line.

\noindent
Then,~on
$\bbR^2 \subset \Tb$,
the map
$\alpha'$
is given by the formula
\begin{eqnarray*}
((t_0\!:\!t_1), \!(u_0\!:\!u_1)) &\mapsto& \textstyle \gamma_0(\frac{t_1}{t_0}) + (\frac{t_1}{t_0} \!+\! \mi \frac{u_1}{u_0})b \\
&=& \textstyle \big( 1 : \gamma_{0,1}(\frac{t_1}{t_0}) \!+\! (\frac{t_1}{t_0} \!+\! \mi \frac{u_1}{u_0})b_1 : \gamma_{0,2}(\frac{t_1}{t_0}) \!+\! (\frac{t_1}{t_0} \!+\! \mi \frac{u_1}{u_0})b_2 \big) \\
&=& \textstyle \big( t_0u_0 : t_0u_0 \gamma_{0,1}(\frac{t_1}{t_0}) \!+\! (t_1u_0 \!+\! \mi t_0u_1)b_1 : \\[-2.5mm]
&&\textstyle \hspace{5.8cm} : t_0u_0 \gamma_{0,2}(\frac{t_1}{t_0}) \!+\! (t_1u_0 \!+\! \mi t_0u_1)b_2 \big) ,
\end{eqnarray*}
which immediately extends to
$\Tb \!\setminus\! \{((0\!:\!1),(0\!:\!1))\} = \Tb \!\setminus\! \{(\infty,\infty)\}$
and defines a
$C^1$-map
on this open~subset.

Thus,~it remains to show that
$\alpha'$
is
almost~$C^1$.
We~choose
$\id\colon\Tb' \to \Tb$
as the auxiliary map, so that it suffices to verify that
$\alpha' \!\circ\! \id\colon \Tb' \to \Pb^2(\bbC)$
is a
\mbox{$C^1$-map}
near~$(\infty,\infty)$.
For~this, it needs to be shown that
$\alpha' \!\circ\! \varphi_{\infty,\infty}^{-1}$
is~$C^1$
in a neighbourhood
of~$(0,0)$.
But~$\varphi_{\infty,\infty}^{-1}$
sends
$(x,y)$
to
$\big((x(x^2+y^2)\!:\!1),(y(x^2+y^2)\!:\!1)\big)$,
so that this map is given~by
\begin{eqnarray*}
(x,y) &\mapsto& \textstyle \big( xy(x^2+y^2)^2 : xy(x^2+y^2)^2 \gamma_{0,1}(\frac1{x(x^2+y^2)}) \!+\! (y \!+\! \mi x)(x^2+y^2)b_1 : \\[-1.5mm]
&& \textstyle \hspace{4.5cm} : xy(x^2+y^2)^2 \gamma_{0,2}(\frac1{x(x^2+y^2)}) \!+\! (y \!+\! \mi x)(x^2+y^2)b_2 \big) \\
&=& \big( xy(x^2+y^2)^2 : (y \!+\! \mi x)(x^2+y^2)b_1 : (y \!+\! \mi x)(x^2+y^2)b_2 \big) \\
&=& \big( xy(y-\mi x) : b_1 : b_2 \big) \,,
\end{eqnarray*}
at least in a neighbourhood of the~origin. Note~here that
$\gamma_0 = (\gamma_{0,1}, \gamma_{0,2})$
has compact support. Finally,~the function defined by
$xy(y-\mi x)$
is clearly continuously differentiable
at~$(0,0)$.\smallskip

\noindent
{\em Second case.}
$\Gamma$
is a curve encircling a polygon.

\noindent
Here,~for
$\alpha'$,
one has a formula that is valid
on~$\Pb^1(\bbR) \times \bbR \subset \Tb$,
\begin{eqnarray*}
((t_0\!:\!t_1), \!(u_0\!:\!u_1)) &\mapsto& \textstyle \overline\gamma(t_0\!:\!t_1) + \mi \frac{u_1}{u_0}b \\
&=& \textstyle \big( 1 : \overline\gamma_1(t_0\!:\!t_1) \!+\! \mi \frac{u_1}{u_0}b_1 : \overline\gamma_2(t_0\!:\!t_1) \!+\! \mi \frac{u_1}{u_0}b_2 \big) \\
&=& \textstyle \big( u_0 : u_0 \overline\gamma_1(t_0\!:\!t_1) \!+\! \mi u_1 b_1 : u_0 \overline\gamma_2(t_0\!:\!t_1) \!+\! \mi u_1 b_2 \big) \,.
\end{eqnarray*}
From~this, one sees that
$\alpha'$
extends as a
$C^1$-map
to~$\Pb^1(\bbR) \times \Pb^1(\bbR)$.
Finally,~let us note that, by our construction, the identity map
$\id\colon \Tb \to \Pb^1(\bbR) \times \Pb^1(\bbR)$
is~$C^1$.%
}%
\eop
\end{lem}

\begin{rem}
According~to its definition, given in (\ref{constr_Stephan_defl}) and~(\ref{proj_ext}), the torus
$\alpha'$
is constant
on~$\Pb^1(\bbR) \times \{\infty\}$.
Contracting~this subset, the
$2$-torus goes over into a
$2$-sphere
with two points~identified. Thus,~one might possibly work with spheroids instead of~tori.

However,~for the sphere, it is seemingly much harder to explicitly write down a differential structure that makes the resulting
map~$C^1$,
at least in the case of a deformed~line. And~this is what we need in view of Proposition~\ref{cup_prod_integral}, in order to write the periods as~integrals.
\end{rem}\medskip

\paragraph{{\bf Step 2.}~%
Lifting to the double cover, i.e.\ to the
$K3$
surface\/~$X$,
except at\/
$\frakt_1, \ldots, \frakt_n$.}

\begin{nota}
We~write
$\frakt_1, \ldots, \frakt_n$
for the points
$(t_1,0), \ldots, (t_n,0)$
on~$\Tb$.
\end{nota}

\begin{lem}
Assume~that\/
$b$
is generic in the sense of Assumption~(\ref{b_generic}) in the case that\/
$\Gamma$
is a curve encircling a~polygon. I.e.,~that
$$b \not\in \bbR(A_{12}, -A_{11}) \cup \ldots \cup \bbR(A_{62}, -A_{61}) \,.$$
Then~the torus\/
$\alpha'$
meets the branch locus\/
$V(l_1 \cdots l_6)$
only at\/
$\frakt_i$,
for\/
$i = 1, \ldots, n$.
I.e.,~only in the three to five double~points.\smallskip

\noindent\looseness-1
{\bf Proof.}
{\em
First~of all,
$\gamma'(t,u) \in V(l_1 \cdots l_6)$
implies
that~$u = 0$.
Indeed,~$(x,y) \in V(l_i)$
yields
$A_{i1} \Im x + A_{i2} \Im y = 0$.
On~the other hand,
$\Im \gamma'(t,u) = ub = (ub_1, ub_2)$,
so that 
$\gamma'(t,u) \in V(l_i)$
is possible only when
$u(A_{i1} b_1 + A_{i2} b_2) = A_{i1} ub_1 + A_{i2} ub_2 = 0$.
%$u(A_{i1} b_1 + A_{i2} b_2) = 0$.
As~$b \not\in \bbR(A_{i2}, -A_{i1})$,
one
has~$A_{i1} b_1 + A_{i2} b_2 \neq 0$
and hence
$u=0$,
as~required. Thus,~the assertion is shown
for~$\alpha' |_{\bbR^2}$.
In~the case of a curve encircling a polygon, the same argument applies to
$t = \infty$,~too.

Moreover,~all points of
$\Tb$
that are not yet covered are mapped under
$\alpha'$
to the point
of~$\Pb^2(\bbC)$
on the line at infinity, in the direction of the
vector~$\mi b \sim b$.
This~can not be the point at infinity of any of the
lines~$V(l_i)$,
for
$i = 1, \ldots, 6$.
In~fact, one has
$b \not\in \bbR(A_{12}, -A_{11}) \cup \ldots \cup \bbR(A_{62}, -A_{61})$
in either case. For~deformed lines, this is due to Assumption~A.\ref{A3}.
}
\eop
\end{lem}

\begin{prop}
\label{lift}
\begin{abc}
\item
Let\/
$\alpha'\colon \Tb \to \Pb^2(\bbC)$
be a torus as~above. In~the case that\/
$\Gamma$
is a curve encircling a polygon, assume\/
$b \in \bbR^2 \setminus \big( \bbR(A_{12}, -A_{11}) \cup \ldots \cup \bbR(A_{62}, -A_{61}) \big)$.
Then the restriction\/
$\alpha' |_{\Tb \setminus \{\frakt_1,\ldots,\frakt_n\}}$
allows a continuous lift\/
$\alpha''$
to the
$K3$~surface\/~$X$.
\item
In~either case,
$\alpha''$
is automatically
almost\/~$C^1$.
\end{abc}\medskip

\noindent
{\bf Proof.}
{\em
a)
{\em First step.}
Reduction to a statement on fundamental groups.

\noindent
Put~$X_0 \subset X$
to be the preimage
of~$\Pb^2(\bbC) \!\setminus\! V(l_1 \cdots l_6)$
under the natural~map.
As~$X_0$
avoids the branch locus and, in particular, the locus that is blown up, the restriction
$\pi_0\colon X_0 \to \Pb^2(\bbC) \!\setminus\! V(l_1 \cdots l_6)$
of~$\pi$
is a covering projection \cite[Chapter~2, Section~1]{Sp}.

The~assertion is that
$\alpha_0 := \alpha' |_{\Tb \setminus \{\frakt_1,\ldots,\frakt_n\}} \colon \Tb \!\setminus\! \{\frakt_1, \ldots, \frakt_n\} \to \Pb^2(\bbC) \!\setminus\! V(l_1 \cdots l_6)$
allows a continuous
lift~$\alpha''$,
as indicated in the diagram~below,\vspace{-2mm}
$$
\xymatrix{
 & X_0 \ar@{->}[d]^{\pi_0} \\
\Tb \!\setminus\! \{\frakt_1, \ldots, \frakt_n\} \ar@{->}[r]^{\alpha_0\;\;\;} \ar@{.>}[ru]^{\;\;\;\;\;\;\alpha''} & \Pb^2(\bbC) \!\setminus\! V(l_1 \cdots l_6) \,.
}\vspace{-2mm}
$$
In~order to establish this claim, according to~\cite[Chapter~2, Theorem~4.5]{Sp}, we have to show that
$(\alpha_0)_{\#} \pi_1(\Tb \!\setminus\! \{\frakt_1, \ldots, \frakt_n\},\cdot) \subseteq (\pi_0)_{\#} \pi_1(X_0,\cdot)$.

For~this, let us note that
$\pi_1(\Tb \!\setminus\! \{\frakt_1, \ldots, \frakt_n\},\cdot)$,
the fundamental group of the
\mbox{$2$-torus}~minus
$n$~points,
is a free group on
$(n+1)$~generators.
A~generating system is provided by the homotopy classes of
$n$~small
loops
$\nu_1, \ldots, \nu_n\colon S^1 \to \Tb \!\setminus\! \{\frakt_1, \ldots, \frakt_n\}$,
each of which encircles exactly one of the points
$\frakt_1, \ldots, \frakt_n$,
together with the homotopy class of the closed path
$\mu_1$
running through
$\Pb^1(\bbR) \times \{\infty\}$
and the homotopy class of the closed path
$\mu_2$
running through
$\{\infty\} \times \Pb^1(\bbR)$.\medskip

\noindent
{\em Second step.}
The small loops.

\noindent
Let~us first show that
$(\alpha_0)_{\#} [\nu_k] \in (\pi_0)_{\#} \pi_1(X_0,\cdot)$,
for~$k = 1, \ldots, n$.
This~claim is clearly equivalent to the liftability of the loop
$\alpha_0 \!\circ\! \nu_k\colon S^1 \to \Pb^2(\bbC) \!\setminus\! V(l_1 \cdots l_6)$
to~$X_0$.
As~the particular covering projection
$\pi_0\colon X_0 \to \Pb^2(\bbC) \!\setminus\! V(l_1 \cdots l_6)$
is given by
$w^2 = l_1 \cdots l_6$,
this just means that the~loop
\begin{align*}
\underline\nu_k = (l_1 \cdots l_6) \!\circ\! \alpha_0 \!\circ\! \nu_k\colon S^1 &\longrightarrow \bbC \!\setminus\! \{0\} \,,\\
\quad z\, &\;\mapsto\; (l_1 \cdots l_6) (\alpha_0(\nu_k(z)))
\end{align*}
lifts under the covering projection
$\bbC \!\setminus\! \{0\} \to \bbC \!\setminus\! \{0\}, w \mapsto w^2$.
And~the latter is well-known to be possible if and only if the winding number
of~$\underline\nu_k$
is~even.

We may assume that
$\nu_k\colon S^1 \to \Tb \!\setminus\! \{\frakt_1, \ldots, \frakt_n\}$
is given by
$z \mapsto (t_k + \varepsilon \Re z, \varepsilon \Im z)$,
for some real number
$\varepsilon > 0$.
Then the composition
$\alpha_0 \!\circ\! \nu_k\colon S^1 \to \Pb^2(\bbC) \!\setminus\! V(l_1 \cdots l_6)$~is
\begin{eqnarray*}
\alpha_0 \!\circ\! \nu_k\colon (\cos \varphi, \sin \varphi) &\mapsto& \gamma(t_k + \varepsilon \cos \varphi) + \varepsilon \mi \sin \varphi \!\cdot\! b \\
&=& \gamma(t_{i_k,j_k} + \varepsilon \cos \varphi) + \varepsilon \mi \sin \varphi \!\cdot\! b \\
&=& x_{i_k,j_k} + \varepsilon \cos \varphi \!\cdot\! b_{i_k,j_k} + \varepsilon \mi \sin \varphi \!\cdot\! b \,.
\end{eqnarray*}
Hence,~for
$a = 1,\ldots,6$,
the triple composition
$l_a \!\circ\! \alpha_0 \!\circ\! \nu_k\colon S^1 \to \bbC \!\setminus\! \{0\}$
is given~by
\begin{eqnarray*}
l_a \!\circ\! \alpha_0 \!\circ\! \nu_k\colon (\cos \varphi, \sin \varphi) &\mapsto& l_a'(x_{i_k,j_k} + \varepsilon \cos \varphi \!\cdot\! b_{i_k,j_k} + \varepsilon \mi \sin \varphi \!\cdot\! b) \\
&=& l_a'(x_{i_k,j_k}) + \widetilde{l}_a(\varepsilon \cos \varphi \!\cdot\! b_{i_k,j_k} + \varepsilon \mi \sin \varphi \!\cdot\! b) \\
&=& l_a'(x_{i_k,j_k}) + \varepsilon \cos \varphi \!\cdot\! \widetilde{l}_a(b_{i_k,j_k}) + \varepsilon \mi \sin \varphi \!\cdot\! \widetilde{l}_a(b) \,.
\end{eqnarray*}
Here,~for
$a \neq i_k, j_k$,
one has that
$x_{i_k,j_k} \not\in V(l_a)$,
i.e.\ that
$l_a'(x_{i_k,j_k}) \neq 0$.
This shows that the winding number of
$l_a \!\circ\! \alpha_0 \!\circ\! \nu_k$
is~$0$,
at least as long as
$\varepsilon$
is sufficiently~small.

On the other hand, for
$a = i_k$
or~$a = j_k$,
clearly
$x_{i_k,j_k} \in V(l_a)$,
and hence
$l_a'(x_{i_k,j_k}) = 0$.
Moreover,~the constants
$\smash{\widetilde{l}_a(b_{i_k,j_k})}$
and
$\smash{\widetilde{l}_a(b)}$
are both~nonzero. Indeed,~for the first, this follows from Assumption~A.\ref{A1}, while, for the second, this is either Assumption~A.\ref{A3} or Assumption~(\ref{b_generic}). Therefore,~the winding number of
$l_a \!\circ\! \alpha_0 \!\circ\! \nu_k$
must be
$1$
or~$(-1)$.
Altogether,~the winding number of
$\underline\nu_k = (l_1 \cdots l_6) \!\circ\! \alpha_0 \!\circ\! \nu_k$
is an element
of~$\{-2,0,2\}$
and, in particular, even, as~claimed.\medskip

\noindent
{\em Third step.}
The large closed paths.

\noindent
One~has that
$\alpha_0 \!\circ\! \mu_1$
is a constant map. In the case of a deformed line,
$\alpha_0 \!\circ\! \mu_2$
is a constant map,~too.
Therefore,~it suffices to consider
$\mu_2$
in the case of a curve encircling a~polygon. One~has to verify that
$(\alpha_0)_{\#} [\mu_2] \in (\pi_0)_{\#} \pi_1(X_0,\cdot)$,
which is equivalent to the liftability of the closed path
$\alpha_0 \!\circ\!\mu_2\colon \Pb^1(\bbR) \cong S^1 \to \Pb^2(\bbC) \!\setminus\! V(l_1 \cdots l_6)$
to~$X_0$.
For~the particular covering projection
$\pi_0\colon X_0 \to \Pb^2(\bbC) \!\setminus\! V(l_1 \cdots l_6)$,
in exactly the same way as above, this just means to show that the winding number~of
\begin{align*}
\underline\mu_2 = (l_1 \cdots l_6) \!\circ\! \alpha_0 \!\circ\! \mu_2\colon \Pb^1(\bbR) &\longrightarrow \bbC \!\setminus\! \{0\} \,,\\
\quad (u_0\!:\!u_1) &\;\mapsto\; (l_1 \cdots l_6) (\alpha_0(\mu_2(u_0\!:\!u_1)))
\end{align*}
is~even.

Assume that
$\mu_2\colon \Pb^1(\bbR) \to \Tb \!\setminus\! \{\frakt_1, \ldots, \frakt_n\}$
is given by
$(u_0\!:\!u_1) \mapsto (\infty, (u_0\!:\!u_1))$,
for~$\infty = (0\!:\!1)$,
as~before. Then
$\alpha_0 \!\circ\! \mu_2\colon \Pb^1(\bbR) \cong S^1 \to \Pb^2(\bbC) \!\setminus\! V(l_1 \cdots l_6)$~is
\begin{eqnarray*}
\alpha_0 \!\circ\! \mu_2\colon (u_0\!:\!u_1) &\mapsto& \textstyle \overline\gamma(\infty) + \mi \frac{u_1}{u_0} b \\
&=& \textstyle \big( \overline\gamma_1(\infty) + \mi \frac{u_1}{u_0} b_1 : \overline\gamma_2(\infty) + \mi \frac{u_1}{u_0} b_2 : 1 \big) \\
&=&\textstyle \big( \overline\gamma_1(\infty) \!\cdot\! u_0 + \mi b_1 \!\cdot\! u_1 : \overline\gamma_2(\infty) \!\cdot\! u_0 + \mi b_2  \!\cdot\! u_1 : u_0 \big) \,.
\end{eqnarray*}
Hence,~for
$a = 1,\ldots,6$,
the triple composition
$l_a \!\circ\! \alpha_0 \!\circ\! \mu_2\colon \Pb^1(\bbR) \to \bbC \!\setminus\! \{0\}$
is given~by
$$l_a \!\circ\! \alpha_0 \!\circ\! \mu_2\colon (u_0\!:\!u_1) \mapsto l_a(\overline\gamma_1(\infty), \overline\gamma_2(\infty), 1) \!\cdot\! u_0 + \mi l_a(b_1, b_2, 0) \!\cdot\! u_1 \,.$$
Here,~both constants,
$l_a(\overline\gamma_1(\infty), \overline\gamma_2(\infty), 1)$
and~$l_a(b_1, b_2, 0)$,
are~nonzero. Indeed,~for the first, this follows from Assumption~A.\ref{A0}, while, for the second, this is just Assumption~(\ref{b_generic}). Thus,~the winding number of
$l_a \!\circ\! \alpha_0 \!\circ\! \mu_2$
is either
$1$
or~$(-1)$.
Altogether,~the winding number of
$\smash{\underline\mu_2 = (l_1 \cdots l_6) \!\circ\! \alpha_0 \!\circ\! \mu_2}$
is even, as~required. Its~absolute value is bounded
by~$6$.\smallskip

\noindent
b)
The final assertion is an immediate consequence of Lemma~\ref{C1_toP2_toro}, together with the fact that
$\pi_0\colon X_0 \to \Pb^2(\bbC) \!\setminus\! V(l_1 \cdots l_6)$
is a local~diffeomorphism.
}
\eop
\end{prop}\medskip

\paragraph{{\bf Step 3.}~%
Widening the holes of the domain\/
$\Tb \setminus \{\frakt_1,\ldots,\frakt_n\}$}
\leavevmode\smallskip

\noindent
Let us first construct a few auxiliary functions.

\begin{iii}
\item[$\bullet$ ]
We put
$$
\psi\colon [0,\infty) \longrightarrow [0,\infty) \,, \quad r \mapsto \left\{\!
\begin{array}{ll}
0                                 & \text{ if } r \leq \frac12 \,, \\
10(r-\frac12)^2 - 12(r-\frac12)^3 & \text{ if } \frac12 \leq r \leq 1 \,, \\
r                                 & \text{ if } 1 \leq r \,.
\end{array}
\right.
$$
As~is easily checked, the
function~$\psi$
is monotonically increasing
and~$C^1$.
\item[$\bullet$ ]
Working~in polar coordinates and
using~$\psi$,
we~define
$$\Psi\colon \bbR^2 \longrightarrow \bbR^2 \,, \quad (r,\varphi) \mapsto (\psi(r),\varphi) \,.$$
This~is a
$C^1$-map
being the identity
outside~$U_1(0)$
and contracting
$\smash{U_\frac12(0)}$
into the~origin. Consequently,~for
$p \in \bbR^2$,
$$\underline\Psi_p\colon \bbR^2 \longrightarrow \bbR^2 \,, \quad x \mapsto p + \Psi(x-p)$$
is~$C^1$,
equal to the identity
outside~$U_1(p)$,
and contracts
$\smash{U_\frac12(p)}$
to~$p$.
\item[$\bullet$ ]
As~$\smash{\underline\Psi_p}$
differs from
$\id$
only within a compact subset
of~$\bbR^2$,
it extends uniquely to a self-map
$$\Psi_p\colon \Tb \to \Tb$$
that is
again~$C^1$,
contracts
$\smash{U_\frac12(p)}$
to~$p$
and is the identity outside
$\smash{U_1(p)}$.
\end{iii}

\begin{constr}[Widening the holes]
For
$\alpha''\colon \Tb \setminus \{\frakt_1,\ldots,\frakt_n\} \to X$,
put
$$\alpha''' := \alpha'' \!\circ\! (\Psi_{\frakt_1} \!\circ\! \cdots \!\circ\! \Psi_{\frakt_n}) |_{\Tb \setminus (\overline{U_{1/2}(\frakt_1)} \cup \ldots \cup \overline{U_{1/2}(\frakt_n)})} \colon \Tb \setminus \big( \overline{U_{\frac12}(\frakt_1)} \cup \ldots \cup \overline{U_{\frac12}(\frakt_n)} \big) \to X \, .$$
Note that the order of the maps
$\Psi_{\frakt_j}$
is irrelevant here, due to Assumption~A.\ref{A2}.
\end{constr}

\begin{prop}
\label{ext_boundary}
The~almost\/
$C^1$-map
$$\smash{\alpha'''\colon \Tb \setminus \big( \overline{U_{\frac12}(\frakt_1)} \cup \ldots \cup \overline{U_{\frac12}(\frakt_n)} \big) \to X}$$
allows an extension to the manifold with boundary\/
$\smash{\Tb \setminus \big( U_{\frac12}(\frakt_1) \cup \ldots \cup U_{\frac12}(\frakt_n) \big)}$
that is again
almost\/~$C^1$.
It~is, in fact, a\/
\mbox{$C^1$-map}
on\/~$\Tb \setminus \{(\infty,\infty)\}$.\smallskip

\noindent
{\bf Proof.}
{\em
By~the constructions above, on an annulus
around~$\frakt_k$
of inner radius
$\frac12$
and outer
radius~$1$,
the projection
$\smash{\pi \circ \alpha'''\colon \Tb \setminus \big( \overline{U_{\frac12}(\frakt_1)} \cup \ldots \cup \overline{U_{\frac12}(\frakt_n)} \big) \to \Pb^2(\bbC)}$
is given~by
$$(t_k + r \cos\varphi, r \sin\varphi) \mapsto x_{i_k,j_k} + \psi(r) \cos\varphi \!\cdot\! b_{i_k,j_k} + \mi \psi(r) \sin\varphi \!\cdot\! b \in \bbC^2 \subset \Pb^2(\bbC) \,.$$
The~exceptional curve of the blowup
$\Bl_{x_{i_k,j_k}}(\Pb^2(\bbC))$
of~$\Pb^2(\bbC)$
in~$x_{i_k,j_k}$
is covered by two open subschemes, both of which being affine
\mbox{$2$-spaces}.
The~projections down to
$\Pb^2(\bbC)$
are given by
$(x',y') \mapsto x_{i_k,j_k} + (x',x'y')$
and
$(x',y') \mapsto x_{i_k,j_k} + (x'y',y')$.
Let~us work in the first chart, the other one being~analogous.

Then,~on the annulus above, the lift of
$\pi \circ \alpha'''$
to~$\smash{\Bl_{x_{i_k,j_k}}(\Pb^2(\bbC))}$
is given~by
\begin{align}
(t_k \!+\! r \cos\varphi, r \sin\varphi)\! &\textstyle\mapsto \!\big( \psi(r) \cos\varphi \!\cdot\! (b_{i_k,j_k})_1 \!+\! \mi \psi(r) \sin\varphi \!\cdot\! b_1, \frac{\psi(r) \cos\varphi \cdot (b_{i_k,j_k})_2 + \mi \psi(r) \sin\varphi \cdot b_2}{\psi(r) \cos\varphi \cdot (b_{i_k,j_k})_1 + \mi \psi(r) \sin\varphi \cdot b_1} \big) \nonumber\\
&\textstyle= \big( \psi(r) \!\cdot\! (\cos\varphi \!\cdot\! (b_{i_k,j_k})_1 + \mi \sin\varphi \!\cdot\! b_1), \frac{\cos\varphi \cdot (b_{i_k,j_k})_2 + \mi \sin\varphi \cdot b_2}{\cos\varphi \cdot (b_{i_k,j_k})_1 + \mi \sin\varphi \cdot b_1} \big) \,.  \label{blow_coord}
\end{align}
Here, of course, the denominator in the second coordinate might vanish for certain values
of~$\varphi$.
For~these values, however, the numerator is nonzero and therefore other chart~applies. Indeed,~by Assumptions~A.\ref{A1}, A.\ref{A3}, as well as~(\ref{b_generic}), the vectors
$b_{i_k,j_k}$
and~$b$
are both real and different from the zero~vector.

For~$U \subset \{\varphi \in [0,2\pi) \mid \cos\varphi \!\cdot\! (b_{i_k,j_k})_1 + \mi \sin\varphi \!\cdot\! b_1 \neq 0\}$
an open interval, the map~(\ref{blow_coord}) clearly extends as a
\mbox{$C^1$-map}
from
$(\frac12, 1) \times U$
to~$[\frac12, 1) \times U$.
In~total, this shows that
$$\pi \circ \alpha''' |_{U_1(\frakt_k) \setminus \overline{U_\frac12(\frakt_k)}}\colon U_1(\frakt_k) \setminus \overline{U_\frac12(\frakt_k)} \longrightarrow \Pb^2(\bbC)$$
allows an extension to
$\smash{U_1(\frakt_k) \setminus U_\frac12(\frakt_k)}$
that lifts as a
\mbox{$C^1$-map}
to~$\smash{\Bl_{x_{i_k,j_k}}(\Pb^2(\bbC))}$.

We~still have to eliminate the
projection~$\pi$,
i.e.\ to lift to the
$K3$~surface~$X$.
For~this, recall that the double cover
$\pi'\colon X \to \Bl_{x_{12}, \ldots, x_{56}}(\Pb^2(\bbC))$
is a local diffeomorphism outside the branch locus. Moreover,~the branch locus
on~$\smash{\Bl_{x_{i_k,j_k}}(\Pb^2(\bbC))}$
consists of the union of six disjoint~curves. Locally~near
$x_{i_k,j_k}$,
it is the union of the two curves, given by
$\smash{\widetilde{l}_{i_k}(1,y') = 0}$
and~$\smash{\widetilde{l}_{j_k}(1,y') = 0}$.
Thus,~to complete the proof, it suffices to show that the image
of~$\smash{U_\frac12(\frakt_k)}$
does not meet either of the two~curves.

For~this, let us note first that, by our assumptions,
$\smash{\widetilde{l}_{i_k}(1,y') = 0}$
or
$\smash{\widetilde{l}_{j_k}(1,y') = 0}$
is possible only when
$y'$
is~real. Moreover,~in the generic case that the vectors
$b_{i_k,j_k}$
and~$b$
are linearly~independent,
$$\textstyle \frac{\cos\varphi \cdot (b_{i_k,j_k})_2 + \mi \sin\varphi \cdot b_2}{\cos\varphi \cdot (b_{i_k,j_k})_1 + \mi \sin\varphi \cdot b_1} \in \bbR$$
implies that
$\cos\varphi = 0$
or~$\sin\varphi = 0$.
I.e.,~that
$$\textstyle y' = \frac{(b_{i_k,j_k})_2}{(b_{i_k,j_k})_1} \quad\text{or}\quad y' = \frac{b_2}{b_1} \,.$$
But~$\smash{\widetilde{l}_{i_k}(1,\frac{b_2}{b_1}) = 0}$
or
$\smash{\widetilde{l}_{j_k}(1,\frac{b_2}{b_1}) = 0}$
would mean nothing but
$b \in \bbR(A_{i_k,2}, -A_{i_k,1})$
or
$b \in \bbR(A_{j_k,2}, -A_{j_k,1})$,
which is excluded by Assumption A.\ref{A3} or (\ref{b_generic}),~respectively. And
$$\textstyle \widetilde{l}_{i_k}\big( 1,\frac{(b_{i_k,j_k})_2}{(b_{i_k,j_k})_1} \big) = 0 \quad\text{or}\quad \widetilde{l}_{j_k}\big( 1,\frac{(b_{i_k,j_k})_2}{(b_{i_k,j_k})_1} \big) = 0$$
would yield that
$b_{i_k,j_k} \in \bbR(A_{i_k,2}, -A_{i_k,1})$
or
$b_{i_k,j_k} \in \bbR(A_{j_k,2}, -A_{j_k,1})$,
in contradiction with Assumption~A.\ref{A1}.

Finally,~in the exceptional case that the vectors
$b_{i_k,j_k}$
and~$b$
are linearly dependent, one has that
$\smash{\frac{\cos\varphi \cdot (b_{i_k,j_k})_2 + \mi \sin\varphi \cdot b_2}{\cos\varphi \cdot (b_{i_k,j_k})_1 + \mi \sin\varphi \cdot b_1}}$\vspace{.3mm}
is
constantly~$\smash{\frac{b_2}{b_1}}$.
But~both
$\smash{\widetilde{l}_{i_k}(1,\frac{b_2}{b_1}) = 0}$
and
$\smash{\widetilde{l}_{j_k}(1,\frac{b_2}{b_1}) = 0}$
are contradictory to our assumptions, as~above.
}%
\eop
\end{prop}\medskip

\paragraph{{\bf Step 4.}~%
Extending to the whole
of\/~$\Tb$.}

\begin{prop}
\begin{iii}
\item
The~almost\/
$C^1$-map\/
$\smash{\alpha'''\!\colon \Tb \setminus \big( \overline{U_{\frac12}(\frakt_1)} \cup \ldots \cup \overline{U_{\frac12}(\frakt_n)} \big) \to X}$~allows
an extension\/
$\alpha\colon \Tb \to X$
to the whole
of\/~$\Tb$
that is again
almost\/~$C^1$.
\item
The~map\/
$\alpha$
may be chosen in such a way that, for
each\/~$k$,
the image\/
$\smash{\alpha\big( \overline{U_{\frac12}(\frakt_k)} \big)}$
is completely contained in an exceptional~curve.
\end{iii}\medskip

\noindent
{\bf Proof.}
{\em
We~start with the extension
$\smash{\alpha''''\colon \Tb \setminus \big( U_{\frac12}(\frakt_1) \cup \ldots \cup U_{\frac12}(\frakt_n) \big) \to X}$
that is provided by Proposition~\ref{ext_boundary}. Then,
$$\alpha''''(\partial\overline{U_{\frac12}(\frakt_k)}) \subset E_{i_k,j_k} \,,$$
for each
$k \in \{1, \ldots, n\}$,
according to our construction.
Moreover,~for every point
$\smash{x \in \partial\overline{U_{\frac12}(\frakt_k)}}$
and every tangent vector
$v \in T_x \Tb$
that is perpendicular to the boundary and points to the interior of
$\smash{\Tb \setminus \big( U_{\frac12}(\frakt_1) \cup \ldots \cup U_{\frac12}(\frakt_n) \big)}$,
one has
$$T\alpha''''(v) = \vec{0} \in T_{\alpha''''(x)} X \,.$$

Moreover,~as
$\smash{E_{i_k,j_k} \cong \Pb^1(\bbC)}$
is simply-connected, the modification of~\cite[Corollary~17.8.1]{BT} for the
\mbox{$C^1$-case}
yields a
$C^1$-homotopy
$\smash{H\colon \partial\overline{U_{\frac12}(\frakt_k)} \times [0,\frac12] \to E_{i_k,j_k}}$
connecting
$$\smash{H_{\partial\overline{U_{\frac12}(\frakt_k)} \times \{\frac12\}} = \alpha'''' |_{\partial\overline{U_{\frac12}(\frakt_k)}}}$$
with a constant~map. Applying~a suitable monotonously increasing, bijective
$C^1$-self-map
of~$[0,\frac12]$,
if necessary, one may assume that
$$\smash{TH(v) = \vec{0} \in T_{H(x,\frac12)} E_{i_k,j_k}}
\quad\text{and}\quad
\smash{TH(x,v) = \vec{0} \in T_{H(x,0)} E_{i_k,j_k}}$$
for every tangent vector
$\smash{v \in T_{(x,\frac12)}(\partial\overline{U_{\frac12}(\frakt_k)} \times [0,\frac12])}$
or
$\smash{v \in T_{(x,0)}(\partial\overline{U_{\frac12}(\frakt_k)} \times [0,\frac12])}$,
respectively, that is perpendicular to the boundary and points to the~inside.

Putting,~finally,
$$\alpha(tx) := H(x,t)$$
for
$\smash{x \in \partial\overline{U_{\frac12}(\frakt_k)}}$
and~$t \in [0,\frac12]$,
and gluing the maps together, one finds an extension
$\alpha$
of~$\alpha'''$
with all the properties~required.
\eop
}
\end{prop}\medskip

\paragraph{{\bf Step 5.}~%
Projecting down
to\/~$X'$.}
\leavevmode\smallskip

\noindent
We~obtain the torus
$\alpha_\Gamma\colon \Tb \to X'$
or
$\alpha_{\Gamma,b}\colon \Tb \to X'$,
respectively, as the composition
$\bl \circ \alpha$.\medskip

\begin{rem}
Modifying Definition~\ref{almost_C1} in an obvious way, one defines the concept of an {\em
almost~$C^\infty$-map}.
Starting~with a
$C^\infty$-parametrisation of the
\mbox{$1$-manifold}
$\Gamma$,
it is certainly possible to work out the constructions just presented entirely within the setup of
almost~$C^\infty$-maps.
We~prefer the use of
almost~$C^1$-maps,
as this approach keeps Steps~3 and~4 slightly~simpler.
\end{rem}

\subsubsection*{Periods as improper integrals}

\begin{proo_periods}
According~to our construction, we may suppose that the lift
$\alpha\colon \Tb \to X$
of the torus
$\alpha_\Gamma$
or
$\alpha_{\Gamma,b}$
considered is an almost
$C^1$~map.
Therefore,~Proposition~\ref{cup_prod_integral} applies and shows that
$(c_\alpha, [\omega]) = \int_\Tb \alpha^* \omega$.

Moreover,~the exceptional curves
$E_{i_1,j_1}, \ldots, E_{i_n,j_n}$
are holomorphic curves
on~$X$.
Hence,~as
$\omega$
is a
\mbox{$(2,0)$-form},
one has that the restrictions
$\omega |_{E_{i_k,j_k}}$
are equal to the null~form. In~particular,
$\smash{\int_{E_{i_k,j_k}} \omega = 0}$.
Pulling~back, one finds that
$\smash{\int_{U_\frac12(\frakt_k)} \alpha^* \omega = 0}$.
Since,~in addition,
$\partial U_\frac12(\frakt_k)$
is a Lebesgue null set, this~yields
$$(c_\alpha, [\omega]) = \!\!\int\limits_{\Tb \setminus (\overline{U_{1/2}(\frakt_1)} \cup \ldots \cup \overline{U_{1/2}(\frakt_n)})} \!\!\!\!\!\!\!\!\!\!\!\!\!\!\!\!\alpha^* \omega = \int\limits_{\Tb \setminus (\overline{U_{1/2}(\frakt_1)} \cup \ldots \cup \overline{U_{1/2}(\frakt_n)})} \!\!\!\!\!\!\!\!\!\!\!\!\!\!\!\!\alpha''' {}^* \omega = \int\limits_{\Tb \setminus \{\frakt_1,\ldots,\frakt_n\}} \!\!\!\!\!\!\!\!\alpha''{}^* \omega = \int\limits_{\Tb \setminus \{\frakt_1,\ldots,\frakt_n\}} \!\!\!\!\!\!\!\!\alpha''{}^* \omega' \,.$$
Furthermore,~noticing~that
$\bbR^2 \setminus \{\frakt_1,\ldots,\frakt_n\}$
differs from
$\Tb \setminus \{\frakt_1,\ldots,\frakt_n\}$
by a Lebesgue null set, we conclude~that
$$(c_\alpha, [\omega]) = \!\!\!\int\limits_{\bbR^2 \setminus \{\frakt_1,\ldots,\frakt_n\}} \!\!\!\!\!\!\!\!(\alpha'' |_{\bbR^2 \setminus \{\frakt_1,\ldots,\frakt_n\}} ){}^* \omega' = \int\limits_{\bbR^2} (\alpha'' |_{\bbR^2 \setminus \{\frakt_1,\ldots,\frakt_n\}} ){}^* \omega' \,.$$
Finally,~a direct calculation of the differential
\mbox{$2$-form}
$(\alpha'' |_{\bbR^2 \setminus \{\frakt_1,\ldots,\frakt_n\}} ){}^* \omega'$
in terms of the coordinates
$t$
and~$u$,
making use of formulae~(\ref{global_form}) and~(\ref{constr_Stephan}), establishes the~claim.%
\eop
\end{proo_periods}

\subsubsection*{Independence of the choice
of\/~$b$
in the case of a curve encircling a polygon}

\begin{proo_indep}
{\em First step.}
Preparations.

\noindent
By~Corollary~\ref{lift_coh_modP}, it suffices to show that
$\smash{\alpha_{\Gamma,\underline{b}}\colon \Tb \to X'}$
is homotopic either to
$\smash{\alpha_{\Gamma,\underline{\underline{b}}}}$
or to
$\smash{\alpha_{\Gamma,\widetilde{\underline{\underline{b}}}}\colon \Tb \to X'}$.
Moreover,~in doing so, it is clearly sufficient to assume that
$\smash{\underline{b}}$~and~$\smash{\underline{\underline{b}}}$
lie in neighbouring sectors of the set
$\bbR^2 \setminus \big( \bbR(A_{12}, -A_{11}) \cup \ldots \cup \bbR(A_{62}, -A_{61}) \big)$
of admissible~vectors.

We~let
$s_0 \in (0,1)$
be the unique real number such that
$\smash{b_0 := \underline{b} + s_0(\underline{\underline{b}} - \underline{b})}$
is not an admissible vector, cf.~Figure~\ref{adm_vec}. Furthermore,~we denote by
$i_0 \in \{1, \ldots, 6\}$
the index, for which
$\smash{b_0 \in \bbR(A_{i_0,2}, -A_{i_0,1})}$.

\begin{figure}[H]
\centerline{
\begin{overpic}[scale=0.8]{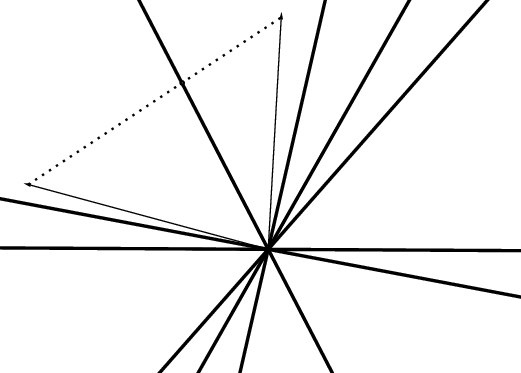}
 \put(110,138){$\underline{b}$}
 \put(4,75){$\underline{\underline{b}}$}
 \put(58,113){$b_0$}
\end{overpic}}\vspace{-.3cm}
\caption{The set of admissible vectors}
\label{adm_vec}\vspace{-.2cm}
\end{figure}

Then,~at first, the map
\begin{align*}
\underline{H}\colon \bbR^2 \times [0,1] &\longrightarrow \bbC^2 \subset \Pb^2(\bbC) \,, \\
(t,u,s) &\;\mapsto\; \gamma(t) + \mi u[\underline{b} + s(\underline{\underline{b}} - \underline{b})] \,,
\end{align*}
clearly allows an extension to a continuous~map
$$H\colon \Tb \times [0,1] \longrightarrow \Pb^2(\bbC) \,,$$
which is a homotopy connecting
$\smash{\alpha'_{\underline{b}}}$
with~$\smash{\alpha'_{\underline{\underline{b}}}}$.
Moreover,~one readily sees that
$H^{-1}(V(l_1 \cdots l_6)) = B$,~for
$B \subset \Tb \times [0,1]$
the~subset
$$B := \big(\{\frakt_1\} \!\times\! [0,1]\big) \cup \ldots \cup \big(\{\frakt_n\} \!\times\! [0,1]\big) \cup \!\!\!\!\bigcup_{\atop{k=1,\ldots,n}{l_{i_0}(\gamma(t_k)) = 0}} \!\!\!\!\!\!\!\big(\{t_k\} \!\times\! \Pb^1(\bbR) \!\times\! \{s_0\}\big) \,.$$
I.e.,~the map
$H$
hits the branch locus
$V(l_1 \cdots l_6)$
exactly for the arguments
$(t,u,s) \in B$.

There~are now two cases to be distinguished. Either,
$V(l_{i_0})$
is one of the edges of the polygon encircled
by~$\Gamma$.
Then
$l_{i_0}(\gamma(t_k)) = 0$
happens exactly two times, say
for~$k_1$
and~$k_2$.
In~this case,
$\smash{\gamma(t_{k_1}) = x_{i_{k_1},j_{k_1}}}$
and~$\smash{\gamma(t_{k_2}) = x_{i_{k_2},j_{k_2}}}$
are the two vertices of the
edge~$V(l_{i_0})$.
Which~means that
$i_0$
is one of the indices
$\smash{i_{k_1}}$
and~$\smash{j_{k_1}}$,
as well as one of the indices
$\smash{i_{k_2}}$
and~$\smash{j_{k_2}}$.
We~take the ordering such that
$t_{k_1} < t_{k_2}$.
Moreover,~we assume, without loss of generality, that the parametrisation
$\gamma$
is of the kind that
$\gamma((t_{k_1}, t_{k_2}))$
is the direct path
on~$\Gamma$
from
$\smash{x_{i_{k_1},j_{k_1}}}$
to~$\smash{x_{i_{k_2},j_{k_2}}}$,
not meeting any other branch point, cf.~Figure~\ref{segment}~below.

Or,~$V(l_{i_0})$
is not an edge of the polygon.
Then~$l_{i_0}(\gamma(t_k)) = 0$
does not happen, for
any~$k$.\medskip

\begin{figure}[H]
\centerline{
\begin{overpic}[scale=0.7]{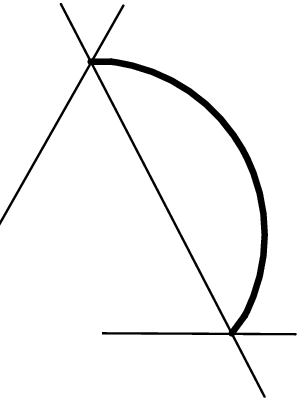}
 \put(86,-9){$V(l_{i_0})$}
 \put(44,28){$\gamma(t_{k_1})$}
 \put(-1,114){$\gamma(t_{k_2})$}
 \put(86,80){$\gamma((t_{k_1}, t_{k_2})) \subset \Gamma$}
\end{overpic}}\vspace{-.1cm}
\caption{The direct path along $\Gamma$ from $\gamma(t_{k_1})$ to~$\gamma(t_{k_2})$}
\label{segment}\vspace{-.2cm}
\end{figure}

\noindent
{\em Second step.
Lifting to the double
cover~$X'$.}

\noindent
As~usual, the continuity of a lift
of~$H$
to~$X'$,
$$
\xymatrix{
 &&& X_0 \hsmash{\;\subset X'} \ar@{->}[d]^{\pi_0} \\
(\Tb \!\times\! [0,1]) \!\setminus\! B \ar@{->}[rrr]^{H|_{(\Tb \times [0,1]) \setminus B}} \ar@{.>}[rrru]^{H'} &&& \Pb^2(\bbC) \!\setminus\! V(l_1 \cdots l_6) \,,
}\vspace{0mm}
$$
is governed by the image of the fundamental group
$\pi_1((\Tb \!\times\! [0,1]) \!\setminus\! B, \cdot)$
under
$(H|_{(\Tb \times [0,1]) \setminus B})_\#$.

By~\cite[Chapter~2, Theorem~4.5]{Sp}, there would be a continuous lift if and only~if
$$(H|_{(\Tb \times [0,1]) \setminus B})_\# \pi_1((\Tb \times [0,1]) \!\setminus\! B, \cdot) \subseteq (\pi_0)_\# \pi_1(X_0, \cdot) \,.$$
So~let us consider the fundamental group
$\pi_1((\Tb \!\times\! [0,1]) \!\setminus\! B, \cdot)$.
It~is generated by the homotopy classes
$[\nu_1], \ldots, [\nu_n]$
of small loops encircling
$\{\frakt_1\} \!\times\! [0,1]$, \ldots, $\{\frakt_n\} \!\times\! [0,1]$,
respectively, together with the homotopy class
$[\mu_1]$
of a closed path running through
$\Pb^1(\bbR) \times \{\infty\} \times \{0\}$,
and the homotopy class
$[\mu_2]$
of the closed path running through
$\{\infty\} \times \Pb^1(\bbR) \times \{0\}$.
In~the more interesting case that
$V(l_{i_0})$
is one of the edges of the polygon encircled
by~$\Gamma$,
two further generators need to be taken into consideration, namely the homotopy classes
$[\lambda_{k_1}]$
and
$[\lambda_{k_2}]$
of small loops
$\lambda_{k_1}, \lambda_{k_2}\colon S^1 \to (\Tb \!\times\! [0,1]) \!\setminus\! B$~encircling
$$\{t_{k_1}\} \!\times\! \Pb^1(\bbR) \!\times\! \{s_0\} \quad\text{and}\quad \{t_{k_2}\} \!\times\! \Pb^1(\bbR) \!\times\! \{s_0\} \,,$$
respectively.

One~indeed has
$(H|_{(\Tb \times [0,1]) \setminus B})_\# ([\nu_k]) \subseteq (\pi_0)_\# \pi_1(X_0, \cdot)$,
for~$k = 1,\ldots,n$,
as well as
$(H|_{(\Tb \times [0,1]) \setminus B})_\# ([\mu_1]) \subseteq (\pi_0)_\# \pi_1(X_0, \cdot)$
and
$(H|_{(\Tb \times [0,1]) \setminus B})_\# ([\mu_2]) \subseteq (\pi_0)_\# \pi_1(X_0, \cdot)$.
The~calculations showing this are essentially the same as those performed in the second and third steps of the proof of Proposition~\ref{lift}.a). Thus,~there is no need to repeat them~here.

As~a conclusion, one finds that
$H$
allows a continuous
lift~$H'\colon \Tb \times [0,1] \to X'$
in the case when
$V(l_{i_0})$
is not an edge of the~polygon.
In~particular, in this case,
$\smash{\alpha_{\Gamma,\underline{b}}}$
is homotopic either to
$\smash{\alpha_{\Gamma,\underline{\underline{b}}}}$
or to
$\smash{\alpha_{\Gamma,\widetilde{\underline{\underline{b}}}}}$,
as~claimed.\medskip

\noindent
{\em Third step.
A discontinuous lift.}

\noindent
On~the other~hand,
$(H|_{(\Tb \times [0,1]) \setminus B})_\# ([\lambda_k]) \subseteq (\pi_0)_\# \pi_1(X_0, \cdot)$
is {\em false\/} for both indices,
$k = k_1$
and~$k_2$.
In~order to see this, let us restrict considerations
to~$[\lambda_{k_1}]$,
the arguments for the other one being exactly the~same. Similarly~to the above, we have to show that the winding number of
$\underline\lambda_{k_1} := (l_1 \cdots l_6) \!\circ\! H|_{(\Tb \times [0,1]) \setminus B} \!\circ\! \lambda_{k_1} \colon S^1 \to \bbC \!\setminus\! \{0\}$
is~odd.

Choose a real number
$u_0 \neq 0$
and let
$\lambda_{k_1}\colon S^1 \to (\Tb \times [0,1]) \!\setminus\! B$
be given by
$$(\cos \varphi, \sin \varphi) \mapsto (t_{k_1} + \varepsilon \cos \varphi, u_0, s_0 + \varepsilon \sin \varphi) \,,$$
for some real number
$\varepsilon > 0$.
Then
$H|_{(\Tb \times [0,1]) \setminus B} \!\circ\! \lambda_{k_1}\colon S^1 \to \Pb^2(\bbC) \!\setminus\! V(l_1 \cdots l_6)$~is
\begin{eqnarray*}
H|_{(\Tb \times [0,1]) \setminus B} \!\circ\! \lambda_{k_1}\colon (\cos \varphi, \sin \varphi) &\mapsto& \gamma(t_k + \varepsilon \cos \varphi) + \mi u_0 [\underline{b} + (s_0 + \varepsilon \sin \varphi)(\underline{\underline{b}} - \underline{b})] \\
 &=& \gamma(t_k + \varepsilon \cos \varphi) + \mi u_0[b_0 + \varepsilon \sin \varphi(\underline{\underline{b}} - \underline{b})] \\
 &=& x_{i_{k_1}, j_{k_1}} \!\!+\! \mi u_0 \!\cdot\! b_0 \!+\! \varepsilon \cos \varphi \!\cdot\! b_{i_{k_1}, j_{k_1}} \!\!+\! \varepsilon \mi \sin \varphi \!\cdot\! u_0 (\underline{\underline{b}} - \underline{b}) \,.
\end{eqnarray*}
Hence,~for
$a = 1,\ldots,6$,
the triple composition
$l_a \!\circ\! H|_{(\Tb \times [0,1]) \setminus B} \!\circ\! \lambda_{k_1}\colon S^1 \to \bbC \!\setminus\! \{0\}$
is given~by
\begin{align*}
\smash{l_a \!\circ\! H|_{(\Tb \times [0,1]) \setminus B}} \!\circ\! \lambda_{k_1}&\colon (\cos \varphi, \sin \varphi) \mapsto \\[-1mm]
  & \smash{l_a(x_{i_{k_1}, j_{k_1}} + \mi u_0 \!\cdot\! b_0 + \varepsilon \cos \varphi \!\cdot\! b_{i_{k_1}, j_{k_1}} + \varepsilon \mi \sin \varphi \!\cdot\! u_0(\underline{\underline{b}} - \underline{b}))} \\
={} & \smash{l_a'(x_{i_{k_1}, j_{k_1}}) \!+\! \mi u_0 \!\cdot\! \widetilde{l}_a(b_0) + \varepsilon \cos \varphi \!\cdot\! \widetilde{l}_a(b_{i_{k_1}, j_{k_1}}) + \varepsilon \mi \sin \varphi \!\cdot\! u_0 \widetilde{l}_a(\underline{\underline{b}} - \underline{b}) \,. \!}
\end{align*}
Here,~for
$a \neq i_0$,
the real number
$\widetilde{l}_a(b_0)$
is~nonzero.
As~$l_a'(x_{i_{k_1}, j_{k_1}}) \in \bbR$,
this means that the winding number of
$\smash{l_a \!\circ\! H|_{(\Tb \times [0,1]) \setminus B} \!\circ\! \lambda_{k_1}}$
is~$0$,
at least as long as
$\varepsilon$
is sufficiently~small.

On~the other hand,
$l_{i_0}'(x_{i_{k_1}, j_{k_1}}) = 0$
and
$\smash{\widetilde{l}_{i_0}(b_0) = 0}$.
Moreover,~the constants
$\smash{\widetilde{l}_{i_0}(b_{i_{k_1}, j_{k_1}})}$
and
$\smash{\widetilde{l}_{i_0}(\underline{\underline{b}} - \underline{b}))}$
are both~nonzero. This~is due to Assumption~A.\ref{A1} for the first and obvious for the~second. Therefore,~the winding number of
$\smash{l_{i_0} \!\circ\! H|_{(\Tb \times [0,1]) \setminus B} \!\circ\! \lambda_{k_1}}$
must be
$1$
or~$(-1)$.
Altogether,~the winding number of the triple composition
$\smash{\underline\lambda_{k_1} = (l_1 \cdots l_6) \!\circ\! H|_{(\Tb \times [0,1]) \setminus B} \!\circ\! \lambda_{k_1}}$
is
$1$
or~$(-1)$
and, in particular, odd, as~claimed.

As~a conclusion, one finds that
$H\colon \Tb \times [0,1] \to \Pb^2(\bbC)$
does not allow any continuous lift
to~$X'$.
There~is, however, a discontinuous~lift
$$H'\colon \Tb \times [0,1] \longrightarrow X'$$
that is continuous
on~$(\Tb \!\times\! [0,1]) \!\setminus\! \big( (t_{k_1}, t_{k_2}) \!\times\! \Pb^1(\bbR) \!\times\! \{s_0\} \big)$.
The~discontinuity is of the kind that,
for~$(t,u) \in \Tb$
such that
$t \in (t_{k_1}, t_{k_2})$,
the limit
$\smash{\alpha_{\Gamma,b_0}^+ (t,u) := \!\lim\limits_{s\to s_0+0}\!\!\! H'(t,u,s)}$
differs from
$\smash{\alpha_{\Gamma,b_0}^- (t,u) := \lim_{s\to s_0-0} H'(t,u,s)}$
by the
involution~$\zeta$.

In~particular,
$H'$
provides us with two homotopies, one connecting
$\smash{\alpha_{\Gamma,\underline{b}}}$
with
$\smash{\alpha_{\Gamma,b_0}^-}$
and another, connecting
$\smash{\alpha_{\Gamma,\underline{\underline{b}}}}$
or
$\smash{\alpha_{\Gamma,\widetilde{\underline{\underline{b}}}}}$
with~$\smash{\alpha_{\Gamma,b_0}^+}$.\medskip\pagebreak[3]

\noindent
{\em Fourth step.
Straightening the\/
\mbox{$1$-manifold}~$\Gamma$.}

\noindent
It~is thus our remaining task to show that
$\smash{\alpha_{\Gamma,b_0}^+}$
and
$\smash{\alpha_{\Gamma,b_0}^-\colon \Tb \to X'}$
are homotopic to each~other. For~this, consider the~following,
\begin{align*}
\underline{D}\colon \Pb^1(\bbR) \!\times\! \bbR \!\times\! [0,1] &\longrightarrow  \bbC^2 \subset \Pb^2(\bbC) \,, \\
(t,u,s) &\;\mapsto\;
\left\{
\begin{array}{ll}
\textstyle s \big( \frac{t_{k_2}-t}{t_{k_2}-t_{k_1}} \gamma(t_{k_1}) + \frac{t-t_{k_1}}{t_{k_2}-t_{k_1}} \gamma(t_{k_1}) \big) \\
\hspace{3.1cm} {} + (1-s) \gamma(t) + \mi u \!\cdot\! b_0 & \;\text{ if } t \in (t_{k_1}, t_{k_2}) \,, \!\!\\
\gamma(t) + \mi u \!\cdot\! b_0 & \;\text{ otherwise} \,.
\end{array}
\right.
\end{align*}
The~map
$\underline{D}$
clearly allows an extension to a continuous map
$D\colon \Tb \times [0,1] \to \Pb^2(\bbC)$,
which is a homotopy between the maps
$\alpha'_\Gamma$
and~$\alpha'_{\Gamma'}$,
for
$\Gamma'$
a straightening
of~$\Gamma$,
as indicated in the figure~below.\vspace{-3mm}

\begin{figure}[H]
\centerline{
\begin{overpic}[scale=0.7]{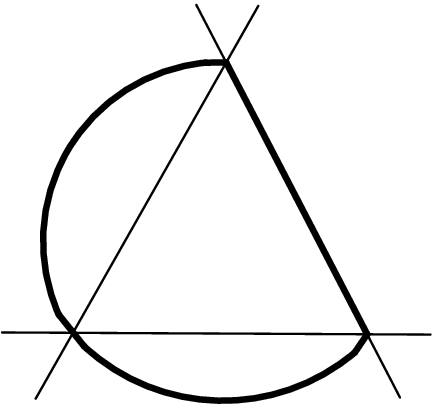}
 \put(126,-9){$V(l_{i_0})$}
 \put(124,28){$\gamma(t_{k_1})$}
 \put(79,114){$\gamma(t_{k_2})$}
\end{overpic}}\vspace{-.1cm}
\caption{A straightening of $\Gamma$}
\vspace{-.2cm}
\end{figure}
Moreover,~one has
$D^{-1}(V(l_1 \cdots l_6)) = B$,
for
$$B := \bigcup_{\atop{k=1,\ldots,n}{k \neq k_1,k_2}} \!\!\!\!\big( \{\frakt_k\} \!\times\! [0,1] \big) \,\cup \!\!\bigcup_{k = k_1,k_2} \!\!\!\! \big( \{t_k\} \!\times\! \Pb^1(\bbR) \!\times\! [0,1] \big) \cup \big( (t_{k_1}, t_{k_2}) \!\times\! \Pb^1(\bbR) \!\times\! \{1\} \big) \,.$$
\begin{figure}[H]
\centerline{
\begin{overpic}[scale=0.7]{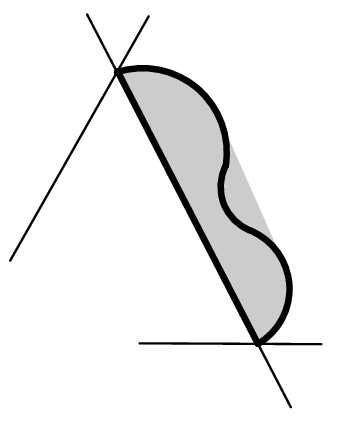}
 \put(89,-3){$V(l_{i_0})$}
 \put(60,17){$\gamma(t_{k_1})$}
 \put(7,117){$\gamma(t_{k_2})$}
 \put(76,105){$\gamma((t_{k_1}, t_{k_2}))$}
\end{overpic}}\vspace{-.1cm}
\caption{The convex hull of the line segment and the direct path
$\gamma((t_{k_1}, t_{k_2}))$}
\label{conv_hull}\vspace{-.2cm}
\end{figure}
Indeed,~the inclusion
``$\supseteq$''
is obvious. In~order to verify the reverse
inclusion~``$\subseteq$'',
we first observe that the convex hull of the union of the line segment
$[\gamma(t_{k_1}), \gamma(t_{k_2})]$
with the direct path
$\gamma((t_{k_1}, t_{k_2}))$
does not contain any point
of~$V(l_1 \cdots l_6)$
in its~interior.
In~fact, a
line~$V(l_a)$
through such a point would necessarily intersect
$\gamma((t_{k_1}, t_{k_2}))$
somewhere, in contradiction to the assumption that this is a direct path from one vertex of the encircled polygon to the next, cf.\ Figure~\ref{conv_hull},~above.

Now~let
$(t,u,s) \in \Tb \!\times\! [0,1]$
be any point such that
$D(t,u,s) \in V(l_1 \cdots l_6)$.
I.e., such that
$l_a(D(t,u,s)) = 0$
for
some~$a \in \{1,\ldots,6\}$.
In~the particular case
that~$t \in (t_{k_1}, t_{k_2})$,
the reasoning just given shows that
\begin{equation}
\label{real_D}
0 = \Re l_a(D(t,u,s)) = l_a(\Re D(t,u,s))
\end{equation}
is possible only
for~$s=1$,
as~claimed.

On~the other hand, for
$\smash{t \not\in (t_{k_1}, t_{k_2})}$,
equation (\ref{real_D}) means that
$0 = l_a(\gamma(t))$,
which, according to the general assumptions on the
\mbox{$1$-manifold}~$\Gamma$,
together with Notation~\ref{ramif_pts}.ii), may happen only
for~$t = t_1, \ldots, t_n$.
Moreover,~in this case,
$$0 = \Im l_a(D(t,u,s)) = l_a(\Im D(t,u,s)) = u \!\cdot\! \widetilde{l}_a(b_0) \,,$$
which is true only for
$u = 0$
or~$a = i_0$.
But~$l_{i_0}(\gamma(t)) = 0$
happens only for
$\smash{t = t_{k_1}}$
and~$\smash{t = t_{k_2}}$,
which completes the argument for the 
inclusion~``$\subseteq$''.

The~domain
$(\Tb \!\times\! [0,1]) \!\setminus\! B$
turns out to be disconnected into two~components. We~construct a continuous lift
$D'\colon \Tb \!\times\! [0,1] \to X'$
of~$D$
as~follows.
For~$t \not\in (t_{k_1}, t_{k_2})$,
we
put~$\smash{D'(t,u,s) := \alpha_{\Gamma,b_0}^+(t,u)}$.
For~$(t,u,s) \in B$,
there is a unique lift
of~$D(t,u,s)$,
anyway. Finally,
$\smash{\big( (t_{k_1}, t_{k_2}) \!\times\! \Pb^1(\bbR) \!\times\! [0,1] \big) \!\setminus\! B}$
is simply connected. Thus,~there are two continuous lifts, interchanged by the
involution~$\zeta$,
and we may take the one that agrees
with~$\smash{\alpha_{\Gamma,b_0}^+}$
on~$\Tb \times \{0\} \cong \Tb$.
Then~$D'$
is automatically continuous on the whole
of~$\Tb \!\times\! [0,1]$,
and a homotopy connecting
$\smash{\alpha_{\Gamma,b_0}^+}$
with~$\smash{\alpha_{\Gamma',b_0}}$.

It~is obvious that a completely analogous construction yields a homotopy connecting
$\smash{\alpha_{\Gamma,b_0}^-}$
with~$\smash{\alpha_{\Gamma',b_0}}$.
The proof is hence~complete.%
\eop
\end{proo_indep}

\setlength\parindent{0mm}
\end{document}